%% file: paper.tex
\documentclass[pdftex,12pt,a4paper]{amsart}
\usepackage[pdftex]{graphicx}
\usepackage{tikz, xcolor}
\DeclareGraphicsExtensions{.png,.pdf}
\usepackage{caption,subcaption}
\graphicspath{{figs/}}
\usepackage[english]{babel}
\usepackage[utf8]{inputenc}
\usepackage{amssymb}
\usepackage{amsmath}
\usepackage{amsthm}
\usepackage{color}
\usepackage{url}
\usepackage{bm}
\usepackage[hidelinks]{hyperref}
\usepackage{cleveref}

\input{macros.tex}

\theoremstyle{definition}
\newtheorem{remark}{Remark}

\newcommand{\enlletres}{%
   \renewcommand{\theenumi}{\alph{enumi}}%
   \renewcommand{\labelenumi}{(\theenumi)}%
}

\title[Structure of the sets of connections]{
   On the structure of the sets of connections between libration
   points and some applications
}

\author[M. Barcelona ]
{Miquel Barcelona }
\address{Departament de Matem\`atiques,
 Universitat Aut\`onoma de Barcelona. Av. de l'Eix Central, Edifici C, 
 08193 Bellaterra (Barcelona), Spain.}
\email{miquel.barcelona@uab.cat}
\author[A. Haro]{Àlex Haro}
\address{Departament de Matem\`atiques i Inform\`atica, Universitat de Barcelona,
Gran Via 585, 08007 Barcelona, Spain \& Centre de Recerca Matemàtica, Edifici C,
Campus Bellaterra, 08193 Bellaterra (Barcelona), Spain.}
\email{alex.haro@ub.edu}
\author[J.M. Mondelo]{J.M. Mondelo}
\address{Departament de Matemàtiques, Universitat Autònoma de
Barcelona, Av.~de l'Eix Central, Edifici C, 08193 Bellaterra (Barcelona),
Spain.}
\email[Corresponding author]{josemaria.mondelo@uab.cat}

\makeatletter
\@namedef{subjclassname@2020}{\textup{2020} Mathematics Subject Classification}
\makeatother
\thanks{ The authors have been supported by the Spanish grants
   PID2020-118281GB-C31, PID2021-125535NB-I00, and PID2025-174015NB-I00
   (MCIU/AEI/FEDER, EU), by the Catalan grants 2021 SGR 00113, 2021 SGR 01072,
   and 2025 SGR 00322 (AGAUR, Generalitat de Catalunya), and by the Air Force
   Office of Scientific Research under award number FA8655-24-1-7059.  This work
   has also received funding from the Spanish State Research Agency through the
   Severo Ochoa and María de Maeztu Program for Centers and Units of Excellence
in R\&D (CEX2020-001084-M).  }
\date{}

\begin{document}

\begin{abstract}
This paper presents a methodology for the global description of the sets of
heteroclinic and homoclinic connections between the center manifolds of the
collinear libration points of the spatial circular Restricted Three-Body
Problem. The approach builds on high-order parameterizations of the
center-stable and center-unstable invariant manifolds and
uses them to extract global information on how connections are
organized across broad energy ranges. It involves: (a) constructing
\emph{connectivity maps} at the level of invariant tori, which encode which
tori in the neighborhoods of $L_1$ and $L_2$ are connected
through heteroclinic and homoclinic trajectories; (b) refining this description
to the level
of individual connections by computing traces of torus--center manifold
connections
and using them to determine the number and distribution of distinct
connecting trajectories between selected pairs of tori.
Results are given for the Earth--Moon and Sun--Earth problems, together
with
examples illustrating how the resulting connectivity information can be used to
design heteroclinic and homoclinic pseudo-trajectories with different purposes.
\end{abstract}

\maketitle
\begingroup
\def\thempfn{$\ast$\ }
\footnotetext{%
Corresponding author.
Address: Departament de Matem\`atiques,
Av.~de l'Eix Central, Edifici C,
Universitat Aut\`onoma de Barcelona,
08193 Bellaterra,
Barcelona,
Spain.
\emph{Phone:} (+34) 93 581 3104
\emph{FAX:} (+34) 93 581 2790
\emph{E-mail address:} {\tt josemaria.mondelo@uab.cat}
}
\endgroup

\section{Introduction}

Homoclinic and heteroclinic connections play an important role in the
description of dynamical systems from a global point of view. They are also of
special interest in applications such as astrodynamics. The pioneering works
of \cite{esa1,esa2,esa3,esa4} unveiled the importance of invariant objects and
their asymptotic invariant manifolds in trajectory design for libration point
space missions, for which the basic model is the circular, spatial restricted
three-body problem (RTBP; see, e.g.~\cite{szebehely}).

The interest in the computation of homoclinic and heteroclinic connections
between invariant objects in the neighborhoods of $L_1$ and $L_2$ of the RTBP
was boosted by the celebrated work~\cite{2000KoLoMaRo}. In that work, a
connection joining two Lyapunov periodic orbits (p.o.) of the Sun-Jupiter planar
RTBP, one around $L_1$ and the other around $L_2$, was used to explain the
(apparently) erratic motion of comet Oterma. This provided a clear dynamical
mechanism showing how such connections can link the neighborhoods of $L_1$ and
$L_2$ and sparked sustained interest in the systematic computation of such
connections. Following this work, several contributions focused on connections
between p.o., including \cite{2006CaMa, 2009BaMoOll, 2013BaMoOll}, as well as
through rigorous numerics~\cite{2003WiZg,2005WiZg}.  Works addressing
connections between invariant tori include \cite{2000GoMa, 2004Goetal,
2005GoMaMo, 2007aArMa, 2008aDeMaRo, 2012CaDoeHuLeOld}, and, more recently,
\cite{2023McCaHo, 2024OwBa, 2024BaHaMo, 2024HeSche}.  References
\cite{2025KuAndLla, 2024LiMaZhu, 2026LiMaZhuGa} compute heteroclinic connections
between tori of periodic extensions of the RTBP.

From an application point of view, a few missions have flown trajectories
explicitly designed to pass near heteroclinic connections; in particular,
\href{https://science.nasa.gov/mission/genesis/}{Genesis} \cite{1997HoBaWiLo}
and \href{https://science.nasa.gov/mission/themis-artemis/}{Themis-Artemis}
\cite{2011aSweetBroetal} used transfers closely related to connections of
invariant tori.
The continued interest in computing heteroclinic connections is
strongly motivated by the desire to facilitate a systematic understanding of
transport mechanisms and the design of complex mission trajectories.  For this
reason, it is important to advance both methodology and the global description
of families of connections.

A key aspect of the spatial RTBP is that the center manifolds of $L_1$ and $L_2$
contain large families of invariant tori, together with lower-dimensional
invariant objects such as periodic orbits, which are among the most abundant
invariant structures near these equilibria and exhibit a rich geometric variety.
From the point of view of applications, this variety is potentially relevant for
communication and observation constraints.  Making it usable requires the
systematic description of families of connections between tori rather than the
study of isolated examples.  While a number of contributions have illustrated
computational methods on selected connections or have followed families only
partially, the first work (to our knowledge) that computes the \emph{global} set
of heteroclinic connections between the center manifolds of $L_1$ and $L_2$ at a
fixed energy level is~\cite{2007aArMa} (in the spatial Hill's problem), using a
semi-analytical approach based on Lindstedt--Poincar\'e expansions \cite{2005Ma}
combined with a reduction-to-the-center-manifold technique \cite{esa3,1999JoMa}.
Reference~\cite{2024BaHaMo} improved the methodology computationally and
provided the full families of heteroclinic connections on several energy levels,
with an approach also semi-analytical but based on the parameterization method,
see \cite{mamotreto}.  More recently, reference~\cite{2024HeSche} was able to
compute a full family of heteroclinic connections in an energy level not
reachable by semi-analytical expansions around the libration points by a purely
numerical approach based on collocation (GMOS);
see~\cite{2001GoMo,2018BaOlSche}.

Our starting point is \cite{2024BaHaMo}. In this reference, the use of the
parameterization method for the computation of the center-stable and
center-unstable manifolds of the $L_1$, $L_2$ points, combined with a meshing
strategy of isoenergetic slices of the center manifold that does not require the
numerical integration of the reduced equations, enabled the fairly automatic and
efficient computation of global sets of connections up to energy levels not
reachable by Linstedt-Poincaré expansions, including those where quasi-Halo
orbits exist. A disadvantage of not using Lindstedt--Poincaré expansions, though,
is that connections are obtained as intersection points of a section, without
some dynamical information often needed in practice:
given two
invariant tori, how many connections do they admit, how do these connections
vary along families (and within a given family), or what role do they play in
organizing the dynamics on the center manifolds. All this information on how the
connections organize dynamically at the level of invariant tori is of major
importance considering its exploitation for trajectory design.

The goal of the present paper is to recover this dynamical information.
To this end, we consider the connection surfaces, which act as departure and
arrival gates for trajectories connecting the center manifolds. These surfaces
are the fundamental geometric objects of our approach, allowing us to
characterize connections involving individual invariant tori through their
intersections with the corresponding connection surfaces. Building on this
description, we introduce the connectivity maps: global maps that encode which
invariant tori are connected.

While our
main focus is on heteroclinic connections between the vicinities of $L_1$ and
$L_2$, the methodology presented in \cite{2024BaHaMo} and in the current work
can be adapted with minor changes to
compute homoclinic sets and their associated connectivity maps.
We have added such computations in order to obtain a more
complete picture of the transport dynamics in the neighborhoods of both
libration points and provide additional examples of how to use them to design
specific pseudo-trajectories.

The connectivity maps are related to the scattering map, which was first
introduced in~\cite{2000DeLlaSe}. The scattering map has become central in the
study of transport phenomena, including its use in proving theoretically the
existence of Arnold diffusion. This tool has been applied to the
design of low-cost transfer trajectories in ~\cite{2006aCaDeMaRo,2008aDeMaRo}
(for the planar RTBP and Hill's problem) and~\cite{2016aDeGiRo,2026aDeGiRo} (for
the spatial RTBP).
The use of the connectivity maps presented here would enable the extension of
these studies to a broader range of energy levels and, more generally, provides
a framework
that can be used to investigate transport and diffusive mechanisms in other
settings as well.

These maps provide information at the level of invariant tori. In addition to
these maps, we also explore the connectivity at the level of individual
connections, namely, given two invariant tori, determine how many connections
exist between them and effectively locate these connections. This is done by
first computing traces of torus vs.~opposite center manifold connections, which
also allows
us to observe different topological configurations.  In particular, by providing
a systematic description of full energy levels, we extend the viewpoint of
\cite{2024OwBa}, where the number of heteroclinic connections between selected
pairs of tori is determined numerically by a knot-theoretical argument.

The paper is organized as follows. In Section~\ref{sec:prelim}, 
we recall some facts about the spatial RTBP and describe the parameterization of
the center-stable and center-unstable manifolds obtained in \cite{2024BaHaMo},
together with the strategy to compute global connection sets from them.
We also fix notation for the rest of the paper and introduce the connection sets
used throughout. In Section~\ref{sec:maps}, we introduce connectivity maps at
the level of tori, which organize the heteroclinic and homoclinic sets and
describe which invariant tori are connected. We also include examples of
trajectory design using these maps. Finally, in Section~\ref{sec:torus2torus},
we present the methodology and results that provide information on the number of
individual connections between pairs of invariant tori. We conclude
Section~\ref{sec:torus2torus} by
presenting an example that illustrates how the number of connections
between a pair of invariant tori can vary when they evolve along their
respective families. We also include several
appendices addressing computational details and a discussion on how the
connections between planar Lyapunov p.o.~lead to changes in the topology of the
connection sets. 

All the explorations presented in this work have been
carried out on an HP Z4 G5 workstation equipped with a 10-core Intel Xeon
w5-2445 processor,
under Debian GNU/Linux 12 with the Xfce
4.18 desktop environment. The source code was written in C, compiled with GCC
12.2.0, and linked against LAPACK 3.11.0. The implementation also uses OpenMP
$4.5$ in most of the computations. In addition, all the numerical integrations
involved in the subsequent sections were performed using a Runge-Kutta-Fehlberg
method of orders 7 and 8 with relative tolerance of $10^{-14}$.

The connection sets used in this work have been gathered in
\cite{2026aBaHaMoData} which can be used under the software from chapter 2 of
\cite{mamotreto}, available at \url{http://www.maia.ub.es/dsg/param/}.

\section{Preliminaries}
\label{sec:prelim}

The present work builds on the results of \cite{2024BaHaMo}. 
We first review some facts about the RTBP. Then we describe the
parameterizations computed there of the center-stable manifold of $L_1$ and the
center-unstable manifold of $L_2$, together with
some of their properties. After recalling some facts on the dynamics on the
center manifolds of $L_1$ and $L_2$, we end the section by describing the sets
of heteroclinic connections computed in \cite{2024BaHaMo}, together with some
additional new sets and related notation.

\subsection{The RTBP}
The Spatial Circular Restricted Three-Body Problem (RTBP) models the dynamics of
a small body of negligible mass under the gravitational influence of two other
bodies, the primaries, which are assumed to move in circular orbits around
their common center of mass.

Let $m_1$ and $m_2$ denote the masses of the primary bodies, with $m_1>m_2$,
and let the mass parameter be defined as $\mu:=m_2/(m_1+m_2) \in [0, 1/2)$.
By considering a uniformly rotating coordinate system, known as the synodic
system, which rotates with the same period as the primaries and has the
center of mass at the origin, both primaries remain fixed on the horizontal
axis.  It is also possible to define a suitable rescaling in space and time
such that the distance between the primaries is equal to 1 and their orbital
period is $2\pi$; see \cite{szebehely} for more details.
In this way, the primaries $m_1$ and $m_2$ are located at $(\mu,0,0)$ and
$(\mu-1,0,0)$, respectively, and their masses become $1-\mu$ and $\mu$,
respectively. Introducing the momenta as $x_4 = \dot x_1 - x_2$, $x_5 = \dot x_2
+ x_1 $ and $x_6 = \dot x_3$, the dynamics of the third body
is described by the autonomous Hamiltonian system defined by
\[
   H(x) = \frac{1}{2}(x_4^2 + x_5^2 + x_6^2 ) - x_1x_5 + x_2x_4 -
   \frac{1-\mu}{r_1} - \frac{\mu}{r_2},
\]
where $x \in \R^6$, $r_1^2=(x_1-\mu)^2+x_2^2+x_3^2$ and
$r_2^2=(x_1-\mu+1)^2+x_2^2+x_3^2$. From now on, we will refer to the value of
the Hamiltonian as {\em the energy}. 

The corresponding system of ordinary differential equations can be expressed
as
\begin{equation}
   \dot x = X_H(x) := J \nabla H(x),
   \label{eq:vf}
\end{equation}
where $J$ is the $6 \times 6$ skew-symmetric matrix 
\[
   J = 
   \begin{bmatrix}
      0 & I \\ -I & 0
   \end{bmatrix}
\]
and $I$ is the $3 \times 3$ identity matrix.

The RTBP exhibits
some symmetries: if
\[
   x(t)=(x_1(t),x_2(t),x_3(t),x_4(t),x_5(t),x_6(t))
\]
is a solution of the system, then 
\begin{gather}
   \label{eq:rev}
   R:x(t) \rightarrow (x_1(-t),-x_2(-t),x_3(-t),-x_4(-t),x_5(-t),-x_6(-t)), \\
   S:x(t) \rightarrow (x_1(t),x_2(t),-x_3(t),x_4(t),x_5(t),-x_6(t))
   \label{eq:sym}
\end{gather}
are also solutions. These symmetries will be useful later on
in the exploitation of the numerical results.

\begin{figure}[htbp]
   \centering
   \input{rtbp.tex}
   \caption{Schematic representation of the RTBP. The dynamics of the third body
   is not restricted to the $x_1x_2$-plane even though it is not shown for
   simplicity purposes.}
   \label{fig:rtbp}
\end{figure}
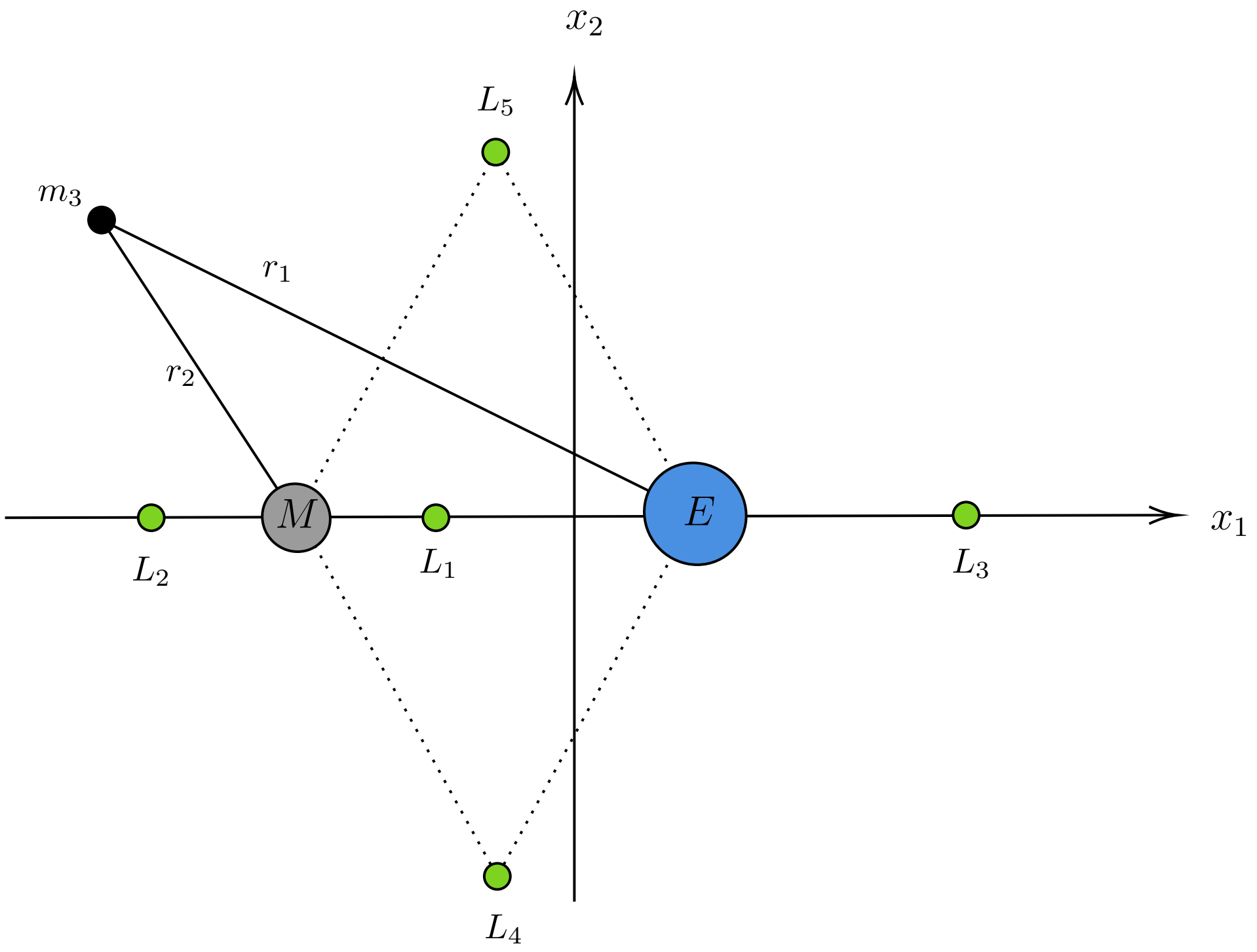

It is well known that the RTBP has five equilibrium points, denoted by $L_j$
for $j=1,\dots,5$.  The first three points, $L_1$, $L_2$ and $L_3$,
discovered by Euler, lie on the $x_1$-axis and are therefore called collinear
points.
Their corresponding $x_1$ coordinates are
given by
$\mu-1+\gamma_1$, $\mu-1-\gamma_2$, and $\mu+\gamma_3$, respectively,
where $\gamma_j$ is the only positive root of the corresponding Euler's
quintic equation,
\begin{equation}\label{eq:quintic}
\begin{aligned}
   \gamma_j^5 \mp (3-\mu)\gamma_j^4 + (3-2\mu)\gamma_j^3 - \mu \gamma_j ^2
   \pm 2\mu \gamma_j - \mu &= 0,  \quad &j&=1,2,\\
   \gamma_j^5 + (2+\mu)\gamma_j^4+(1+2\mu)\gamma_j^3 - (1-\mu)\gamma_j^2 -
   2(1-\mu)\gamma_j - (1-\mu) &= 0, \quad &j&= 3,  
\end{aligned}
\end{equation}
which
provides the distance to the closest primary for each of the collinear
points. Meanwhile, $L_4$ and $L_5$, discovered by Lagrange, are known as the
triangular points because they form equilateral triangles with the primaries.
See Figure~\ref{fig:rtbp}.

The linear behavior of the triangular points depends on the value of $\mu$. In
the particular case where the primaries are the Earth and the Moon, which is
one of our cases of interest, the mass parameter is given by\footnote{As
obtained from DE406 JPL ephemeris file.}
\[
   \mu=1.215058560962404 \times 10^{-2},
\]
whereas, when the primaries are the Sun and the Earth, we have\footnote{In
this case, the DE441 JPL ephemeris file has been used.}
\[
   \mu = 3.003480594542193\times 10^{-6}.
\]
In both cases these two points are linearly stable (see \cite{szebehely}).  On
the other hand, the collinear points are all unstable and of type center
$\times$ center $\times$ saddle, for any value of $\mu$.  Therefore, by the
center manifold theorem (see e.g.~\cite{guckenheimer-holmes}), the collinear
points have a four-dimensional center manifold and one-dimensional stable and
unstable manifolds. 

From now on, we will focus our attention on $L_1$ and $L_2$. Let us denote
their associated set of eigenvalues by
\[
   \Spec
   DX_H(L_{j}) = \{\pm\bmi\omega_p^{(j)}, 
   \pm\bmi\omega_v^{(j)},
   \pm\lambda^{(j)}\}, \quad
   j=1,2,
\]
where $\bmi$ denotes the imaginary unit. Corresponding eigenvectors are denoted
by
\[
v_p^{(j)} \in E\bigl(\bmi\omega_p^{(j)}\bigr), \quad v_v^{(j)} \in
E\bigl(\bmi\omega_v^{(j)}\bigr), \quad v_+^{(j)} \in E\bigl(\lambda^{(j)}\bigr), \quad
v_-^{(j)} \in E\bigl(-\lambda^{(j)}\bigr),
\]
where $E(\lambda)$ denotes the space of eigenvectors corresponding to the
eigenvalue $\lambda$.

Lyapunov's center
theorem \cite{siegel-moser,meyer-hall} applies to each pair of
complex-conjugate eigenvalues, ensuring the existence of two one-parameter
families of periodic orbits (p.o.), each one spanning a two-dimensional
manifold tangent to $\langle v_c^{(j)},\bar v_c^{(j)}\rangle$, for
$c\in\{p,v\}$, at the collinear points. These families, called planar and
vertical due to their tangent directions, extend across the center manifold
and, for the dynamical regime considered in this paper, can be parametrized by
the energy. The bifurcations of these orbits give rise to other families of
periodic orbits such as the
Northern and Southern Halo families, symmetric through symmetry
\eqref{eq:sym}, that are born simultaneously at the first $1:1$ bifurcation of
the planar Lyapunov family. Their denomination comes from the fact that, in the
$x_2$-$x_3$ projection, they look like a halo that, in the northern
(resp.~southern) case, has a larger length contained in the $x_3>0$
(resp.~$x_3<0$) half-plane than in the $x_3<0$ (resp.~$x_3>0$) one.

Surrounding these periodic orbits, families of invariant tori can be found.
The Lissajous family of invariant tori is the one whose birth from the
equilibrium point is predicted by KAM theory \cite{1997JoVi-normal,eliasson}.
They can be thought of as arising from the interaction of the two linear
oscillations around the equilibrium point, and they are named after the fact
that their tori are filled with Lissajous-type trajectories. For energies prior to
the bifurcation of the Halo family of p.o., the Lissajous family of tori can
also be thought of as connecting the planar and vertical Lyapunov families of p.o.
Another family of tori is the quasi-Halo one, that is formed by tori that
surround the Halo family of p.o., originated by the (partial) stability
transfer from the planar Lyapunov family of p.o.~to the Halo one after its 1:1
bifurcation. For the energies considered in this paper, the families of
p.o.~and tori mentioned can be considered to foliate the center manifold of the
collinear points. The cantorian gaps predicted by KAM theory are almost
undetectable when working numerically in double precision.

All the families of periodic orbits and tori mentioned inherit the hyperbolic
direction of the equilibrium point, and thus have stable and
unstable manifolds outside the center manifold. Their intersection gives rise
to homoclinic and heteroclinic connections, that are the subject of this
work.

\subsection{Parameterizations of center-stable and center-unstable manifolds}
\label{ssec:param}

The parameterization method for invariant manifolds of fixed points
used in \cite{2024BaHaMo}
produces Taylor expansions of parameterizations of
$\calW^{(1)}:=\calW^{cs}(L_1)$, $\calW^{(2)}:=\calW^{cu}(L_2)$, together with
Taylor expansions of the differential equations of the RTBP reduced to each of
these manifolds, expressed in the coordinates of the parameterizations.

For $j=1,2$, denoting by $W^{(j)}:\calU^{(j)}\subset\R^5\rightarrow\R^6$ the
parameterization obtained for $\calW^{(j)}$, its (truncated) Taylor expansion
is of the form
\begin{equation}\label{eq:param}
   W^{(j)}(s)
   =
   L_j+
   P^{(j)}s
   +
   \sum_{\stackrel{\scriptstyle m\in\N^5}{2\leq|m|\leq K}}
   W^{(j)}_ms^m,
\end{equation}
where $L_j\in\R^6$ is the libration point and, in the non-linear terms,
$W^{(j)}_m\in\R^6$, $m=(m_1,\dots,m_5)$, $|m|=m_1+\dots+m_5$,
$s=(s_1,\dots,s_5)^\top$, $s^m=s_1^{m_1}\dots s_5^{m_5}$ and $K$ the order of the
parameterization. In the linear term, $s$
is as before and $P^{(j)}\in\R^{6\times5}$
has the columns
\begin{equation}\label{eq:linparam}
   P^{(j)}
   =
   \begin{pmatrix}
      \\
      2\Re(v_p^{(j)})
      &
      -2\Im(v_p^{(j)})
      &
      2\Re(v_v^{(j)})
      &
      -2\Im(v_v^{(j)})
      &
      v_{\mp}^{(j)}
      \\
      \\
   \end{pmatrix}
   ,
\end{equation}
where $v_{\mp}$ is $v_-$ for $j=1$, $v_+$ for $j=2$.
The eigenvectors $v_p^{(j)},v_v^{(j)},v_\pm$ are taken normalized as
to have norm equal to the solution of Euler's quintic equation
\eqref{eq:quintic}, that is, $\|v_p^{(j)}\|=\|v_v^{(j)}\|=\|v_\pm\|=\gamma_j$.

Denote by $\dot s=f^{(j)}(s)$, with
$f^{(j)}:\calU^{(j)}\subset\R^5\rightarrow\R^5$, the RTBP equations reduced to
$\calW^{(j)}$ in parameter coordinates $s\in\calU^{(j)}$. This means that, if
the corresponding flow is denoted as $\phi^{(j)}_t(s)$, the following
invariance equation is satisfied:
\begin{equation}\label{eq:invfl}
   \Phi_t\bigl(W^{(j)}(s)\bigr)
   =
   W^{(j)}\bigl(\phi^{(j)}_t(s)\bigr)
   ,
\end{equation}
where $\Phi_t$ is the flow associated to \eqref{eq:vf}.
The parameterization method as used in \cite{2024BaHaMo} also produces
(truncated) Taylor expansions for $f^{(j)}$, of the form\footnote{Actually,
the Taylor expansions of $W^{(j)}$ and $f^{(j)}$ computed in \cite{2024BaHaMo}
have complex coefficients, because the parameterization method is much more
easily formulated (and implemented) in this way. When evaluating them on real
parameters, symmetries were exploited in order to prevent computational
overhead. For clarity, we have preferred to write the equivalent real
expressions here.}
\begin{equation}\label{eq:intdyn}
   f^{(j)}(s)=
      \begin{pmatrix}
         0 & -\omega_p & 0 & 0 & 0 \\
         \omega_p & 0 & 0 & 0 & 0 \\
         0 & 0 & 0 & -\omega_v & 0\\
         0 & 0 & \omega_v & 0 & 0 \\
         0 & 0 & 0 & 0 & \lambda_\mp
      \end{pmatrix}
      s
      +
      \sum_{\stackrel{\scriptstyle m\in\N^5}{2\leq|m|\leq K}}
      f^{(j)}_ms^m,
\end{equation}
with $f^{(j)}_m\in\R^5$ and $s,m$ as before. Equation \eqref{eq:invfl}
establishes the invariance of the manifolds parameterized by $W^{(j)}$.
For $j=1$ (resp.~$j=2$), 
it follows 
from the expression \eqref{eq:linparam} of the linear
part of \eqref{eq:param} that the manifold parameterized by $W^{(j)}$ is
tangent to the center-stable (resp.~center-unstable) eigenspace of $L_j$,
whereas from the linear part of \eqref{eq:intdyn} follows that the internal
dynamics of $\calW^{(j)}$ corresponds to that of a center-stable
(resp.~center-unstable) manifold.

The parameterization method is based on solving the infinitesimal version of
the invariance equation \eqref{eq:invfl} order by order. When doing this,
there is some freedom in the choice of the terms of the Taylor expansions
\eqref{eq:param}, \eqref{eq:intdyn} that leads to several parameterization
styles. A discussion can be found in \cite{mamotreto}. Reference
\cite{2024BaHaMo} introduced a {\em mixed-uncoupling} style that led to
parameterizations $W^{(j)}$ and reduced vector fields $f^{(j)}$ satisfying the
following properties:
\begin{enumerate}\enlletres
   \item\label{en:propestil1}
   Points of the form $(s_1,s_2,0,0,0)$ describe the planar family of Lyapunov
   p.o.~around $L_j$. That is, for any $(s_1,s_2,0,0,0)\in\calU^{(j)}$, the
   set of points
   \[
      \{W^{(j)}\bigl(\phi^{(j)}_t(s_1,s_2,0,0,0)\bigr)\}_{t\in\R}
      =
      \{\Phi_t\bigl(W^{(j)}(s_1,s_2,0,0,0)\bigr)\}_{t\in\R}
   \]
   is a planar Lyapunov p.o.~around $L_j$.
\item\label{en:propestil2}
   Points of the form $(0,0,s_3,s_4,0)$ describe the vertical family of
   Lyapunov p.o.~around $L_j$.
\item\label{en:propestil3}
   Points of the form $(s_1,s_2,s_3,s_4,0)$ describe the center manifold of
   $L_j$.
\item\label{en:propestil4}
   For any $s=(s_1,s_2,s_3,s_4,s_5)$, the point $s$ is in the fiber of
   $(s_1,s_2,s_3,s_4,0)$. This means that
   \[
      \phi_t(s_1,s_2,s_3,s_4,s_5)
      -
      \phi_t(s_1,s_2,s_3,s_4,0)
      \stackrel{t\rightarrow\pm\infty}{\makebox[3em]{\rightarrowfill}}
      0,
   \]
   where $\pm\infty$ is $+\infty$ for $j=1$, $-\infty$ for $j=2$. From the
   invariance equation \eqref{eq:invfl}, we have the same property for the
   corresponding manifold points:
   \[
      W^{(j)}\bigl(\phi_t(s_1,s_2,s_3,s_4,s_5)\bigr)
      -
      W^{(j)}\bigl(\phi_t(s_1,s_2,s_3,s_4,0)\bigr)
      \stackrel{t\rightarrow\pm\infty}{\makebox[3em]{\rightarrowfill}}
      0,
   \]
\end{enumerate}
It follows from property \eqref{en:propestil4} above that, if a point with
parameters $(s_1,\linebreak[1]s_2,\linebreak[1]s_3,\linebreak[1]s_4,0)$ is on
an object (say a periodic orbit or an invariant torus) of $\calW^{c}(L_1)$
(resp.~$\calW^{c}(L_2)$), a point with parameters $(s_1,s_2,s_3,s_4,s_5)$ is
in the stable (resp.~unstable) fiber of this object. In
particular, they have the same energy.

The computed truncated Taylor expansions of the parameterizations of the
center-stable and center-unstable manifolds are divergent; however they provide
good local approximations within a (not small) neighborhood of the libration
points.
From now on, and following \cite{2024BaHaMo}, we choose the order of the
parameterizations to be $K=30$ and limit the coordinate on the fiber so as to
satisfy $|s_5| \leq 10^{-3}$. With these choices, we ensure that the
parameterizations are sufficiently accurate for the energy values of the
Sun-Earth and Earth-Moon cases that will be considered in what follows.

\subsection{Some aspects of the dynamics of $\calW^c(L_j)$}
\label{ssec:dynwc}

As seen in the previous section, the parameterization $W^{(1)}$
(resp.~$W^{(2)}$) of the center-stable (resp.~center-unstable) manifold of $L_1$
(resp.~$L_2$)
allows us to recover the center manifold of $L_1$ (resp.~$L_2$) by putting
$s_5=0$. Then the dynamics of the center manifold can be recovered by direct
numerical integration of the reduced flow, $\phi_t^{(j)}$. In this section, we
will recall some known properties of these dynamics in order to refer to it in
later sections. Equivalent explorations were first done in
\cite{esa3,1999JoMa}.

Let us define the restriction of the parameterization $W^{(j)}$ to the center
manifold as
\[
   W_c^{(j)} (c) := W^{(j)} (s_1,s_2,s_3,s_4,0), 
\]
where $c=(s_1,s_2,s_3,s_4) \in \Vcalj \subset \R^4$,
and
$\Vcalj$ is an open subset on which $W_c^{(j)}$ provides an accurate
parameterization of the center manifold.
By the uncoupling
property \eqref{en:propestil4},
the first four components of $f^{(j)}$ do not depend on $s_5$. Therefore, we 
denote the corresponding reduced flow restricted to the center manifold by
$\varphi_t^{(j)} : D \subset \R \times \Vcalj \rightarrow \Vcalj$, which is
defined as
\[
   \varphi_t^{(j)} (c) :=  \bigl(\phi_t^{(j)}(s_1,s_2,s_3,s_4,0)\bigr)_{i=1}^4.
\]

Note that, in principle, thanks to the invariance
equation \eqref{eq:invfl}, trajectories in the center manifold of $L_j$ can
be generated numerically through the two following (mathematically) equivalent
expressions:
\[
   \Phi_t\bigl(W_c^{(j)}(c)\bigr)
   =
   W_c^{(j)}\bigl(\varphi_t^{(j)}(c)\bigr).
\]
The left-hand expression involves numerically integrating the RTBP equations
\eqref{eq:vf}. This is much faster than using the right-hand expression, which
involves numerical integration of the first 4 components of $\dot s=f^{(j)}(s)$,
with $f^{(j)}$ given by
a 4-variate Taylor series of order 30 with 4-dimensional vector
coefficients\footnote{The series is actually 5-variate as in
\eqref{eq:intdyn}, but, when $s_5=0$, our implementation skips both the terms
with $m_5\neq0$ and the evaluation of the last component of $f^{(j)}$. Even
so, the total number of terms is 1\,623\,160.}. On the
other hand, 
direct integration of the RTBP equations becomes numerically inaccurate for large
$t$
because rounding errors would be propagated by 
the (strong) hyperbolicity of
$L_j$. 
The right-hand
expression, however, can be evaluated for arbitrarily large $t$, with the
computational cost already mentioned. Since this cost is affordable, the
right-hand expression allows us to investigate the dynamics of the center
manifold by direct numerical integration, whereas the left-hand expression
does not.

A classical approach for the study of the dynamics in the center manifold,
introduced in \cite{esa3,1999JoMa}, is to consider isoenergetic Poincaré
sections. Since the center manifold is $4D$, isoenergetic Poincaré sections
are planar. The dynamics is then visualized as a sequence of 2D plots, each
corresponding to a particular energy value.  Here, we will discuss two
specific isoenergetic Poincaré sections relevant to this work, but that
qualitatively represent all the possible scenarios in the range of energies
in which the Taylor expansions of $W^{(j)}$, $f^{(j)}$ are valid.

\begin{figure}[htb]
   \centering
   \includegraphics[width=.49\textwidth]{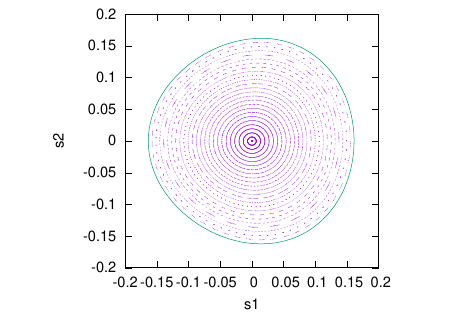}
   \includegraphics[width=.49\textwidth]{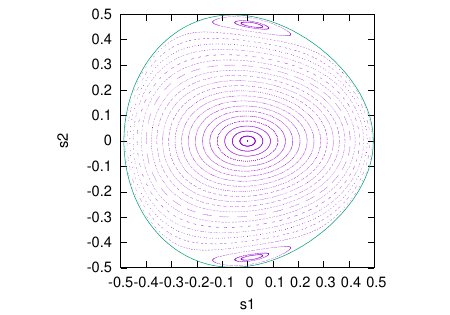}
   \caption{\label{fig:poinc} Isoenergetic $\poincpi$ Poincaré sections
   of $\calW^c_{-1.585}(L_2)$ and $\calW^c_{-1.585}(L_1)$, in manifold coordinates.}
\end{figure}

An isoenergetic slice of $\calW^c(L_j)$ of energy $h$ is effectively defined
as \[ \calW_h^c(L_j):=\calW^c(L_j)\cap\{x\in\R^6:H(x)=h\}, \] and is 3D. A
Poincaré section of it will reduce the restricted flow to a planar map, whose
phase portrait is easily visualized. Such phase portraits are represented in
Figure~\ref{fig:poinc}\footnote{These Poincaré maps are planar and
area-preserving, and therefore show typical features of such maps (see
e.g.~\cite{arrowsmith-place}).}. The left plot corresponds to a Poincaré
section of $\calW^c_{-1.585}(L_2)$, whereas the right one corresponds to
$\calW^c_{-1.585}(L_1)$. In both cases, the Poincaré section is $\poincpi$,
that corresponds to considering only cuts that go from $s_4<0$ to $s_4>0$. The
plots are actually generated from 25 points of the form
$s=(s_1,s_2,s_3,s_4,0)$, with $s_1=0$, $s_2$ taking equally spaced positive
values between $0$ and the green curve, $s_4=0$ and $s_3$ obtained by solving
$H\bigl(W^{(j)}(s)\bigr)=-1.585$. Then, $\dot s=f^{(j)}(s)$ is numerically
integrated from these 25 initial conditions, but plotting only the cuts with
$\poincpi$. For each initial condition, a total of 200 cuts are
effectively computed and displayed in the plots of Figure~\ref{fig:poinc}. The
outer green curve is not made of Poincaré cuts.  Instead, it is a (continuous)
trajectory of the reduced flow, made of points with $s_3=s_4=0$.  It is
therefore, by property (\ref{en:propestil1}) of the mixed-uncoupled
parameterization style used (see the previous Section), the planar Lyapunov
orbit of the energy level.  So, technically, this green curve is not part of
the Poincaré map, but it is convenient to display it since it is the boundary
of the region where the Poincaré map is defined.

\begin{figure}[htb]
\includegraphics{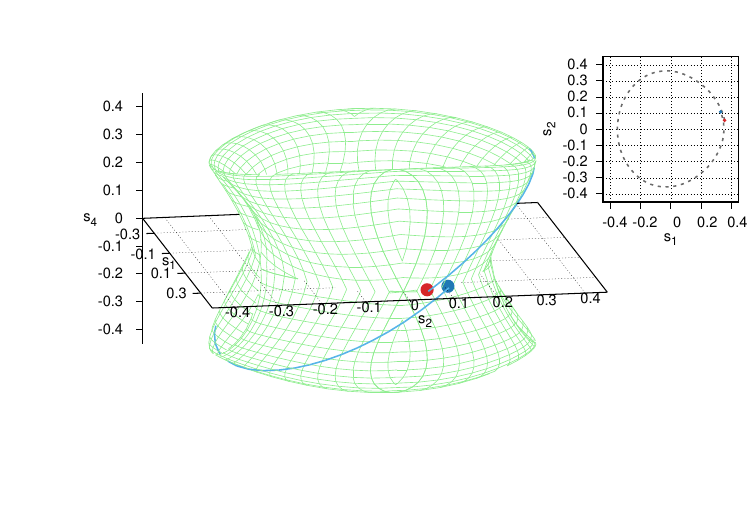}
\caption{
   Representation of the dynamics on the center manifold of $L_1$ projected
   onto the $s_1s_2s_4$ coordinates.  The motion on the invariant torus (green)
   evolves from the red point to the blue one in the Poincaré section
   $\poincpi$, advancing with an average angular rate $\rho$.
   \label{fig:trajtor}
}
\end{figure}

In both plots of Figure \ref{fig:poinc}, the fixed points with $s_1=s_2=0$
correspond, according to property (\ref{en:propestil2}), to the vertical
Lyapunov p.o.~of the energy level. The surrounding invariant curves correspond
to traces (through the Poincaré section) of invariant 2D tori of the Lissajous
family. Along each of these invariant curves, consecutive Poincaré cuts evolve
counterclockwise, at an average angular rate that is known as {\em rotation
number}, that we will denote as $\rho$. The rotation number varies from an
invariant curve to a neighboring one. Namely, in Figure~\ref{fig:poinc} left,
$\rho$ ranges from $44.7 \times 10^{-3}$ at the vertical Lyapunov
p.o.~($s_1=s_2=0$) to $45.7 \times 10^{-3}$ at the planar Lyapunov p.o.~(green
curve).  Outside the Poincaré section, from each Poincaré cut to the next,
trajectories also go counterclockwise in the $(s_1,s_2)$ projection by nearly
completing a revolution. See Figure~\ref{fig:trajtor}.

In Figure~\ref{fig:poinc} right,
in addition to the fixed point with $s_1=s_2=0$, two additional fixed points
can be found. The one with $s_2>0$
corresponds to the Southern
Halo p.o.~of the energy level, whereas the one with $s_2<0$
corresponds to the Northern Halo p.o.\footnote{In $L_2$ $\poincpi$ Poincaré
sections, the fixed point with $s_2>0$ corresponds to the Northern Halo orbit,
and the one with $s_2<0$ to the Southern one.} The invariant curves that
surround these additional fixed points are traces of 2D quasi-Halo tori
surrounding the Halo p.o. In this case, consecutive Poincaré cuts evolve
clockwise.

\begin{figure}[htb]
\centering
\includegraphics{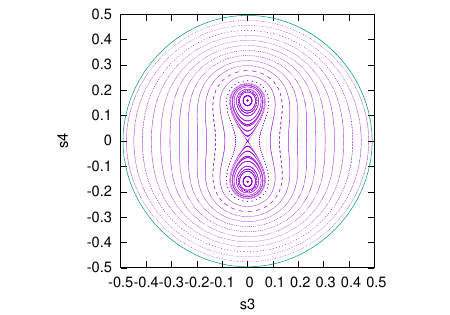}
\caption{\label{fig:poinc2}Isoenergetic Poincaré section of
$\calW^c_{-1.585}(L_1)$ with $\{s_2=0\}$ in manifold coordinates.}
\end{figure}

Two homoclinic connections of the planar p.o.~that are inside the center
manifold act as separatrices between quasi-periodic motion of Lissajous and
quasi-Halo type. These connections are better visualized in a different
Poincaré section.  Figure~\ref{fig:poinc2} displays the phase portrait of the
$\{s_2=0, s_1>0\}$ Poincaré map of $\calW^c_{-1.585}(L_1)$, that corresponds
to consider only cuts from $s_2<0$ to $s_2>0$. Here, the point with
$s_3=s_4=0$ corresponds to the planar Lyapunov p.o.~of the energy level, and
the two homoclinic connections previously mentioned are outlined by the
8-shaped curve that goes through $s_3=s_4=0$. It is worth mentioning that,
with respect to the flow $\varphi_t^{(1)}$, these homoclinic connections are
(continuous) trajectories connecting the planar Lyapunov p.o.~with itself.
From the point of view of the Poincaré map, however, they are discrete
trajectories connecting the fixed point $s_3=s_4=0$ of the Poincaré map with
itself. The invariant manifolds of this fixed point of the Poincaré map that
give rise to these connections are continuous curves that give rise to a
homoclinic tangle that is ``too smashed'' to be easily
visualized\footnote{This is actually related to exponentially small splitting
phenomena. For a discussion on homoclinic tangles of
area-preserving maps see e.g.~\cite{guckenheimer-holmes,arrowsmith-place}},
although a zoom around $s_3=s_4=0$ in the plot of Figure~\ref{fig:poinc2}
starts to slightly outline it. The traces of the Halo p.o.~and the quasi-Halo
tori are the invariant curves inside the 8-shaped curve, with northern
families in $s_4<0$ and southern ones in $s_4>0$. The invariant curves
surrounding the 8-shaped curve are traces of Lissajous tori. The outer green
curve is a continuous trajectory of the flow $\varphi_t^{(1)}$ made of points with
$s_1=s_2=0$, and is therefore the vertical Lyapunov orbit of the energy level.
As before, it is not part of the Poincaré map but the boundary of the region
where the Poincaré map is defined.

\subsection{Computing sets of heteroclinic connections}
\label{sec:sets}

Reference \cite{2024BaHaMo} introduced a numerical methodology for computing
discretizations of global sets of
heteroclinic connections between isoenergetic slices of $\calW^c(L_2)$ and
$\calW^c(L_1)$. More precisely, for several values of $h$, (finite) sets of
heteroclinic connections were computed that completely span the set
\begin{equation}\label{eq:varhomisoen}
   \Gamma_h^{2,1}:=\calW^{cu}(L_2)\cap\calW^{cs}(L_1)\cap\{H=h\}.
\end{equation}
The remainder of this section summarizes the methodology introduced there, which
will be used throughout the present work. This methodology is also formulated so
as to cover the homoclinic case. In particular, we extend the results of
\cite{2024BaHaMo} with additional connecting sets and we introduce related
notation necessary for the exploitation of these results in the following
sections.

Denote, thanks to property (\ref{en:propestil3}) of the mixed-uncoupling
parameterization style (see the previous Section), the parametric
representation of an isoenergetic slice of the center manifold as
\[
   S^{(j)}_h=\{c\in\R^4:H\bigl(W_c^{(j)}(c)\bigr)=h\}.
\]
Reference \cite{2024BaHaMo} computes meshes $\hat S_h^{(j)}$ of $S^{(j)}_h$.
In order to find connections (either heteroclinic or homoclinic) as points
rather than trajectories, we introduce a Poincaré section located midway
between $L_1$ and $L_2$, $\Sigma:=\{x\in\R^6 :x_1=\mu-1\}$.
For $l,m\in\{1,2\}$, we seek to determine the set
\[
   \Gamma_{h,k}^{l,m}
   :=
   \calW^{cu}(L_l)\cap\calW^{cs}(L_m)\cap\{H=h\}\cap{\Sigma}
   ,
\]
where $k:=k_u+k_s$, and $k_s$ (resp.~$k_u$) denotes the number of intersections
of the stable (resp.~unstable) manifold with $\Sigma$. We will denote the
corresponding connections as $k$-heteroclinic (or homoclinic). Observe that the
total number of crossings of the Poincaré section along a connection  is $k-1$,
whereas the number of loops around the second primary is $(k-2)/2$ in the
heteroclinic case and $(k-1)/2$ in the homoclinic one.  When $k$ corresponds to
the minimal value for which $\calW^{cu}(L_l)$ and $\calW^{cs}(L_m)$ intersect,
we may omit it from the notation for simplicity and refer to the corresponding
connections as {\em primary connections}.  Further details on primarity are
given in Appendix~\ref{app:primarity}.

So as to compute $\Gamma_{h,k}^{l,m}$, the manifolds $\calW^{cu}(L_l)$ and
$\calW^{cs}(L_m)$ are numerically globalized through the flow $\Phi_t$,
starting from the points 
\begin{equation}
  \{W^{(j)}(c,\xi^{(j)}),c\in \hat S_h^{(j)}\}, 
  \label{eq:seclocws}
\end{equation}
for\footnote{The $v_\pm$ eigenvectors can be chosen so that the manifold
branches we want to follow correspond to $\xi^{(j)}$ positive.}
$\xi^{(j)}=10^{-3}$. 
Let us denote as $\hat\Gamma^{(j)}_{h,k}$ the
$\Sigma$-sections of the trajectories starting from the points in
\eqref{eq:seclocws}. An all-against-all search for minimum
distances\footnote{The $k$-$d$ tree algorithm is able to do this in $O(N\log
N)$ operations, where $N$ is the number of points in $\hat\Gamma_{h,k}^{(j)}$,
$j=1,2$. On a search between two sets of 6-dimensional $4 \cdot 10^{6}$
points, our parallelized implementation, that uses the
\href{https://www.alglib.net/}{ALGLIB} library, achieves a total runtime of
2.126 seconds on an AMD Ryzen 7 3700X CPU with 16 threads.  The brute-force
strategy takes 202.38 seconds on an NVIDIA RTX A4000 GPU with a total of 6144
FP32 cores.} between points in $\hat\Gamma_{h,k}^{(1)},\hat\Gamma_{h,k}^{(2)}$
provides initial approximations of connections that can be refined into true
ones through Newton iterations. Careful remeshings lead to fine grids of the
global set $\Gamma^{l,m}_{h,k}$. 
The above procedure produces fine discretizations of the sets
$\Gamma_{h,k}^{l,m}$. Further details can be found in
\cite{2024BaHaMo}.

In the present work, these discretizations are taken as the starting point for
the new developments. Specifically, we represent them by finite sets
$\hat R_{h,k}^{l,m}$ of tuples
$(c^{(l)},T^{(l)},c^{(m)},T^{(m)})$
satisfying 
$c^{(j)}\in S_h^{(j)}$, $T^{(j)}\in\R$, $j=1,2$, and
\begin{equation}
   \Phi_{T^{(l)}}\bigl(W^{(l)}(c^{(l)},\xi^{(l)})\bigr)
   =
   \Phi_{T^{(m)}}\bigl(W^{(m)}(c^{(m)},\xi^{(m)})\bigr)
   \in\Sigma
   ,
   \label{eq:conncond}
\end{equation}
with $\xi^{(l)}=\xi^{(m)}=10^{-3}$. The set $\hat
R_{h,k}^{l,m}$ can be regarded as (and actually is) a numerical dataset --- or
equivalently, a ``results file'' --- containing all the computed connection
tuples satisfying \eqref{eq:conncond}.

The corresponding reverse set $\hat R_{h,k}^{m,l}$ 
can be directly obtained by applying the
reversibility property in \eqref{eq:rev} to the original set $\hat R_{h,k}^{l,m}$.

\begin{figure}[tbp]
   \centering
   \includegraphics{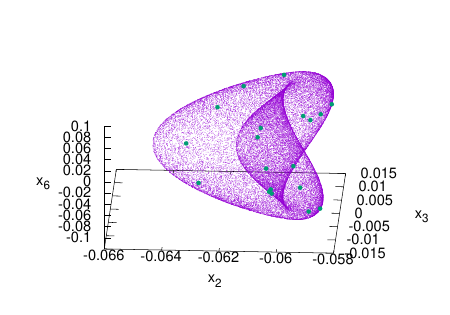}
   \includegraphics{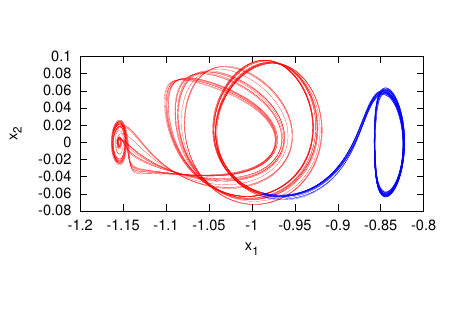}
   \caption{Heteroclinic connections in the Earth-Moon RTBP for $h=-1.585$. Left:
   numerical discretization of the set
   $\Gamma_{-1.585,6}^{2,1}$ in phase space. Right: heteroclinic trajectories
   corresponding to the highlighted
   green points in $\Gamma_{-1.585,6}^{2,1}$, connecting $L_2$ to $L_1$, as defined by
   the set $R_{-1.585,6}^{2,1}$. 
   }
   \label{fig:Gammah}
\end{figure}

Reference \cite{2007aArMa} was the first to visualize an isoenergetic slice of a
heteroclinic manifold connecting the neighborhood of $L_1$ and $L_2$ as a
topological sphere, in the context of Hill's problem. For the Earth-Moon RTBP,
the left panel of Figure~\ref{fig:Gammah} 
shows a discretization of $\Gamma^{2,1}_{h,k}$ for 
$h=-1.585$ and $k=6$, 
obtained from the dataset $\hat R_{h,k}^{2,1}$ computed in
\cite{2024BaHaMo}. This set also has
the topology of a sphere.  Some of the corresponding
heteroclinic connections are shown in the right panel of Fig.~\ref{fig:Gammah}.

Rather than working directly with $\Gamma^{l,m}_{h,k}$, we consider the sets
\[
      \hat{\mathcal{D}}_{h,k}^{l,m} := \Pi_1 \hRlmhk \subset \calV^{(l)},
      \qquad
      \hat{\mathcal{A}}_{h,k}^{l,m} := \Pi_3 \hRlmhk \subset \calV^{(m)}, 
\]
where $\Pi_i$ denotes the projection onto the $i$-th element of the tuples.
These sets describe the departure and arrival points of the 
connections in the corresponding isoenergetic slice $S_h^{(j)}$ of the center
manifolds.
We refer to them as the {\em departure and arrival connection surfaces} (CS),
respectively. 

Although all the connection surfaces computed in \cite{2024BaHaMo} were found to
be topological spheres, the methodology is not restricted to this topology.
Indeed, the presence of tangencies, changes in intersection multiplicity with
$\Sigma$, or the coexistence of multiple connection mechanisms may lead to more
general topologies. These phenomena are investigated in
Appendix~\ref{app:primarity}.

\begin{remark}
Unlike in the heteroclinic case, the departure and arrival connection
surfaces $\hat{\mathcal{D}}_{h,k}^{l,m}$ and
$\hat{\mathcal{A}}_{h,k}^{l,m}$ coincide in the homoclinic case,
due to the symmetries of the RTBP and the structure of the
parameterization.
\end{remark}

The sets $\hRlmhk$, $\hat{\mathcal{D}}_{h,k}^{l,m}$,
$\hat{\mathcal{A}}_{h,k}^{l,m}$ are finite by construction. The $\hRlmhk$ set
defines a (discrete) one-to-one map between $\hat{\mathcal{D}}_{h,k}^{l,m}$
and $\hat{\mathcal{A}}_{h,k}^{l,m}$ (or vice-versa) that sends $c^{(l)}\mapsto
c^{(m)}$ according to Eq.~\eqref{eq:conncond}. As can be seen in Figure 11 of
\cite{2024BaHaMo}, it turns out that, in all the cases considered, for a
suitable choice of indexes $\{i_1,i_2,i_3,i_4\}\subset\{1,2,3,4\}$, the sets
$\hat{\mathcal{D}}_{h,k}^{l,m}$ and $\hat{\mathcal{A}}_{h,k}^{l,m}$, when
projected onto the coordinates $c_{i_1}, c_{i_2}, c_{i_3}$, are convex and, in
this particular case, sufficiently close to true spheres as to allow
parameterization using spherical coordinates from the center of their bounding
boxes. This fact
allows us to obtain accurate approximations of (continuous) 2D surfaces
$\mathcal{D}_{h,k}^{l,m}$ (resp.~$\mathcal{A}_{h,k}^{l,m}$), meshed by the
$\hat{\mathcal{D}}_{h,k}^{l,m}$ sets (resp.~$\hat{\mathcal{A}}_{h,k}^{l,m}$),
and one-to-one maps between them that are evaluated by interpolating the
tuples in $\hat R_{h,k}^{l,m}$. The details can be found in
Appendix~\ref{app:remeshed}. This
interpolation strategy also yields a continuous extension $R_{h,k}^{l,m}$ of
$\hat R_{h,k}^{l,m}$, such that, for every $ (c^{(l)}, T^{(l)}, c^{(m)},
T^{(m)}) \in R_{h,k}^{l,m}$, equation~\eqref{eq:conncond} is satisfied.

For convenience, we will abuse the notation and call connection surfaces both
the finite and the continuous sets.
The one-to-one maps obtained by interpolating tuples in $\hRlmhk$ will be
denoted as
\begin{equation}\label{eq:trmaps}
   \tau^{l,m}_{h,k}:\mathcal{D}_{h,k}^{l,m} \rightarrow \mathcal{A}_{h,k}^{l,m}
   ,\qquad
   \left(\tau^{l,m}_{h,k}\right)^{-1}:
   \mathcal{A}_{h,k}^{l,m} \rightarrow \mathcal{D}_{h,k}^{l,m}
\end{equation}
and, following \cite{2013aDeGiRo}, will be denoted as
{\em transition maps}.

\begin{remark}\label{re:revtransmap}
In the heteroclinic case, when $l \ne m$, the inverse 
transition map $(\tau^{l,m}_{h,k})^{-1}$ does not, in general,
coincide with the
reverse transition map $\tau^{m,l}_{h,k}$, which is instead induced by the
reversibility property
\eqref{eq:rev}.
\end{remark}
 
\subsection{Relation with the scattering map}

In the same way that the $\Rlmhk$ sets allow us to define the transition maps
$\tau^{l,m}_{h,k}$, they also allow us to evaluate the corresponding
scattering maps as defined in \cite{2026aDeGiRo}. If $x\in\Gamma^{l,m}_{h,k}$,
and we define $x_-\in\Wcal^c_h(L_l)$ (resp.~$x_+\in\Wcal^c_h(L_m)$) as the
base of the fiber of $x$ within $W^{cu}(L_l)$ (resp.~$W^{cs}(L_m)$), the
scattering map sends $x_-$ to $x_+$. Since $x\in\Gamma^{l,m}_{h,k}$, by
\eqref{eq:conncond} there exists $(c^{(l)},T^{(l)},c^{(m)},T^{(m)})\in\Rlmhk$
such that
\begin{equation}\label{eq:condxR}
   \Phi_{T^{(l)}}\bigl(W^{(l)}(c^{(l)},\xi^{(l)})\bigr)
   =
   \Phi_{T^{(m)}}\bigl(W^{(m)}(c^{(m)},\xi^{(m)})\bigr)
   =x
   .
\end{equation}
By the uncoupling property of the parameterization (property
(\ref{en:propestil4}), see Section~\ref{ssec:param}), we know that
$W^{(l)}(c^{(l)},0) = W_c^{(l)}(c^{(l)}) $ is the base of the fiber of
$W^{(l)}(c^{(l)},\xi^{(l)})$.
By fiber invariance,
$
   \Phi_{T^{(l)}}\bigl(W_c^{(l)}(c^{(l)})\bigr)
$
is the base of the fiber of 
$
   \Phi_{T^{(l)}}\bigl(W^{(l)}(c^{(l)},\xi^{(l)})\bigr),
$
that, by \eqref{eq:condxR}, is equal to $x$. Then,
\[
   x_-=\Phi_{T^{(l)}}\bigl(W^{(l)}(c^{(l)},0)\bigr)
      =
      W^{(l)}\Bigl(\phi^{(l)}_{T^{(l)}}(c^{(l)},0)\Bigr)
      ,
\]
or, equivalently, using the reduced flow restricted to the center manifold,
\[
   x_-=\Phi_{T^{(l)}}\bigl(W_c^{(l)}(c^{(l)})\bigr)
      =
      W^{(l)}_c\Bigl(\varphi^{(l)}_{T^{(l)}}(c^{(l)})\Bigr)
      ,
   \]
where the last equality comes from the invariance equation \eqref{eq:invfl}.
So $x_-$ can be obtained by numerical integration through the reduced flow
$\varphi^{(l)}$ of the point $c^{(l)}$ given by the $\Rlmhk$ set through
Eq.~\eqref{eq:condxR}. An analogous argument shows that
\[
   x_+=W^{(m)}_c\Bigl(\varphi^{(m)}_{T^{(m)}}(c^{(m)})\Bigr)
   .
\]
Hence, in center manifold parameter coordinates, the scattering map $\sigma$
associated with the heteroclinic/homoclinic manifold described by
$\Gamma_{h,k}^{l,m}$ sends $c_-$ to $c_+$, where $x_-=W_c^{(l)}(c_-)$ and
$x_+=W_c^{(m)}(c_+)$. This map can be obtained by composing the backward reduced
flow from the departure CS, the transition map, and the forward reduced flow
from the arrival CS
\[
   \sigma = \varphi_{T^{(m)}}^{(m)} \circ \tau_{h,k}^{l,m} \circ
   \varphi_{-T^{(l)}}^{(l)},
\]
where $T^{(l)}$ and $T^{(m)}$ are the flight times associated with each
connection.

\subsection{Some homoclinic and heteroclinic connection surfaces} 
\label{ssec:some}

Reference \cite{2024BaHaMo} presented the CS of all primary heteroclinic
connections from ${\mathcal W}_h^c(L_2)$ to ${\mathcal W}_h^c(L_1)$ of the
Earth--Moon RTBP for several values of $h$.
In this work we will use two CS of that reference: one with $h=-1.585$, $k=6$
and another with $h=-1.575$, $k=2$.  We have remeshed these sets according to
the strategy of Appendix~\ref{app:remeshed}, in order to be able to use the
interpolation strategy and, in this way,
construct numerically the (continuous) 2D CS $\Dhk{2}{1}$, $\Ahk{2}{1}$ and
the transition maps $\tau_{h,k}^{2,1}$ between them provided by the
heteroclinic connections.

\begin{figure}[!htbp]
   \begin{subfigure}[tb]{.45\textwidth}
     \centering
     \includegraphics{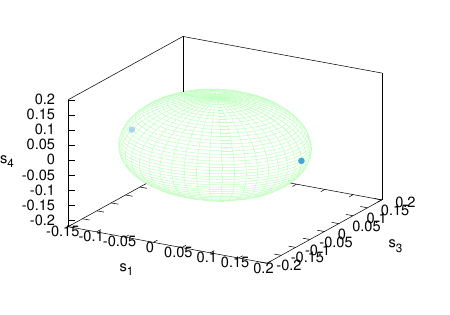}
     \caption{$\CalD{2}{1}{-1.585}{6}$\label{sfig:esfera}}
   \end{subfigure}
   \begin{subfigure}[tb]{.45\textwidth}
     \centering
     \includegraphics{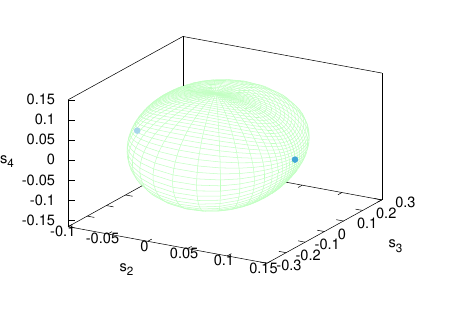}
     \caption{$\CalD{2}{1}{-1.575}{2}$}
     \label{sfig:esfera2}
   \end{subfigure}
   \begin{subfigure}[tb]{.45\textwidth}
     \centering
     \includegraphics{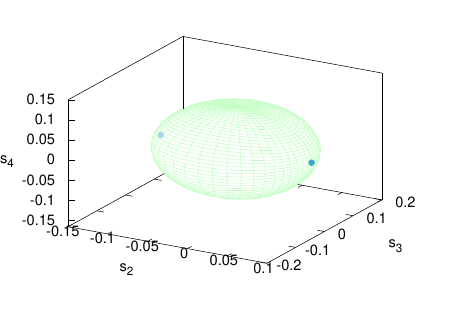}
     \caption{$\CalD{2}{1}{-1.500425}{4}$}
     \label{sfig:esfera3}
   \end{subfigure}
   \caption{Parametrization of the heteroclinic (from $L_2$ to $L_1$) 
   connection surfaces discussed in the text. Points
   represent the heteroclinic connections between planar Lyapunov p.o.
   Corresponding energy levels and parameters are given therein.}
   \label{fig:somehet}
\end{figure}

\begin{figure}[!htbp]
   \begin{subfigure}[tb]{.45\textwidth}
     \centering
     \includegraphics{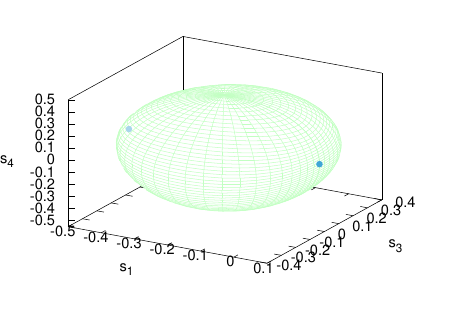}
     \caption{\label{sfig:homoesf}$\CalD{1}{1}{-1.585}{3}$}
     \label{sfig:esfera4}
   \end{subfigure}
   \begin{subfigure}[tb]{.45\textwidth}
     \centering
     \includegraphics{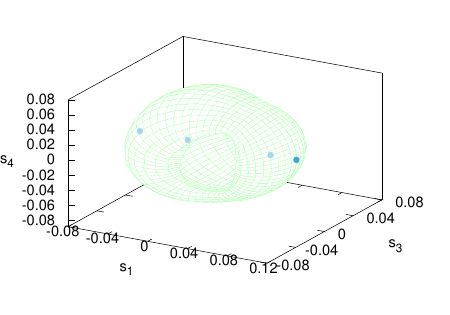}
     \caption{$\CalD{2}{2}{-1.585}{7}$\label{sfig:tor}}
   \end{subfigure}
   \caption{Parametrization of the homoclinic connection surfaces from $L_1$ and
   $L_2$ discussed in the text. Points
   represent the homoclinic connections between planar Lyapunov p.o.
   Corresponding energy levels and parameters are given therein.}
   \label{fig:somehom}
\end{figure}

For the benefit of a more complete discussion, in this work we have computed
some additional CS:
\begin{itemize}
\item
   The primary heteroclinic connections in the Sun-Earth system from $L_2$ to
   $L_1$ at $h=-1.500425$ with $k=4$.
\item
   The primary homoclinic connections between the center manifolds of $L_1$ and
   $L_2$ in the Earth-Moon system at $h=-1.585$ with $k=3$ and $k=7$,
   respectively.
\end{itemize}
All the homoclinic and heteroclinic sets mentioned are shown in
Figures~\ref{fig:somehet} and \ref{fig:somehom} and summarized in
Table~\ref{tab:summsomehomhet}.
Observe that the homoclinic CS shown in
Fig.~\ref{sfig:tor} 
is not a topological sphere but a
topological torus.
A discussion on the reason for this fact has been postponed to
Appendix~\ref{app:primarity}, together with additional discussion on non-primary
connection surfaces.

\begin{table}[ht]
   \centering
   \begin{tabular}{cccccc}
   \hline
   Mass parameter & $L_l$ & $L_m$ & $h$ & $k$ & Figure \\ \hline
   Earth-Moon & $2$ & $1$ & $-1.585$ & $6$ & \ref{sfig:esfera} \\
   Earth-Moon & $2$ & $1$ & $-1.575$ & $2$ & \ref{sfig:esfera2} \\
   Sun-Earth  & $2$ & $1$ & $-1.500425$ & $4$ & \ref{sfig:esfera3} \\
   Earth-Moon & $1$ & $1$ & $-1.585$ & $3$ & \ref{sfig:esfera4} \\
   Earth-Moon & $2$ & $2$ & $-1.585$ & $7$ & \ref{sfig:tor} \\
   \end{tabular}
   \caption{Summary of the computed global sets of connections presented in
   Figure~\ref{fig:somehet} and \ref{fig:somehom}.}
   \label{tab:summsomehomhet}
\end{table}

\section{Connectivity at the level of tori}
\label{sec:maps}

As it has been mentioned in the previous section, the $\Rlmhk$ sets 
provide fine grids of both
sections of isoenergetic heteroclinic (resp.~homoclinic) manifolds,
$\Gamma^{l,m}_{h,k}$, and CS
sets, $\Dhk{l}{m}$ and $\Ahk{l}{m}$, $l,m \in \{1,2\}$. In this way, they
provide connectivity information through heteroclinic (resp.~homoclinic)
connections at the level of points in both $\calW_h^c(L_l)$ and
$\calW_h^c(L_m)$. For applications, it would be interesting to have
connectivity information at the level of tori as well. For example, given a
torus of $\calW^c_h(L_2)$, we would like to know the range of tori of
$\calW^c_h(L_1)$ (resp.~$\calW^c_h(L_2)$) that are reachable through
heteroclinic (resp.~homoclinic) connections. The goal of this section is to
produce such connectivity maps.
Examples of astrodynamical interest using this tool are also given at the end
of the section.

\subsection{Strategy}\label{sec:estrkvk}

Before addressing the computation of these connectivity maps, we first need a
way to identify a torus within an isoenergetic family of $\calW^c_h(L_j)$.  For
that, we will make use of the isoenergetic Poincaré sections $\poincpi$ of the
center manifold depicted in Figure~\ref{fig:poinc}. We distinguish two cases,
depending on whether they contain Lissajous tori only or both Lissajous and
quasi-Halo. These are, in fact, the only possible scenarios within the domain
of validity of the expansions provided by the parameterization method (see
Section~\ref{ssec:param}).

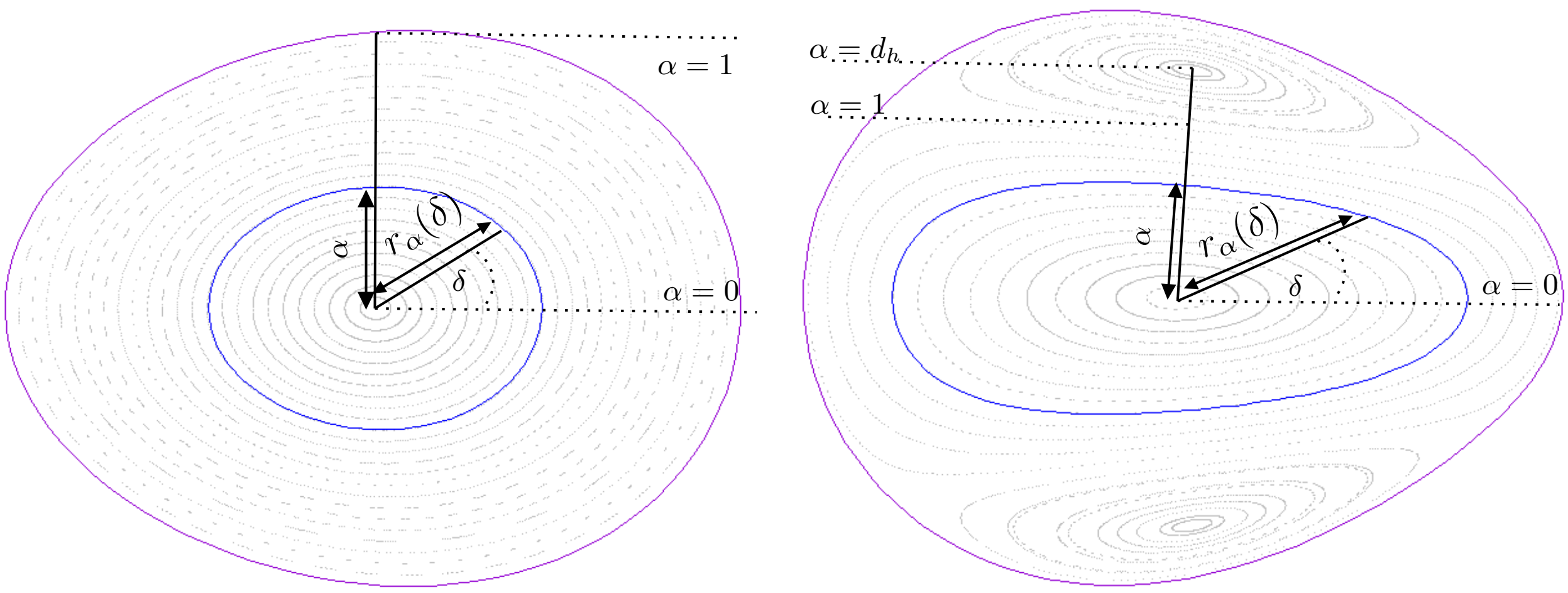
\begin{figure}[htb]
\input{paramsinvcurve.tex}
\caption{\label{fig:alphas}Schematic representation of the parameterization of
the invariant curves in the $\poincpi$ Poincaré section and their
identification with $\alpha$.}
\end{figure}

Each torus is identified by the invariant curve it produces in its
isoenergetic $\poincpi$ Poincaré section. In turn, this
invariant curve is characterized by the distance from the origin to its
intersection with a reference line, defined as follows.  In the {\em
Lissajous-only} case (Fig.~\ref{fig:alphas} left), the reference line is the
(vertical) $s_2$-axis.  In the {\em Lissajous-and-quasi-Halo} case
(Fig.~\ref{fig:alphas} right), it will be the one joining the fixed points
corresponding to the Halo p.o.~and the vertical Lyapunov p.o. We will denote as
$\alpha$ a rescaled version of this distance. The rescaling is chosen so that
$\alpha=1$
corresponds either to the planar Lyapunov p.o., in the Lissajous-only case, or
to the separatrix between Lissajous and quasi-Halo tori, in the
Lissajous-and-quasi-Halo one. We denote by $d_h$ the maximum value that
$\alpha$ can attain in either case. Thus, in the Lissajous-only case, $d_h=1$.
In the Lissajous-and-quasi-Halo case, $d_h$
is the rescaled distance between the fixed points associated with the Halo and
vertical Lyapunov p.o., so that the Halo p.o.~satisfies $\alpha=d_h$.
Consequently, $0<\alpha<1$ corresponds to Lissajous tori, while $1<\alpha<d_h$
identifies quasi-Halo tori.  See Fig.~\ref{fig:alphas}.

Consider now an invariant curve given by $\alpha$ according to the previous
discussion.  Assume that its backbone p.o.~is located at the fixed point
$(s_1^{(0)}, s_2^{(0)})$ on the isoenergetic Poincaré section (if this backbone
p.o.~is the vertical Lyapunov one then $s_1^{(0)}=s_2^{(0)}=0$).  We want to
parameterize it in polar coordinates with the angular variable measured from
the positive $s_1$-axis,
and therefore given by a function $r_\alpha(\delta)$ as
$\{c(\delta,r_\alpha(\delta))\}_{0\leq \delta < 2\pi}\subset\calV^{(j)}$, where
$c(\delta,r)=\bigl(s_i(\delta, r)\bigr)_{i=1}^4$ is defined by
\begin{equation}\label{eq:polars}
   \begin{aligned}
      s_1(\delta,r) &= r\cos\delta + s_1^{(0)}, \\
      s_2(\delta,r) &= r\sin\delta + s_2^{(0)}, \\
      s_4(\delta,r) &= 0, \\
      s_3(\delta,r) &\ \lefteqn{\mbox{found from
      $H\bigl(W_c^{(j)}\bigl(c(\delta,r)\bigr)\bigr)=h$,}}
   \end{aligned}
\end{equation}
see Figure~\ref{fig:alphas}.
By starting from
the intersection of the invariant curve with the reference line (see
Fig.~\ref{fig:alphas}) and computing $\poincpi$ Poincaré sections of the flow
$\varphi_t^{(j)}$, a (finite) sequence
$\{\bigl(\delta_k,r_\alpha(\delta_k)\bigr)\}_k$ is obtained by inverting
\eqref{eq:polars} at the Poincaré cuts. When the sequence $\{\delta_k\}_k$
samples the interval $[0,2\pi]$ densely enough, $r_\alpha(\delta)$ can be
recovered for arbitrary $\delta$ by interpolation.  Once this reconstruction
is available, a truncated Fourier series for $r_\alpha(\delta)$ can be
obtained by choosing $N$ and computing the Discrete Fourier Transform of
$\{r_\alpha(i\frac{2\pi}N)\}_{i=0}^{N-1}$. We could also compute their
rotation number with frequency analysis methods (see
e.g.~\cite{2010GoMoSi,2020SaMe,2021MeSa}).

\begin{remark}
   From this approximation, we could generate a 2D invariant torus for the flow
   with the corresponding frequencies, and apply the methods developed in
   \cite{2021aHaMo}.
\end{remark}

Denoting by $P$ the Poincaré map corresponding to the section $\poincpi$, the
error of such a truncated Fourier series $r_\alpha(\delta)$ can be estimated
as
\begin{equation}\label{eq:este}
   e_\alpha := \max_{ 0 \leq i \leq N-1 } |P
   \bigl(c(\delta_i,r_\alpha(\delta_i))\bigr) - c(\theta_i,
   r_\alpha(\theta_i))|,
\end{equation}
where $\delta_i=i\frac{2\pi}N$ and $\theta_i$ denotes
the polar angle of the first two coordinates of
$P\bigl(c(\delta_i,r_\alpha(\delta_i))\bigr)$.

Finally, in order to design a procedure to look for tori-to-tori connections,
we
require a continuous description with respect to $\alpha$ of the invariant
curves corresponding to the entire isoenergetic families of tori under
consideration.
For that, we first choose a set $\{\alpha_i\}_{i=0}^{M-1}$ of values of
$\alpha$, with $\alpha_i=i/M$, if the invariant curve corresponds to a
Lissajous torus, or $\alpha_i=d_h-i(d_h-1)/M$, if it corresponds to a
quasi-Halo torus.  For each family, we compute the truncated Fourier series
for $r_{\alpha_i}(\delta)$ as described before. Then, for arbitrary
$\alpha\in[0,d_h]\backslash\{1\}$ and $\delta\in[0,2\pi]$, we compute
$r_\alpha(\delta)$ by first interpolating the Fourier coefficients of
$r_\alpha$ from the ones of $r_{\alpha_i}$ for nearby values of $\alpha_i$ of
the same family, and subsequently evaluating the Fourier series obtained.
Denoting again by $r_\alpha(\delta)$ the Fourier series obtained in this way,
the approximation error of this interpolation strategy on the family (either
Lissajous or quasi-Halo) as a whole can be estimated as
\begin{equation}\label{eq:estE}
   E := \max_{ 0 \leq i \leq M-2 }
      e_{(\alpha_i+\alpha_{i+1})/2}
      .
\end{equation}
The parameters $N$ and $M$
will be selected so
as to ensure that, for a given family of invariant tori, the estimates
$\{e_{\alpha_i}\}_{i=0}^{M-1}$ and $E$ fall below a prescribed tolerance.

Now, for $j=1,2$, consider isoenergetic families with energy $h$ of tori
around $L_j$ with $\{s_4=0\}$ invariant curves given by
$
   \{r^{(j)}_\alpha(\delta)\}_{
   \alpha\in[0,d_h^{(j)}]\backslash\{1\},\delta\in[0,2\pi] },
$
according to the previous discussion.
We will denote as
\begin{equation}\label{eq:deftoralfa}
   \{\toralfa jh\alpha\}_{
      \alpha\in[0,d_h^{(j)}]\backslash\{1\},\delta\in[0,2\pi]
   },\quad
   \toralfa jh\alpha\subset\calV^{(j)},
\end{equation}
the corresponding tori.
For $l,m \in \{1,2\}$ and $\alpha^{(l)}\in[0,d_h^{(l)}]\backslash\{1\}$, we
would like to know the tori within $\Wcal^c_h(L_m)$ that
$\toralfa lh{\alpha^{(l)}}$
is connected to.  For that, we proceed as
follows:
\begin{enumerate}
\item
   Take a dense enough, equally-spaced grid of values
   $\{\delta_k\}_{k=0}^{n-1}\subset[0,2\pi)$.
\item\label{en:ctrtr2}
   For $k=0,\dots,n-1$, propagate
   $c\bigl(\delta_k,r^{(l)}_{\alpha^{(l)}}(\delta_k)\bigr)$ through the flow
   $\varphi_t^{(l)}$ until the next Poincaré cut with $\poincpi$. On the way,
   record its possible intersections (one or two) with the CS $\Dhk{l}{m}$.
\item\label{en:ctrtr3}
   For each intersection with $\Dhk{l}{m}$ obtained in the previous
   step, obtain the corresponding point in $\Ahk{l}{m}$ through
   the transition map $\tau_{h,k}^{l,m}$ (see Eq.~\eqref{eq:trmaps}).
\item\label{en:ctrtr4}
   Propagate the points $c$ obtained in the previous step through the flow
   $\varphi_t^{(m)}$ until the first $\poincpi$ Poincaré section. If this
   intersection is given by $c_*=c(\delta_*,r_*)$, find $\alpha^{(m)} \in
   [0,d_h^{(m)}]\backslash\{1\}$ such
   that $r_*=r^{(m)}_{\alpha^{(m)}}(\delta_*)$ through
   inverse interpolation of
   $\{\bigl(\alpha_k,r^{(m)}_{\alpha_k}(\delta_*)\bigr)\}_k$ for
   $\{\alpha_k\}_k$ with $\{r_{\alpha_k}(\delta_*)\}_k$ near $r_*$.
   Then,
   $\toralfa mh{\alpha^{(m)}}$
   is connected to
   $\toralfa lh{\alpha^{(l)}}$.
\end{enumerate}
The numerical evaluation of Poincaré sections with respect to the CS
$\Dhk{l}{m}$
can be performed using the (finite) results set $\Rlmhk$ by means of the
interpolation strategy described in Appendix~\ref{app:globalpoinc}.

The output of iterating these steps over a grid of values $\alpha^{(l)} \in
[0,d_h^{(l)}]\backslash\{1\}$
allows us to approximate
the {\em connectivity map},
that we define as the correspondence
between sets of $\alpha^{(j)}$ values
\[
   F_{h,k}^{l,m} : {\mathcal P}\bigl([0,d_h^{(l)}]\backslash\{1\}\bigr)
   \longrightarrow
   {\mathcal P}\bigl([0,d_h^{(m)}]\backslash\{1\}\bigr),
\]
where $\mathcal P$ denotes the set of all subsets, and $F_{h,k}^{l,m}$ is
defined by the property:
$F_{h,k}^{l,m}({\mathcal
A})={\mathcal B}$ if and only if
for any pair of $(\alpha^{(l)},\alpha^{(m)})\in\mathcal A\times\mathcal B$,
the corresponding tori $\toralfa lh{\alpha^{(l)}}$,
$\toralfa mh{\alpha^{(m)}}$
are heteroclinically/homoclinically connected.
By definition, this correspondence is invertible.
Since, for any pair of tori $\Tcal^{(l)}\subset\Wcal^c(L_l)$,
$\Tcal^{(m)}\subset\Wcal^c(L_m)$ of the families
under consideration, the reversibility property \eqref{eq:rev} turns a
heteroclinic connection from $\Tcal^{(l)}$ to $\Tcal^{(m)}$ into one from
$\Tcal^{(m)}$ to $\Tcal^{(l)}$, it turns out that 
\begin{equation}\label{eq:invcon}
   (F_{h,k}^{l,m})^{-1}=F_{h,k}^{m,l}.
\end{equation}
Recall that this is not true for transition maps (see
Remark~\ref{re:revtransmap}).

\subsection{Connectivity maps}
\label{ssec:conmaps}

In this section, we present connectivity maps computed according to the
discussion in the previous section. They correspond to the mass ratios and
energy values $\{h_i\}_{i=1}^3$ given in Table~\ref{tab:scen}. They have been
chosen in order to
cover the following possibilities: Lissajous tori only both at departure and
on arrival (Case 1); Lissajous tori only at departure (or arrival) but
Lissajous and quasi-Halo on arrival (or departure, Case 2); Lissajous and
quasi-Halo both at departure and on arrival (Case 3).

\begin{table}
   \centering
   \begin{tabular}{c|cccccc}
      Case ($i$) & Mass ratio & Energy ($h_i$) & $\Wcal^c_h(L_2)$ & $\Wcal^c_h(L_1)$ &
      $d_h^{(2)}$ & $d_h^{(1)}$\\
      \hline 
      1 & Sun-Earth & $-1.500425$ & Lissa & Lissa & $1$ & $1$\\
      2 & Earth-Moon & $-1.585$ & Lissa & Lissa+qH & $1$ & $1.035$\\
      3 & Earth-Moon & $-1.575$ & Lissa+qH & Lissa+qH & $1.012$ & $1.17$\\
   \end{tabular}

   \medskip

   \caption{Considered scenarios for the computation of connectivity maps
   according to the
   families of invariant tori involved in each
   isoenergetic center manifold: Lissajous (Lissa) or quasi-Halo (qH).
   }

   \label{tab:scen}
\end{table}

In all the connectivity maps that will be shown, the $N$ discretization
parameter considered in the previous Section has been set to 128 for Lissajous
tori and 1024 for quasi-Halo ones. For higher accuracy, the corresponding
Fourier series have been constructed from the first half of DFT coefficients.
With these choices, we have obtained that the $e_\alpha$ estimate defined in
Eq.~\eqref{eq:este} is under $10^{-8}$ for all the values of $\alpha$
considered. As for the $M$ parameter, numerical experiments have shown that
$M\geq20$ Lissajous tori and $M\geq140$ for quasi-Halo are enough in order to
ensure that the $E$ estimate of Eq.~\eqref{eq:estE} is below $10^{-8}$.
Nonetheless, $M=500$ has been taken in all cases in order to reduce the gaps
around $\alpha=1$ in the connectivity maps, that correspond to the separatrix
between the Lissajous and quasi-Halo families of tori (see
Section~\ref{ssec:dynwc}).

\begin{figure}[tbp]
   \centering
   \includegraphics{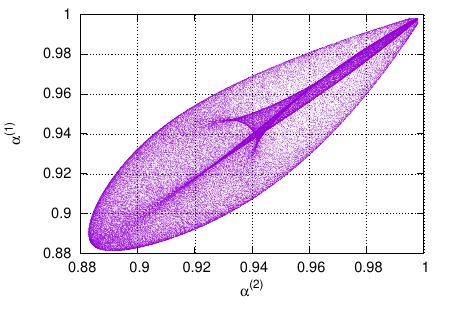}
   \caption{Primary heteroclinic connectivity map $F_{h_1,4}^{2,1}$
   for Case~$1$ in
   Table~\ref{tab:scen}.
   The horizontal axis represents the value of $\alpha$ corresponding to the
   departing invariant tori, while the vertical axis represents the value of
   $\alpha$ corresponding to the arriving invariant tori.}
   \label{fig:kvk_se_-1.500425}
\end{figure}

Figure~\ref{fig:kvk_se_-1.500425} shows the connectivity map of
$4$-heteroclinic connections between $\Wcal_{h_1}^c(L_2)$ and
$\Wcal_{h_1}^c(L_1)$ for Case~1 in Table~\ref{tab:scen} (Sun-Earth mass
ratio), where only Lissajous tori are present in both isoenergetic center
manifolds.
The values of $\alpha^{(j)}$ close to $0.94$ correspond to the most highly
connected tori. The plot is
generated from connections from ${\mathcal W}^c_{h_1}(L_2)$ to ${\mathcal
W}^c_{h_1}(L_1)$.
But, since $(F^{2,1}_{h_1,4})^{-1}=F^{1,2}_{h_1,4}$ (see
Eq.~\eqref{eq:invcon}), the plot can be read either from the horizontal axis to
the vertical one, in order to find $L_2\rightarrow L_1$ connections, or
vice-versa, in order to find $L_1\rightarrow L_2$ connections.

\begin{figure}[tbp]
   \centering
   \includegraphics{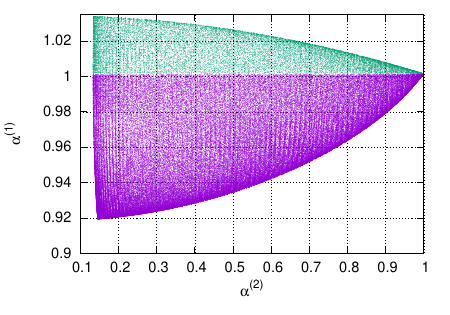}
   \caption{ Primary heteroclinic connectivity map $F_{h_2,6}^{2,1}$ for Case~2 in
   Table~\ref{tab:scen}.
   Purple: Lissajous--Lissajous connections.
   Green: Lissajous--Southern (Northern) quasi-Halo connections.
   }
   \label{fig:kvk_tl_-1.585}
\end{figure}

Figure~\ref{fig:kvk_tl_-1.585} presents the connectivity map corresponding to
heteroclinic connections for Case~2 in Table~\ref{tab:scen}, which involves
both Lissajous and quasi-Halo tori in $\Wcal_{h_2}^c(L_1)$. The Southern
$L_1$ quasi-Halo family is represented by $\alpha^{(1)}>1$. In this setting,
tori near the separatrix between both families emerge as the most connected
ones in $\Wcal_{h_2}^c(L_1)$. Moreover, since $d_h^{(1)}=1.035$, the top left
corner of the plot corresponds to a connection between the $L_1$ Halo p.o.~and
an $L_2$ Lissajous torus. In $\Wcal_{h_2}^c(L_2)$, connectivity improves
as $\alpha^{(2)}$ decreases (that is, as tori approach the vertical Lyapunov
p.o.), except for a small range close to the smallest values of $\alpha^{(2)}$
 in Figure~\ref{fig:kvk_tl_-1.585}, in which connectivity decreases sharply.
Since symmetry~\eqref{eq:sym} leaves Lissajous tori invariant but switches
Northern quasi-Halo tori to Southern ones and vice-versa, for each $L_2$
Lissajous $\rightarrow$ $L_1$ Southern quasi-Halo connection there is a
corresponding $L_2$ Lissajous $\rightarrow$ $L_1$ Northern quasi-Halo one.
Also,
thanks to \eqref{eq:invcon}, Figure~\ref{fig:kvk_tl_-1.585} can be used to
find both $L_2\rightarrow L_1$ and $L_1\rightarrow L_2$ connections.

\begin{figure}[tbp]
   \centering
   \includegraphics{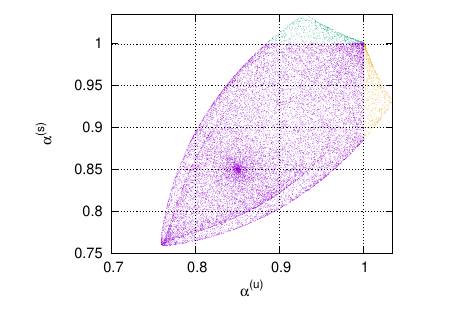}
   \includegraphics{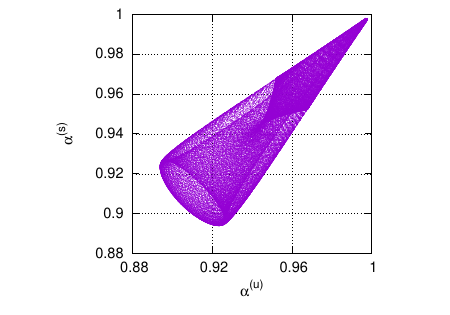}
   \caption{
      Homoclinic connectivity maps for Case 2 in Table~\ref{tab:scen}.
      Left: 
      primary connectivity map $F_{h_2,3}^{1,1}$.
      Right: primary 
      connectivity map $F_{h_2,7}^{2,2}$.
      Purple: Lissajous--Lissajous connections. Green:
      Lissajous--quasi-Halo connections. Orange: quasi-Halo--Lissajous
      connections.}
   \label{fig:homos}
\end{figure}

Figure~\ref{fig:homos} presents results on homoclinic
connections corresponding to Case 2 in Table~\ref{tab:scen} and represented in
Figure~\ref{fig:somehom}. Our motivation in
presenting these results is twofold. On the one hand, it illustrates the fact
that the methodology presented in Section~\ref{sec:estrkvk} works with
homoclinic connections. On the other hand, even when looking to connect tori
around different libration points, considering homoclinic connections together
with heteroclinic ones can enlarge the range of invariant tori around a
libration point that are accessible from tori of the other. In this case, we
refer to parameters $\alpha^{(u)},\alpha^{(s)}$ to indicate departure and
arrival tori in the same $\Wcal_h^c(L_j)$.
Figure~\ref{fig:homos} left involves tori with $\alpha^{(u)},\alpha^{(s)}$
reaching $0.76$, whereas the smallest
$\alpha^{(1)}$ reachable in Figure~\ref{fig:kvk_tl_-1.585} is
$0.92$. An example pseudo-trajectory exploiting this fact is
given in Section~\ref{ssec:traj}. Note that, due to \eqref{eq:invcon}, all the
plots in Figure~\ref{fig:homos} are symmetric with
respect to the $\alpha^{(u)}=\alpha^{(s)}$ line.

Figure~\ref{fig:homos} left has been generated by looking for $3$-homoclinic
connections of $\Wcal^c(L_1)$. The points with $\alpha^{(u)}>1$ or
$\alpha^{(s)}>1$ can be thought to correspond to connections between a
Lissajous torus and either a Southern or Northern quasi-Halo one, because
symmetry \eqref{eq:sym} turns one kind of connection into the other. No
connections between quasi-Halo tori have been found, neither of the same
subfamily (Northern, Southern) or between different ones.
The right plot of Figure~\ref{fig:homos} represents homoclinic
connections of $\Wcal^c(L_2)$, that in this case involve the Lissajous family
only. It can be seen that the range of tori connected is rather limited,
although it could be increased slightly by increasing $k$. A comparison with
Figure~\ref{fig:kvk_tl_-1.585} shows that, by combining two heteroclinic
connections instead of just using homoclinic ones, a much larger change in
$\alpha$ can be obtained. A sample pseudo-trajectory illustrating this fact is
also given in Section~\ref{ssec:traj}.

\begin{figure}[tbp]
   \centering
   \includegraphics{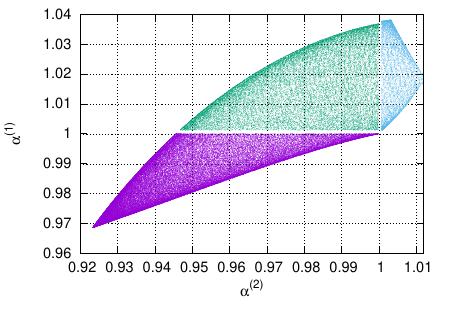}
   \caption{Primary heteroclinic connectivity map $F_{h_3,2}^{2,1}$ for Case~3 in
   Table~\ref{tab:scen}.
   Purple: Lissajous--Lissajous connections.
   Green: Lissajous--Southern (Northern) quasi-Halo connections.
   Blue: Northern (Southern) quasi-Halo--Southern (Northern) quasi-Halo
   connections.
   }
   \label{fig:kvk_tl_-1.575}
\end{figure}

Figure~\ref{fig:kvk_tl_-1.575} shows results on heteroclinic connections for
Case 3 in Table~\ref{tab:scen}, in which both Lissajous and quasi-Halo tori
appear in the isoenergetic slices around both $L_1$ and $L_2$.  In this case,
there are no connections from $L_2$--quasi-Halo to $L_1$--Lissajous tori, nor
within quasi-Halo families of the same type ($L_2$--Northern to
$L_1$--Northern or $L_2$--Southern to $L_1$--Southern). Instead, the entire
range of quasi-Halo orbits in $\Wcal_{h_3}^c(L_2)$ ($d_{h_3}^{(2)}=1.012$, see
Table~\ref{tab:scen}), including the Halo p.o., are connected (to the $L_1$
quasi-Halo family  of the opposite type), while only a limited range of
quasi-Halo orbits in $\Wcal_{h_3}^c(L_1)$ ($d_{h_3}^{(1)}=1.17$, see
Table~\ref{tab:scen}) exhibit connections. The $L_2$ separatrix once again
emerges as the most highly connected. In the connections
used to generate Figure~\ref{fig:kvk_tl_-1.575}, departing quasi-Halos
have always been
southern and ending quasi-Halos
northern.
However, again from symmetry \eqref{eq:sym}, corresponding connections exist
with departing Northern quasi-Halos and ending southern ones. Also
thanks to \eqref{eq:invcon}, for each $L_2\rightarrow L_1$ connection in
Figure~\ref{fig:kvk_tl_-1.575} there is a corresponding $L_1\rightarrow L_2$
with departing and end tori reversed.

\subsection{An example application to trajectory design}
\label{ssec:traj}

By leveraging the dynamics of heteroclinic and homoclinic connections, the
connectivity maps introduced in Section~\ref{ssec:conmaps} can be exploited to
construct pseudo-trajectories, meaning the juxtaposition of trajectories
within $\Wcal_h^c(L_l)$, $l=1,2$, and homo/heteroclinic connection arcs. They
can be thought of as given by finite sequences of points in center manifold
parameters
\[
   c_0^{(l_0)}, c_1^{(l_1)}, \dots, c_N^{(l_N)},
\]
with $l_i\in\{1,2\}$ and $c_i\in\calV^{(l_i)} \subset \R^4$ for $0\leq i\leq N$. The
pseudo-trajectory condition is that, for $0\leq i \leq N-1$,
\begin{align}
   \mbox{either }\exists t_i>0
      &: c_{i+1}^{(l_{i+1})}=\varphi_{t_i}^{(l_i)}(c_i^{(l_i)})
         \quad\mbox{($l_i=l_{i+1}$ in this case)},
      \label{eq:pstrc1}
      \\
   \mbox{or }\exists k_i
   &: c_{i+1}^{(l_{i+1})}=\tau^{l_i,l_{i+1}}_{h,k_i}(c_i^{(l_i)}),
      \label{eq:pstrc2}
\end{align}
where $\varphi_t^{(j)}$ is the flow restricted to the center manifold.
In case \eqref{eq:pstrc2} above, there exists $t_i>0$ such that
$\Phi_{t_i}\bigl(W^{(l_i)}(c_{i}^{(l_i)},\xi^{(l_i)})\bigr)=
W^{(l_{i+1})}(c_{i+1}^{(l_{i+1})},\xi^{(l_{i+i})})$, that is, $t_i$ is the flight
time along the
homo/heteroclinic connection. It can be computed as
$t_i=T^{(l_i)}-T^{(l_{i+1})}$, if
$(c_i^{(l_i)},T^{(l_i)},c_{i+1}^{(l_{i+1})},T^{(l_{i+1})})\in
R^{l_i,l_{i+1}}_{h,k_i}$.
The $k_i$ parameter can be used to distinguish between cases \eqref{eq:pstrc1}
and \eqref{eq:pstrc2} above by setting it to 0 in case \eqref{eq:pstrc1}.
Then, the pseudo-trajectory can be constructed through numerical integration
of the $\varphi^{(1)}_t$, $\varphi^{(2)}_t$, $\Phi_t$ flows from the data
\begin{equation}\label{eq:datapseudo}
   \{c_i^{(l_i)}\}_{i=0}^N,\quad\{k_i\}_{i=0}^{N-1},\quad\{t_i\}_{i=0}^{N-1}
   .
\end{equation}
The total flight time of the pseudo-trajectory is obtained as
$\sum_{i=0}^{N-1}t_i$.

These pseudo-trajectories can be tailored to achieve different objectives
depending on the requirements of a space mission. In this section
we present a systematic framework for constructing such trajectories and we
illustrate it through several representative examples based on the
connectivity maps introduced previously.

Given $\alpha^{(l)}\in[0,d_h^{(l)}]\backslash\{1\}$, the connectivity maps
displayed in Figures~\ref{fig:kvk_se_-1.500425} to \ref{fig:kvk_tl_-1.575}
allow us to determine an interval of values of $\alpha^{(m)}$ for which the
tori $\toralfa lh{\alpha^{(l)}}$, $\toralfa mh{\alpha^{(m)}}$ are
homoclinically/heteroclinically connected.  By choosing $\alpha^{(m)}\in
F^{l,m}_{h,q}(\{\alpha^{(l)}\})$ and iterating this procedure, we can
construct a sequence of connected tori given by
\[
   \alpha_0^{(m_0)}, \alpha_1^{(m_1)}, \dots, \alpha_n^{(m_n)},
\]
with $m_i\in\{1,2\}$ and $\alpha_{i+1}^{(m_{i+1})}\in
F^{m_i,m_{i+1}}_{h,q_i}(\{\alpha_i^{(m_i)}\})$ for some $q_i$. For instance,
direct inspection of Fig.~\ref{fig:kvk_se_-1.500425} shows that, using
$F^{2,1}_{h_1,4}$, by following a staircase-like pattern the following
sequence of connected tori is possible:
\begin{equation}\label{eq:chainalphas}
   \alpha_0^{(2)}\approx0.980,
   \alpha_1^{(1)}\approx0.960,
   \alpha_2^{(2)}\approx0.928,
   \alpha_3^{(1)}\approx0.905,
   \alpha_4^{(2)}\approx0.898,
\end{equation}
where $n=4$, $m_0=m_2=m_4=2$, $m_1=m_3=1$, and $q_i=4$ for $1\leq i\leq4$.  We
would like to construct a pseudo-trajectory that follows these tori, given by
$\{c_i^{(l_i)}\}_{i=0}^N$ with $l_0=l_N=2$, $c_0^{(l_0)}\in\toralfa
2{h_1}{\alpha_0^{(2)}}$
and $c_N^{(l_N)}\in\toralfa 2{h_1}{\alpha_4^{(2)}}$.

We will be able to build such a pseudo-trajectory if, for $q_i$, $l_i$,
$l_{i+1}$ fixed, given $\alpha_i^{(l_i)}\in[0,d_h^{(l_i)}]\backslash\{1\}$ and
an initial condition $c_i^{(l_i)}\in\toralfa
{l_i}h{\alpha_i^{(l_i)}}\subset\calV^{(l_i)}$, we can find
$\alpha^{(l_{i+1})}\in F^{l_i,l_{i+1}}_{h,q_i}(\{\alpha_i^{(l_i)}\})$ within
the given constraints (in our case, being close to the next value in the chain
\eqref{eq:chainalphas}) and a pseudo-trajectory joining $\toralfa
{l_j}h{\alpha^{(l_j)}}$, $\toralfa {l_{j+2}}h{\alpha^{(l_{j+2})}}$ given by
the three-element sequence
$c_j^{(l_j)},c_{j+1}^{(l_{j+1})},c_{j+2}^{(l_{j+2})}$, with $l_j=l_{j+1}=l_i$,
$l_{j+2}=l_{i+1}$ and associated data $\{t_j,t_{j+1}\}$, $\{k_j,k_{j+1}\}$
(see \eqref{eq:datapseudo}) to be determined.  This can be done as follows:
\begin{enumerate}
\item
   Set $c\leftarrow c_i^{(l_i)}$, $k_j\leftarrow0$, $k_{j+1}\leftarrow q_i$,
   $t_j\leftarrow 0$.
\item\label{en:compseud2}
   Propagate $c$ under the restricted flow to the center manifold $\varphi_t^{(l_i)}$
   until the trajectory
   intersects the CS $\Dcal^{l_i,l_{i+1}}_{h,q_i}$. Set $c_{j+1}^{l_{j+1}}$
   equal to this intersection.
\item
   Compute $c_{j+2}^{(l_{j+2})}:=
   \tau^{l_{j+1},l_{j+2}}_{h,k_{j+1}}(c_{j+1}^{(l_{j+1})})$ through the
   interpolation strategy of Appendix~\ref{app:remeshed}.
\item\label{en:compseud4}
   Determine the parameter
   $\alpha^{(l_{i+1})}_{i+1}\in[0,d_h^{(l_{i+1})}]\backslash\{1\}$ of the
   corresponding arrival torus, that is, such that such that
   $c_{j+2}^{(l_{j+2})}\in\toralfa {l_{i+1}}h{\alpha^{(l_{i+1})}_{i+1}}$. This
   can be done through step (\ref{en:ctrtr4}) of the procedure introduced in
   Section~\ref{sec:estrkvk} for the computation of connectivity maps.
\item
   Set $t_j\leftarrow t_j+t$.
\item
   If $\alpha^{(l_{i+1})}_{i+1}$ satisfies the required constraints, set
   $t_{j+1}=T^{(l_{j+1})}-T^{(l_{j+2})}$, where
   $(c_{j+1}^{(l_{j+1})},T^{(l_{j+1})},c_{j+2}^{(l_{j+2})},T^{(l_{j+2})})\in
   R_{h,k_{j+1}}^{l_{j+1},l_{j+2}}$ and stop.
\item
   Otherwise, set $c\leftarrow c_{j+1}^{l_{(j+1)}}$ and go to step
   \eqref{en:compseud2}.
\end{enumerate}
The computational cost of Step~\eqref{en:compseud4} is reduced by extending the
results datasets
$\Rlmhk$ to include the departure and arrival parameters $\alpha^{(l)}$,
$\alpha^{(m)}$. Then, $\alpha^{(l_{i+1})}_{i+1}$ in step~\eqref{en:compseud4}
can be found from interpolation of the data in $\Rlmhk$ instead of following
step~\eqref{en:ctrtr4} in Section~\ref{sec:estrkvk}.

Applying the previous procedure on each segment of the chain
\eqref{eq:chainalphas}, the desired pseudo-trajectory
$\{c_i^{(l_i)}\}_{i=0}^N$ is constructed, with $N=8$ and additional data
$\{k_i\}_{i=0}^7$, $\{t_i\}_{i=0}^7$. The relationship between the data
of the pseudo-trajectory and the sequence of connected tori
can be sketched as follows:\\
\begin{center}
   \input{sketch.tex}
\end{center}
In the lower row, arrows with $k_j=0$ correspond to center manifold arcs,
whereas arrows with $k_j=4$ correspond to connection arcs.

\begin{figure}[tb]
   \centering
   \includegraphics{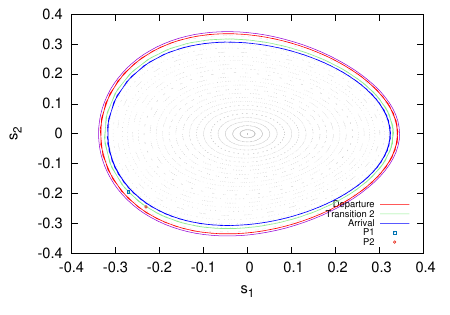}
   \includegraphics{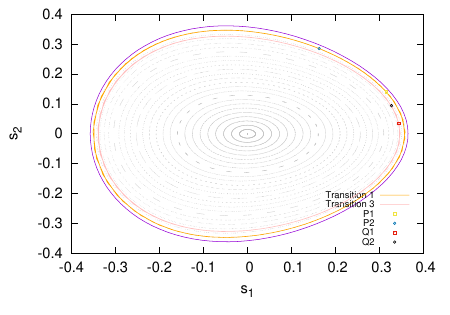}
   \includegraphics{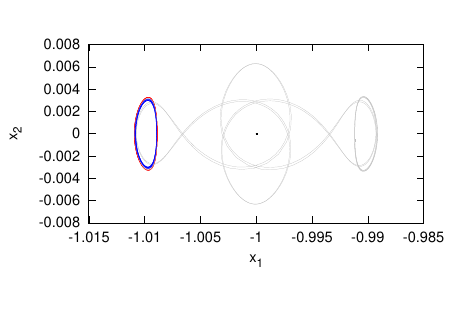}
   \includegraphics{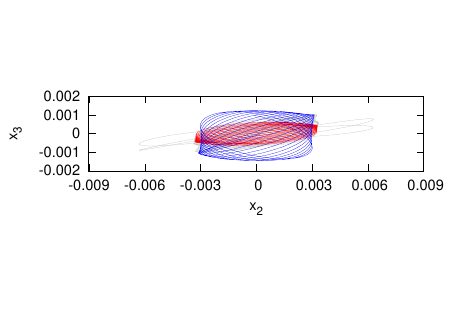}
   \caption{
      Pseudo-trajectory constructed via successive applications of the
      primary connectivity map $F_{h_1,4}^{2,1}$ shown in
      Figure~\ref{fig:kvk_se_-1.500425}. Four heteroclinic excursions increase
      the inclination from $\alpha^{(2)}_0=0.980$ to $\alpha^{(2)}_4=0.898$,
      passing through the intermediate tori $\alpha^{(1)}_1=0.960$,
      $\alpha^{(2)}_2=0.928$ and $\alpha^{(1)}_3=0.905$.  Points $Q_1,P_1$ and
      $Q_2,P_2$ in the $\poincpi$ section indicate the integration intervals
      preceding each transition.
   }
   \label{fig:acrob_h1}
\end{figure}

Numerical integration of each arc with the corresponding flow
($\phi^{(l_j)}_t$ if $k_j=0$, $\Phi_t$ otherwise) allows us to generate the
whole pseudo-trajectory, whose $x_1x_2$ and $x_1x_3$ projections are displayed
in the lower row of Figure~\ref{fig:acrob_h1}. Since this is an example not
related to any particular application, $c_0^{(2)}$ has been directly taken
in order to have $t_0=0$ (that is,
$\tau^{2,1}_{h_1,4}(c_0^{(2)})\in\toralfa1{h_1}{\alpha_1^{(1)}}$).
The integration times of the remaining center manifold arcs are
$t_2=12.361$, $t_4=3.199$,
$t_6=3.168$, that amount to $18.728$ dimensionless
units, corresponding to 2.98 years. The top row of
Figure~\ref{fig:acrob_h1} represents, in the $\poincpi$ Poincaré sections of
$\Wcal^c_{h_1}(L_2)$ and $\Wcal^c_{h_1}(L_1)$, the invariant curves of the
transition tori and, for each one, the fraction of invariant curve covered by
the corresponding center manifold arc (either from $P_1$ to $P_2$ or from $Q_1$
to $Q_2$). The top row of Table~\ref{tab:trajs} summarizes data on this
pseudo-connection. Overall, in $47.584$ dimensionless time units (7.58 years),
$\alpha$ is decreased from
$0.980$ to $0.898$. This corresponds to an increase in vertical amplitude that
can be appreciated in the lower-right corner of Figure~\ref{fig:acrob_h1}.
Due to its temporal length, this pseudo-trajectory is not very practical from
an astrodynamical point of view, but is convenient to illustrate the
construction of pseudo-trajectories with several arcs and also to see the kind
of vertical amplitude changes that can be obtained through heteroclinic
connections in the Sun-Earth case.

\begin{table}
\centering
\renewcommand{\arraystretch}{1.3}
\begin{tabular}{cccccccc}
   Case & $l_0$ & $l_n$ & $\alpha_0^{(l_0)} \rightarrow \alpha_n^{(l_n)}$
        & Chain (maps) & $\sum\limits_{i:k_i=0} t_i$
                       & $\sum\limits_{i:k_i>0} t_i$ & Figure \\
   \hline
   1 & 2 & 2 & $0.980 \rightarrow 0.898$
     & $\tau_{h_1,4}^{2,1}, \tau_{h_1,4}^{1,2},
        \tau_{h_1,4}^{2,1}, \tau_{h_1,4}^{1,2}$
     & $18.728$ & $28.85$ & \ref{fig:acrob_h1} \\
   2 & 2 & 2 & $0.133 \rightarrow 0.992$
     & $\tau_{h_2,6}^{2,1}, \tau_{h_2,6}^{1,2}$
     & $50.26$ & $18.113$ & \ref{fig:acrob_h2} \\
   2 & 1 & 2 & $0.822 \rightarrow 0.432$
     & $\tau_{h_2,3}^{1,1}, \tau_{h_2,6}^{1,2}$
     & $91.286$ & $14.827$ & \ref{fig:acrob_homos} \\
   3 & 2 & 2 & $0.971 \rightarrow 1.010$
     & $\tau_{h_3,2}^{2,1}, \tau_{h_3,2}^{1,2}$
     & $5.819$ & $6.8800$ & \ref{fig:acrob_h3}
\end{tabular}
\caption{Summary data on the pseudo-trajectories computed in this section:
   initial and final torus parameters, chains of connectivity maps,
   accumulated time of center manifold and connection arcs, and corresponding
   figures.}
\label{tab:trajs}
\end{table}

\begin{figure}[tb]
   \centering
   \includegraphics{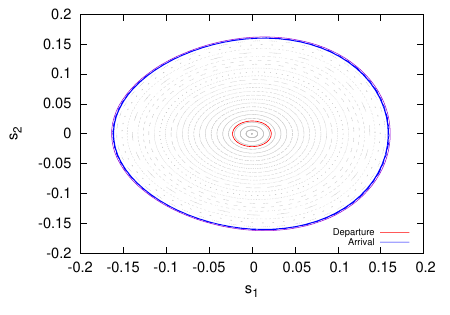}
   \includegraphics{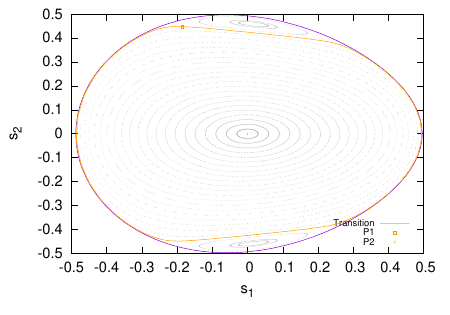}
   \includegraphics{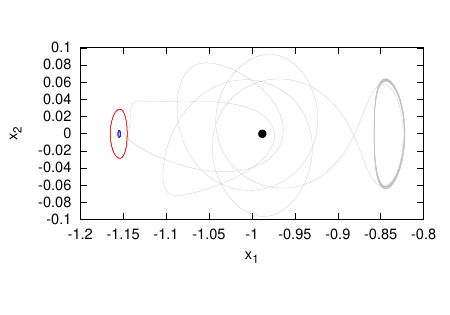}
   \includegraphics{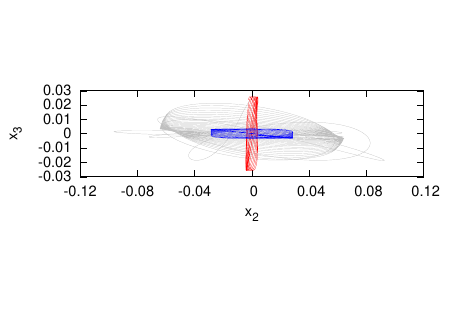}
   \caption{
      Pseudo-trajectory constructed using the primary connectivity map
      $F_{h_2,6}^{2,1}$ shown in
      Figure~\ref{fig:kvk_tl_-1.585}. Two heteroclinic transitions connect
      $\alpha_0^{(2)}=0.133$ to $\alpha_2^{(2)}=0.992$ through the intermediate
      torus $\alpha_1^{(1)}=0.999$.
      Points $P_1$ and $P_2$ in the $\poincpi$ section
      indicate the integration interval
      within the center manifold.
   }
   \label{fig:acrob_h2}
\end{figure}

Our second example is for the Earth-Moon case, namely for case 2 in
Table~\ref{tab:scen}.  The pseudo-trajectory illustrated in
Figure~\ref{fig:acrob_h2}, that uses the connectivity map $F_{h_2,6}^{2,1}$
represented in Figure~\ref{fig:kvk_tl_-1.585}, achieves a substantially larger
variation in $\alpha^{(2)}$: from $\alpha_0^{(2)}=0.133$ to
$\alpha_2^{(2)}=0.992$, using only two heteroclinic connections.
There is a single intermediate torus on $\Wcal_{h_2}^c(L_1)$ with parameter
$\alpha_1^{(1)}=0.999$. The total duration of the trajectory is $68.373$
dimensionless units, corresponding to approximately $297$ days.

\begin{figure}[tb]
   \centering
   \includegraphics{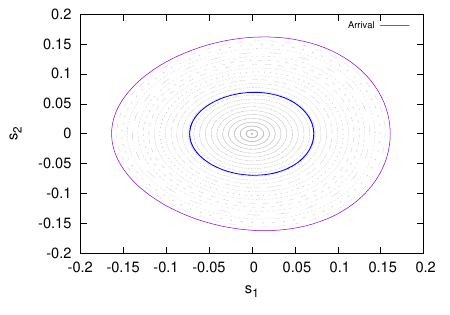}
   \includegraphics{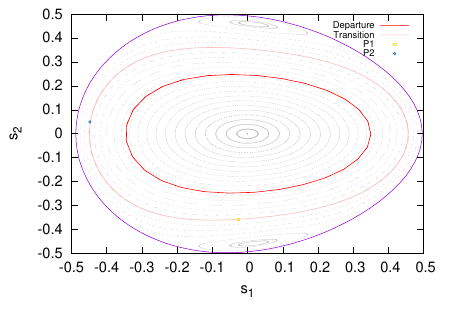}
   \includegraphics{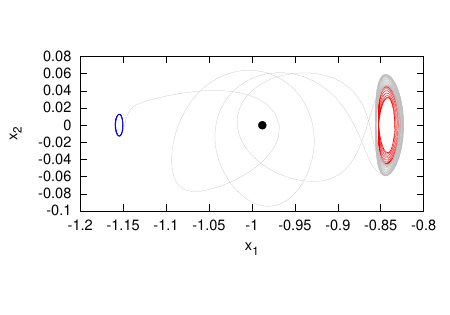}
   \includegraphics{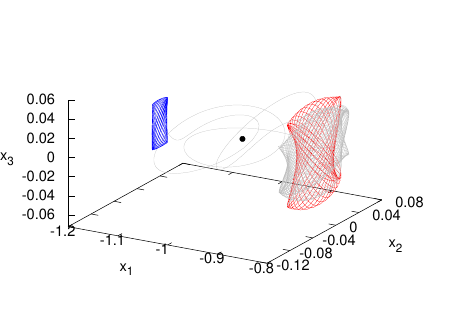}
   \caption{
      Pseudo-trajectory combining one homoclinic and one heteroclinic transition
      to transfer from $\alpha_0^{(1)}=0.822$ on $\Wcal_{h_2}^c(L_1)$ to
      $\alpha_2^{(2)}=0.432$ on $\Wcal_{h_2}^c(L_2)$, using the intermediate
      torus $\alpha_1^{(1)}=0.930$. Points $P_1$ and $P_2$
      indicate the center manifold integration intervals in the $\poincpi$ section.
   }
   \label{fig:acrob_homos}
\end{figure}

As discussed in Section~\ref{ssec:conmaps}, the connectivity map associated
with homoclinic connections in $\Wcal_{h_2}^c(L_1)$, shown in
Figure~\ref{fig:homos} left, can be used to extend the
set of tori around $L_1$ that are heteroclinically connected to tori around
$L_2$. According to Figure~\ref{fig:kvk_tl_-1.585}, the lowest
$\alpha^{(1)}$ value of tori connected to $L_2$ is, roughly, $0.92$.
The pseudo-trajectory displayed in Figure~\ref{fig:acrob_homos}, however,
starts from a torus on $\Wcal_{h_2}^c(L_1)$ with parameter
$\alpha_0^{(1)}=0.822$ and reaches a torus on $\Wcal_{h_2}^c(L_2)$ with
parameter $\alpha_2^{(2)}=0.432$. The trajectory is constructed by first
applying the homoclinic transition map $\tau_{h_2,3}^{1,1}$,
in order to target an intermediate torus in $\Wcal_{h_2}^c(L_1)$ with
parameter $\alpha_1^{(1)}=0.930$. This is followed by a heteroclinic excursion
via the transition map $\tau_{h_2,6}^{1,2}$, which produces the final
connection to $\alpha_2^{(2)}$. The total duration of this pseudo-trajectory
is $106.113$ dimensionless units, corresponding to approximately 462 days.

\begin{figure}[tb]
   \centering
   \includegraphics{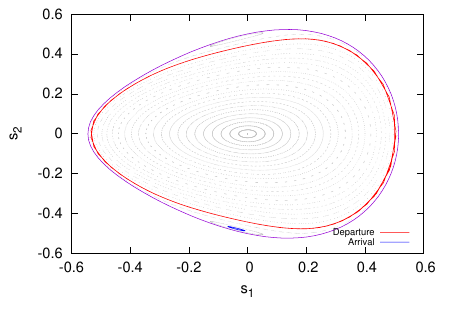}
   \includegraphics{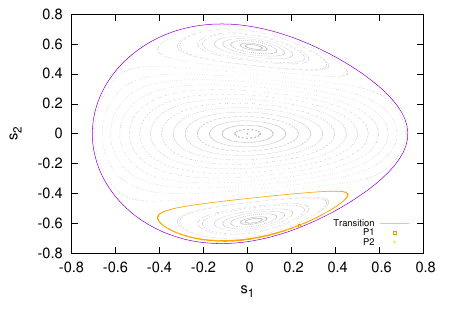}
   \includegraphics{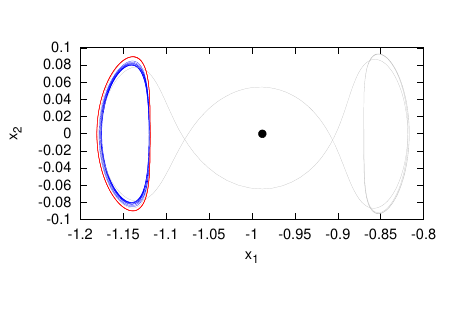}
   \includegraphics{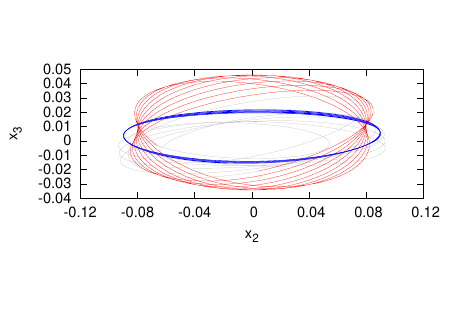}
   \caption{
      Pseudo-trajectory constructed using the primary connectivity map
      $F_{h_3,2}^{2,1}$ shown in Figure~\ref{fig:kvk_tl_-1.575}. Two
      heteroclinic transitions connect $\alpha_0^{(2)}=0.971$ to
      $\alpha_2^{(2)}=1.010$ through the intermediate torus
      $\alpha_1^{(1)}=1.020$.  Points $P_1$ and $P_2$ in $\poincpi$ section
      indicate the integration intervals within the center manifold.
}
   \label{fig:acrob_h3}
\end{figure}

In addition to changes in vertical amplitude, the connectivity maps also
enable transitions between different families of invariant tori within the
center manifold. In particular, the connectivity maps shown in
Figures~\ref{fig:kvk_tl_-1.585} and \ref{fig:kvk_tl_-1.575}
enable transfers between Lissajous and quasi-Halo tori on $\Wcal_{h}(L_j)$
through two heteroclinic excursions.
As an example, we consider the connectivity map
$F^{2,1}_{h_3,2}$ displayed in Figure~\ref{fig:kvk_tl_-1.575}.
The pseudo-trajectory illustrated in Figure~\ref{fig:acrob_h3} connects a
Lissajous torus with parameter $\alpha_0^{(2)}=0.971$ to a quasi-Halo torus
with parameter $\alpha_2^{(2)}=1.010$. The transition is achieved by applying
the transition map $\tau_{h_3,2}^{2,1}$ and its reverse through an
intermediate torus on $\Wcal_{h_3}^c(L_1)$ with parameter
$\alpha_1^{(1)}=1.020$ that allows this transfer.  The total duration of this
trajectory is 12.699 dimensionless units, corresponding to approximately 55
days.

\section{Connectivity
at the level of individual connections}
\label{sec:torus2torus}

In the previous section, we computed connectivity maps
that determine what invariant tori are accessible from a given one through
heteroclinic or homoclinic connections. While these maps reveal whether a pair
of tori are connected, they do not provide information on the number of
individual
connections, their precise locations, or the flight times along the
pseudo-trajectories connecting the tori. As illustrated in the first example
of Section~\ref{ssec:traj}, these times can
be very large. 

Let $l,m \in \{1,2\}$, and consider two-dimensional tori
$\Tcal^{(l)}\subset\Wcal_h^c(L_l)$, $\Tcal^{(m)}\subset\Wcal_h^c(L_m)$. The
$k$-heteroclinic (resp.~$k$-homoclinic) connections between $\Tcal^{(l)}$ and
$\Wcal^c(L_m)$ are
given by the intersection of $\Tcal^{(l)}$ with the departure CS, that we
define as
\[
   \Dhk{l}{m} ( \Tcal^{(l)} ) := \Tcal^{(l)} \cap \Dhk{l}{m}.
\]
On the other hand, the $k$-heteroclinic (resp.~$k$-homoclinic) connections
between $\Tcal^{(m)}$
and $\Wcal^c(L_l)$ are given by the intersection of $\Tcal^{(m)}$ with the
arrival CS, that we define as
\[
   \Ahk{l}{m} ( \Tcal^{(m)} ) := \Tcal^{(m)} \cap \Ahk{l}{m}.
\]
We can
send
these points back to $\Dhk lm$ through the inverse of the
transition map $\tau^{l,m}_{h,k}$, and define
\[
   \tDhk lm(\Tcal^{(m)}):=
      (\tau^{l,m}_{h,k})^{-1}\bigl( \Ahk lm(\Tcal^{(m)}) \bigr).
\]
In this way, the $k$-heteroclinic (resp.~$k$-homoclinic) connections between
$\Tcal^{(l)}$ and
$\Tcal^{(m)}$ can be visualized on $\Dhk lm$ as the intersection
\begin{equation}\label{eq:intcurves}
   \Dhk lm(\Tcal^{(l)}) \cap \tDhk lm(\Tcal^{(m)}).
\end{equation}
Analogously, the departure intersection set $\Dhk{l}{m} (
\Tcal^{(l)} )$
can be sent forward to the arrival CS $\Ahk{l}{m}$ through the 
transition map,
\[
   \tAhk lm(\Tcal^{(l)}):=
   \tau^{l,m}_{h,k}\bigl( \Dhk lm(\Tcal^{(l)}) \bigr),
\]
so that \eqref{eq:intcurves} is equivalently represented on the arrival CS as
\begin{equation}\label{eq:intcurvesarr}
   \Ahk lm(\Tcal^{(m)}) \cap \tAhk lm(\Tcal^{(l)}).
\end{equation}

The goal of this section is to visualize $k$-heteroclinic connections between
pairs of tori by using intersections between tori and connection surfaces in
the way just described.

A straightforward approach to obtain points on the curve $\Curveju$ using last
Section's strategy would be: find the $\alpha^{(l)}$ value corresponding to
$\Tcal^{(l)}$, shoot from equally spaced points in the
$\{c(\theta,r_{\alpha^{(l)}}(\theta))\}_{\theta\in\T}$ curves (see
Eq.~\eqref{eq:polars}), and flow them through $\varphi^{(l)}_t$ in order to
intersect $\Dhk lm$.  The computation of these intersections turns out to be
delicate, since the flow is not always transverse to $\Dhk lm$: there are
tangencies that correspond to non-transverse heteroclinic connections (see
\cite{2008aDeMaRo}) and they give rise to numerical difficulties when computing
intersections with $\Dhk{l}{m}$ nearby.  As a result, gaps are obtained when
plotting the curves in $\Curveju$.

A way to overcome this difficulty is to perform continuation on
Eqs.~\eqref{eq:conncond} but with $c^{(l)}$ constrained to be in
$\Tcal^{(l)}$. For that, a parameterization
$\theta\in\T^2\mapsto\Psi^{(l)}(\theta)\in\calV^{(l)}$ of $\Tcal^{(l)}$ is
needed.  With it available, obtaining $\Dhk lm(\Tcal^{(l)})$ is a matter of
performing continuation on
\begin{equation}\label{eq:conttorWc}
   \Phi_{T^{(l)}}\bigl(W^{(l)}(\Psi^{(l)}(\theta),\xi^{(l)})\bigr)
   =
   \Phi_{T^{(m)}}\bigl(W^{(m)}(c^{(m)},\xi^{(m)})\bigr)
\end{equation}
with respect to $\theta,T^{(l)},c^{(m)},T^{(m)}$. 
When doing this, we will still go through tangency
points that will produce rank deficiency in the system of equations
\eqref{eq:conttorWc}, but a standard predictor--corrector method (see
e.g.~\cite{allgower-georg}) has no difficulty ``jumping over'' these critical
points as long as they are not hit exactly, which is almost impossible
numerically. The specific details on the actual equations used for
continuation are given in Appendix~\ref{app:int2tori}.
The $\Psi^{(l)}$ parameterization can be computed either with a
large-matrix method
(e.g.~\cite{2001GoMo,2018BaOlSche}) or with a
parameterization-based one (e.g.~\cite{2021aHaMo,2022KuAndLla}). But
it is simpler to obtain a Fourier series of a parameterization of it (in
$\calV^{(j)}$ coordinates) from direct numerical integration of the restricted
flow to the center manifold
$\varphi_t^{(l)}$ via a frequency analysis technique
(e.g.~\cite{1999La,2010GoMoSi,2010GoMoSi-2,2020SaMe,2021MeSa}). 
This is possible thanks to the fact that we can perform long-time integrations
in the center manifold using the reduced dynamics (see
Section~\ref{ssec:dynwc}). 
We have used the
methodology in \cite{2010GoMoSi-2}. In this
way,
a parameterization
$\Psi^{(l)}$ and a frequency vector $\omega^{(l)}\in\R^2$ are found
satisfying the invariance equation
\begin{equation}
   \varphi_t^{(l)}\bigl(
      \Psi^{(l)}(\theta)
   \bigr)
   =
      \Psi^{(l)}(\theta+t\omega^{(l)})
   ,
   \quad\forall\theta\in\T^2,
   \label{eq:inv2tori}
\end{equation}
with a tolerance of $10^{-10}$ estimated over a uniform mesh in $\T^2$.
The specific details on how a Fourier series is obtained from this methodology
are given in Appendix~\ref{app:int2tori}. 

Initial conditions for the continuation of
equations \eqref{eq:conttorWc} are found by shooting from the
$\{c(\theta,r_{\alpha^{(l)})}(\theta)\}_{\theta\in\T}$ curves as mentioned
before, but this needs to be done only once for every connected component to
be found in $\Dhk lm(\Tcal^{(l)})$.

Performing numerical continuation on equation~\eqref{eq:conttorWc}
produces a trace in the $\theta$-parameter space that describes the curves in
$\Curvejuk$ in torus coordinates. This trace is useful in predicting the time
of flight necessary within $\Tcal^{(l)}$ between heteroclinic connections, as
will be illustrated in Section~\ref{sec:traces}. On the other hand, the
topology of this trace is related to the topology of $\Dhk lm(\Tcal^{(l)})$
within $\Dhk lm$. This will be useful in discussing the visualization of the
$\Dhk lm(\Tcal^{(l)})$, $\tDhk lm(\Tcal^{(m)})$ sets in
Section~\ref{sec:tor2tor}.

\subsection{Traces of torus-CS intersections in torus parameter space}
\label{sec:traces}
In this section, 
we present numerical results 
for the computation of the curves $\Curveju$ and $\Curveja$, obtained by
numerical continuation of \eqref{eq:conttorWc}.
An absolute tolerance of
$10^{-10}$ is used for the Newton method, and
the invariance equation \eqref{eq:inv2tori} is enforced with a maximum error of
$10^{-10}$ for the tori parameterizations 
(see
Appendix~\ref{app:int2tori}).
Specifically, we will discuss the different topological configurations
that the traces of $\Curveju$ and $\Curveja$ can exhibit in the
$\theta$-parameter space and the connectivity information that can be obtained
in each case.

\subsubsection{Two non-contractible curves}

\begin{figure}[tbp]
   \centering
   \includegraphics{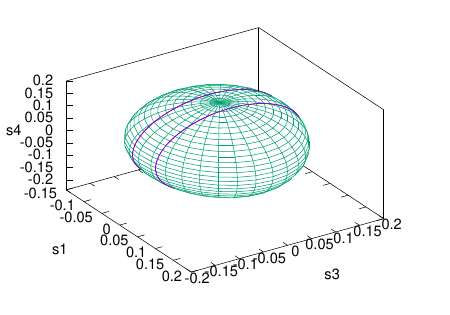}
   \includegraphics{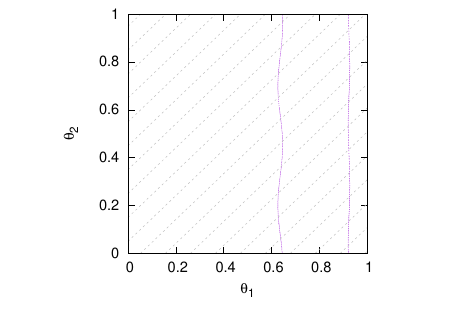}
   \caption{
      Heteroclinic connections between the Lissajous torus
      $\Tcal^{(2)}:=\toralfa 2{h_2}{0.24}$ and the center manifold
      $\Wcal_{h_2}^c(L_1)$.  Left: $\Curvejh{2}{1}{h_2}{6}$
      in purple, with $\CalD{2}{1}{h_2}{6}$ in green. Right: corresponding
      trace of $\Curvejh{2}{1}{h_2}{6}$
      in $\theta$-parameter space. Gray lines indicate the direction of the
      flow on the torus.
   }
   \label{fig:tordinsesfera1}
\end{figure}

The first scenario occurs when the $\theta$-trace of $\Curveju$ or $\Curveja$
consists of two non-contractible curves, as shown in
Figure~\ref{fig:tordinsesfera1} right. In this case, the intersection of the
torus with the connection surface yields two disjoint curves, each homeomorphic
to $S^1$, on the corresponding $\Dhk lm$ or $\Ahk lm$, as illustrated in
Fig.~\ref{fig:tordinsesfera1} left. The torus therefore intersects the
connection surface twice before returning to any section of the form $\{\theta_i
= a \}$, with $a\in\T$. From an applications standpoint, this configuration is
particularly favorable, since, for any of its frequencies, the torus
intersects the CS
at any of the corresponding periods, thus providing multiple opportunities for
heteroclinic connections.

In this scenario, it is possible to estimate the maximum and minimum time
between consecutive intersections of the torus with the CS.  Consider the flow
on a torus with frequency vector $\omega=(\omega_1,\omega_2)$, given by
\eqref{eq:inv2tori}, which is foliated by straight lines of irrational slope
$\omega_2/\omega_1$.
Suppose that the torus, parameterized by $\Psi(\theta)$, intersects the CS
along two non-contractible curves that can be approximated by vertical lines
located at $\theta_1=\theta_a$ and $\theta_1=\theta_b$, with $\theta_b >
\theta_a$.
Assume for simplicity
that $\theta_b-\theta_a < 1 - (\theta_b -
\theta_a)$.
The minimum return time corresponds to the flow travelling from the first intersection
curve at $\theta_a$ to the other at $\theta_b$ in the $\theta_1$-direction, and
is given by
\begin{equation}\label{eq:mintime}
   t_{\min}=\frac{\theta_b-\theta_a}{\omega_1}.
\end{equation}
In contrast, the maximum return time corresponds to
the time required to return to the first intersection curve at $\theta_a$ in the
$\theta_1$-direction after the intersection with $\theta_b$, namely
\[
   t_{\max}=\frac{\theta_a+1-\theta_b}{\omega_1}=\frac{1}{\omega_1}-t_{\min}.
\]
For the case illustrated in Figure~\ref{fig:tordinsesfera1}, we numerically
obtain $\omega_1$ and $\omega_2$ as shown in Table~\ref{tab:freqstors}, with
$\theta_a \approx 0.63$
and $\theta_b \approx 0.92$.  Consequently, the estimated return times are
$t_{\min} = 0.98$ and $t_{\max} = 2.39$ in
nondimensional RTBP units, corresponding roughly to 4 and 10 days, respectively.

\begin{table}[htbp]
\begin{tabular}{c|c|c|c|c|c}
Fig. & $\alpha^{(2)}$ & $h$ & $\omega_1$ & $\omega_2$ &
      $\rho_2(\omega)$ \\
\hline
\ref{fig:tordinsesfera1} &
$0.24$ & $-1.585$ & $0.29670084$ & $0.28406030$ & $-0.04260$ \\
\ref{fig:torpardinsesfera1} &
$0.14$ & $-1.585$ & $0.29672486$ & $0.28403990$ & $-0.04275$ \\
\ref{fig:torqhpardinsesfera1} &
$1.001$ & $-1.575$ & $0.29287653$ & $0.28860622$ & $-0.01458$ \\

\end{tabular}
\caption{\label{tab:freqstors}Parameter values corresponding to
   Figs.~\ref{fig:tordinsesfera1}, \ref{fig:torpardinsesfera1}, and
   \ref{fig:torqhpardinsesfera1}. Each case shows the heteroclinic connections
   between a torus $\Tcal^{(2)}\subset\Wcal_h^{c}(L_2)$ and the center manifold
   $\Wcal^c(L_{1})$.  The torus $\Tcal^{(2)}$ has $\alpha=\alpha^{(2)}$ and
   frequencies $\omega_1,\omega_2$.}
\end{table}

\begin{figure}[tbp]
   \centering
   \includegraphics{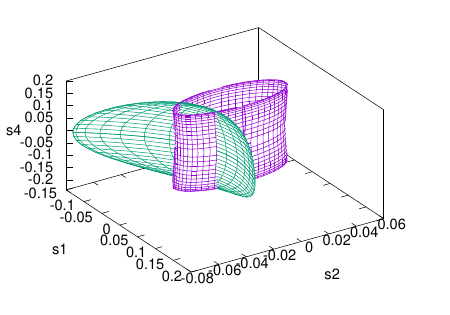}
   \includegraphics{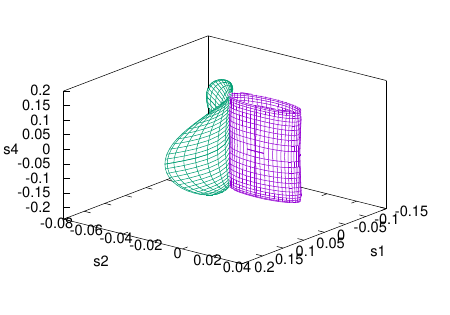}
   \caption{Projection of $\CalD{2}{1}{h_2}{6}$ in green together with the tori
   $\Tcal^{(2)} \in \Wcal_{h_2}(L_2)$ with $\alpha^{(2)}=0.240$ (left) and
   $\alpha^{(2)}=0.140$ (right) in purple.}
   \label{fig:toresf}
\end{figure}

Our numerical experiments indicate that this is the case for tori when they
are close to their backbone p.o.~and this p.o.~intersects transversally the CS
in two points. Nearby tori intersect the CS in two curves that surround these
points, one around each. Trajectories flow inside the tori in a way that they
intersect each of these curves once every period associated to the frequency
close to the one of the backbone p.o., and this is why their $\theta$-traces
show as two non-contractible curves, as in Figure~\ref{fig:tordinsesfera1}
right. In the case of this Figure, the CS is $\calD^{2,1}_{h_2,6}$ and the
backbone orbit is the Lyapunov p.o.~which intersects the CS in two points as
shown in Fig.~\ref{sfig:esfera}. We will see in Section~\ref{sec:tor2tor}
that this also happens for the Halo family, but not always (see
Figure~\ref{fig:esferacas3}). Actually, the torus $\Tcal^{(2)}:=\toralfa
2{h_2}{0.24}$ of Figure~\ref{fig:tordinsesfera1} is not close to the planar
p.o., but to the last torus (lowest $\alpha^{(2)}$ value) for which the
$\theta$-trace is still seen as two non-contractible curves, see
Figure~\ref{fig:toresf} left.
As the torus
grows in size, a change in the topology of the $\theta$-trace takes place that
leads to the scenario that is described next. Further details will be given in
the discussion of Figures~\ref{fig:esferacas1} and \ref{fig:esferacas2} in
next Section.

\subsubsection{Two contractible curves}

\begin{figure}[tb]
   \centering
   \includegraphics{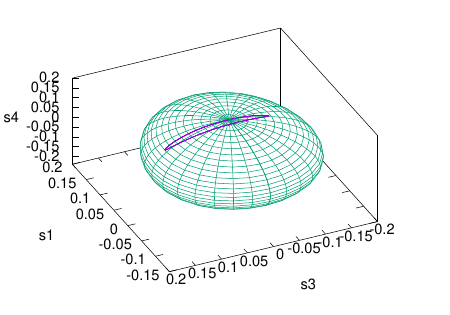}
   \includegraphics{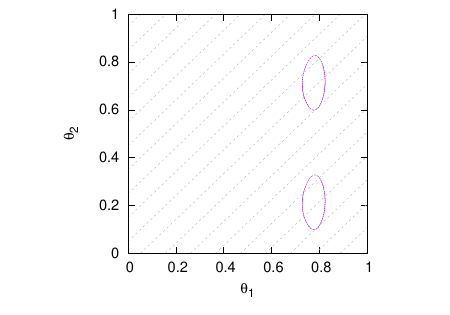}
   \caption{
      Heteroclinic connections between the unstable manifold of a Lissajous
      torus $\Tcal^{(2)} \in \Wcal_{h_2}(L_2)$ with $\alpha^{(2)}=0.14$ and the
      center-stable manifold $\Wcal_{h_2}^{cs}(L_1)$.  
      Left: $\Curvejh{2}{1}{h_2}{6}$, with $\CalD{2}{1}{h_2}{6}$ in green and
      the torus $\Tcal^{(2)}$ in purple. A symmetric curve exists in
      the opposite hemisphere.
      Right: corresponding trace of $\Curvejh{2}{1}{h_2}{6}$ in the
      $\theta$-parameter space. Gray lines indicate the direction of the flow
      on the torus.
   } 
   \label{fig:torpardinsesfera1}
\end{figure}

The second scenario occurs when the trace of $\Curveju$ (resp.~$\Curveja$) in
$\theta$-parameter space consists of two contractible curves, which in turn
give rise to two curves on the departure CS $\Dhk{l}{m}$ (resp.~the arrival CS
$\Ahk{l}{m}$).  An example of this behavior is shown in
Figure~\ref{fig:torpardinsesfera1}, where $\Curveju$ is computed for a
particular Lissajous torus of
$\Wcal_{h_2}^c(L_2)$.

Similarly to the previous scenario, it is possible to estimate the minimum
and maximum times between consecutive intersections of a torus with the
CS. By considering vertical lines $\theta_1=\theta_a$ and $\theta_1=\theta_b$
tangent to the curves of the $\theta$-traces
(Figure~\ref{fig:torpardinsesfera1} right), an upper bound for the minimum
time (understood as the time between two consecutive intersections of a purple
curve in Figure~\ref{fig:torpardinsesfera1}) is given again by
\eqref{eq:mintime}.
Estimating the maximum time between consecutive intersections of a trajectory
with the CS requires additional care. For the sake of discussion, let us
consider the vertical line $\theta_1=0.8$ in the right
plot of Figure~\ref{fig:torpardinsesfera1}. Of all the trajectories of the
torus that go across this line, there are two intervals of values of
$\theta_2$, let us call them $[a_1,a_2]$ and $[b_1,b_2]$ with $a_2<b_1$, for
which the corresponding trajectory intersects one of the closed purple curves
of the $\theta$-trace.
The endpoints of these intervals $a_1$, $a_2$, $b_1$, $b_2$ correspond to
trajectories that are tangent to the CS. These tangencies are responsible for
the need for the continuation strategy of Eq.~\eqref{eq:conttorWc}, due to the
numerical difficulties in shooting directly from the $\poincpi$ sections of
the tori, as it has been mentioned before.
Let $\Delta \theta_2$ be
either $b_1-a_2$ or $a_1+1-b_2$.
Every $1/\omega_1$ units of time, the net increment in the
$\theta_2$-coordinate is given by 
\[
\rho_2(\omega):=\frac{\omega_2}{\omega_1}
   -\round\Bigl(\frac{\omega_2}{\omega_1}\Bigr) , 
\] 
($\round(x)=[x+1/2]$), which represents the effective vertical drift per return
in $\theta_1$. Its actual value can be found in Table~\ref{tab:freqstors}.  A
trajectory starting from $\theta_2=a_2$ or $\theta_2=b_2$ will intersect the
other curve once the accumulated vertical drift reaches $\Delta \theta_2$.
Understanding this as the maximum time between intersections with
$\theta$-traces, it can be estimated as
\[
   t_{\max} = \min_{k \in \N} \biggl\{
      \Bigl|\frac{k}{\omega_1}\Bigr|
            : k|\rho_2(\omega)| \geq \Delta \theta_2
   \biggr\}
   =
   \frac1{|\omega_1|}\Big\lceil
      \frac{\Delta\theta_2}{|\rho_2(\omega)|}
   \Big\rceil .
\]

For the example shown in Figure~\ref{fig:torpardinsesfera1},
the minimum return time is $t_{\min}=0.3249$,
corresponding to an estimated distance
$\theta_b-\theta_a\approx0.0964$ and the value of $\omega_1$ given in
Table~\ref{tab:freqstors}. The maximum return time can be estimated by taking
$\Delta \theta_2 = 0.2710$, which, in this case, is equal to both
$b_1-a_2$ and $a_1+1-b_2$. This yields $k=7$ and then
$t_{\max}=23.591$. These minimum and maximum times
correspond to approximately 1.41 and 102.6 days, respectively.
Large disconnection gaps between
consecutive intersections represent a major limitation when designing
trajectories of the type discussed in Section~\ref{ssec:traj}, as they are the
primary source of long transition times.  Accurate modeling of these
trajectories therefore requires identifying these regions and avoiding them
whenever possible, for instance by exploiting alternative connection manifolds
such as those introduced in Section~\ref{ssec:some}.

We have numerically observed that this scenario arises when the torus
approaches tangency with the connection surface $\Dhk{l}{m}$ or $\Ahk{l}{m}$,
see Figure~\ref{fig:toresf} right.
In the example just discussed, as $\alpha^{(2)}$ decreases,
the near-tangency evolves into a heteroclinic tangency between a torus and the
opposite center manifold, ultimately leading to the extinction of the
corresponding connection manifold. 
All the families of Lissajous tori involved in heteroclinic connections that we
have explored exhibit this same behavior. The details can be found in
Section~\ref{sec:tor2tor}.

\subsubsection{Single contractible curve}

\begin{figure}[tb]
   \centering
   \includegraphics{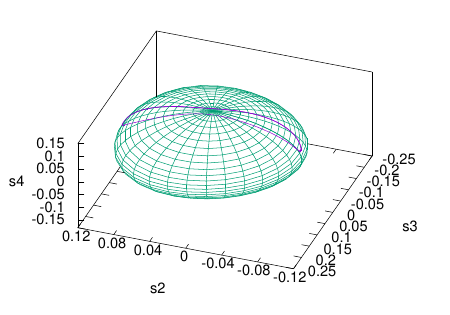}
   \includegraphics{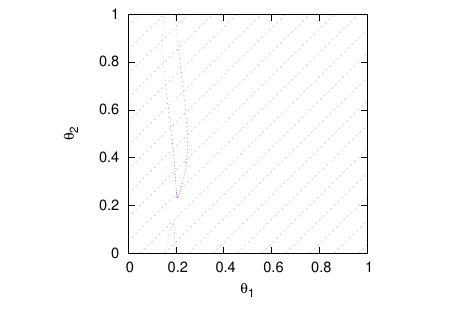}
   \caption{
      Heteroclinic connections between the unstable manifold of a Northern
      quasi-Halo
      torus on $\Wcal_{h_3}(L_2)$ with $\alpha^{(2)}=1.001$ and the
      center-stable manifold $\Wcal_{h_3}^{cs}(L_2)$. 
      Left: $\Curvejh{2}{1}{h_3}{2}$, with $\CalD{2}{1}{h_3}{2}$ in green and
      the torus $\Tcal^{(2)}$ in purple. 
      Right: corresponding trace of $\Curvejh{2}{1}{h_3}{2}$ in the
      $\theta$-parameter space. Gray lines indicate the direction of the flow
      on the torus.
   }
   \label{fig:torqhpardinsesfera1}
\end{figure}

The last scenario occurs when the torus intersects the CS in a single
connected
component, and the corresponding trace of
$\Curveju$ or $\Curveja$
in $\theta$-parameter space consists of a single
contractible curve. An example of this configuration is shown in
Figure~\ref{fig:torqhpardinsesfera1}. In our
computations, this behavior was observed exclusively for quasi-Halo tori. 
The approximation for
times of disconnection
between the torus and the CS
introduced earlier also applies in this case.
From Figure~\ref{fig:torqhpardinsesfera1} and the values of
Table~\ref{tab:freqstors}, 
the minimum return time to the CS is $t_{\min}=0.226$, whereas the maximum
return time is $t_{\max}=27.315$, corresponding to approximately 1 and 119
days, respectively.

\subsection{Visualization of connections between pairs of invariant tori}
\label{sec:tor2tor}

In the previous section, we explored the different scenarios in which
$\Curveju$ and $\Curveja$ arise, illustrating how intersections between a
torus and the CS can present different topologies.  In this section, we turn
our attention to connections between specific pairs of tori, using
\eqref{eq:intcurves} to systematically locate
them.

\subsubsection{Sun-Earth, Lissajous--Lissajous}
\begin{figure}[tb]
   \centering
   \includegraphics[trim=1.5cm 1.5cm 1.5cm 1.5cm,clip]{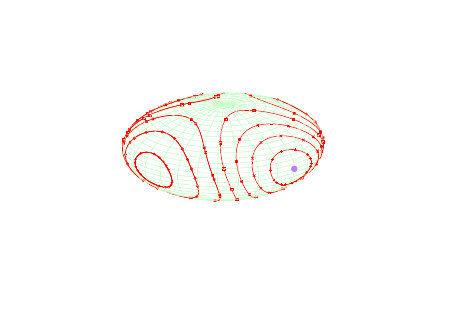}
   \includegraphics[trim=1.5cm 1.5cm 1.5cm 1.5cm,clip]{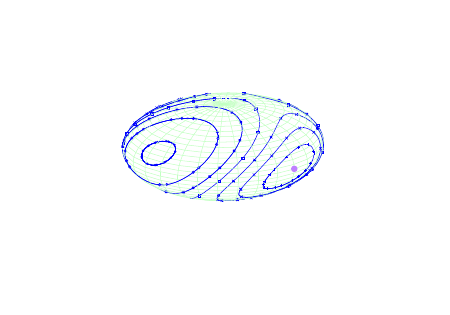}
   \includegraphics[trim=1.5cm 1.5cm 1.5cm 1.5cm,clip]{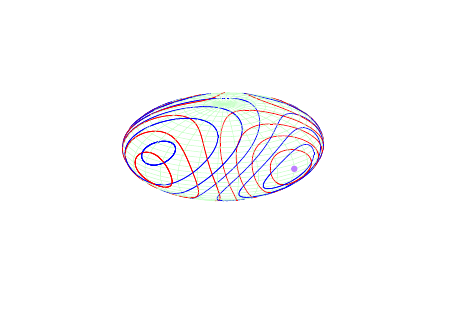}
   \caption{
      For Case 1 in Table~\ref{tab:scen}, 
      $(s_2,s_3,s_4)$ projections of, from left to right,
      $\{\Curvejhk{2}{1}{h_1}{4}\}_i$, $\{\tCurvemjhk{2}{1}{h_1}{4}\}_j$ and
      $\{\Curvejhk{2}{1}{h_1}{4}\}_i \cup \{\tCurvemjhk{2}{1}{h_1}{4}\}_{j}$,
      all represented on the $\CalD{2}{1}{h_1}{4}$ CS, shown in green and
      rendered as opaque. Red curves correspond to Lissajous tori from $L_2$,
      while blue curves correspond to Lissajous tori from $L_1$.  Different
      marker styles correspond to different tori. The heteroclinic connections
      between the Planar Lyapunov p.o.~are shown in purple.
   }
   \label{fig:esferacas1}
\end{figure}

Figure~\ref{fig:esferacas1} presents the results for Case 1 in
Table~\ref{tab:scen}. The left panel displays $\{\Curvejhk{2}{1}{h_1}{4}\}_i$
for a representative set of tori $\Tcal^{(2)}_i\subset\Wcal^c_{h_1}(L_2)$.
Each torus intersects the CS in two curves, that are displayed in
Figure~\ref{fig:esferacas1}, with the same markers. At first, pairs of circles
wrap around the two planar heteroclinic connections, represented as two purple
points (one is hidden behind the CS). They correspond to tori in $L_2$ near
the planar Lyapunov orbit (i.e., $\alpha^{(2)}$ close to 1), producing two
non-contractible curves in $\theta$-parameters, similar to the ones in
Figure~\ref{fig:tordinsesfera1}. For values of $\alpha^{(2)} \approx 0.94$,
the pairs of curves change orientation, resulting in intersections seen as two
contractible curves in $\theta$-parameters, similar to the ones in Figure
\ref{fig:torpardinsesfera1}.
A similar pattern appears in the center panel, that shows the translated
intersections $\{\tCurvemjhk{2}{1}{h_1}{4}\}_j$, for a representative set of
tori $\Tcal^{(1)}_j\subset\Wcal^c_{h_1}(L_1)$. As before, the change in
orientation coincides with a change in the topology of the $\theta$-trace of
the curves $\tCurvemjhk{2}{1}{h_1}{4}$, occurring around $\alpha^{(1)} \approx
0.95$.  
Finally, the right panel shows the superposition of the two previous panels.
Each intersection between a blue and a red curve corresponds to a connection
between a specific pair of tori $\Tcal^{(2)}_i$, $\Tcal^{(1)}_j$.
The numerical computation of these connections can be addressed by a minor
modification of equation \eqref{eq:conncond}. Further details on this
procedure are provided in Appendix~\ref{app:int2tori}. An example that
involves computing individual connections between different pairs of tori in
this setting is presented at the end of this Section, illustrating the
variety of possible outcomes.

\subsubsection{Earth-Moon, Lissajous--Lissajous and quasi-Halo}

\begin{figure}[tb]
   \centering
   \includegraphics[trim=1.5cm 1.5cm 1.5cm 1.5cm,clip]{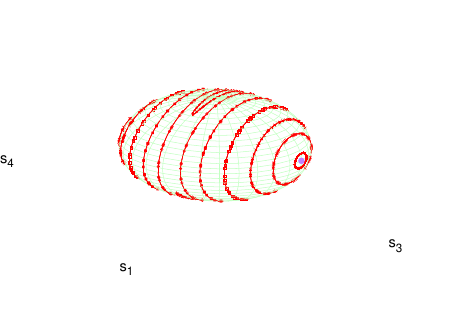}
   \includegraphics[trim=1.5cm 1.5cm 1.5cm 1.5cm,clip]{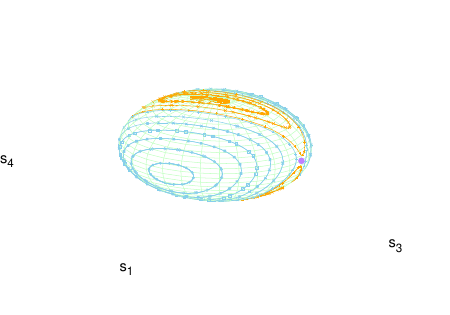}
   \includegraphics[trim=1.5cm 1.5cm 1.5cm 1.5cm,clip]{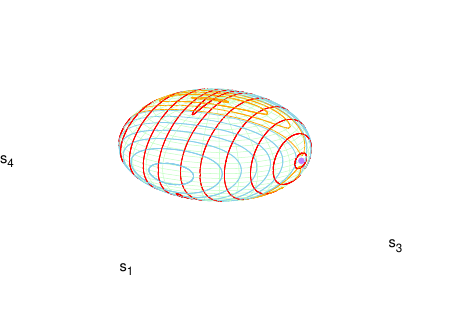}
   \caption{
      For Case 2 in Table~\ref{tab:scen},
      $(s_1,s_3,s_4)$ projections of, from left to right,
      $\{\Curvejhk{2}{1}{h_2}{6}\}_i$, $\{\tCurvemjhk{2}{1}{h_2}{6}\}_j$ and
      $\{\Curvejhk{2}{1}{h_2}{6}\}_i \cup \{\tCurvemjhk{2}{1}{h_2}{6}\}_{j}$,
      all represented on the $\CalD{2}{1}{h_2}{6}$ CS, shown in green and
      rendered as opaque.  Red curves correspond to Lissajous tori from $L_2$,
      while blue and orange curves correspond to Lissajous and quasi-Halo tori
      from $L_1$, respectively. Different marker styles correspond to
      different tori.  The heteroclinic connections between the Planar
      Lyapunov p.o.~are shown in purple.
   }
   \label{fig:esferacas2}
\end{figure}

The symmetry in the qualitative behavior of both libration points breaks down
when moving to the Earth-Moon problem. Figure~\ref{fig:esferacas2} shows the
results for Case 2 of Table~\ref{tab:scen},
in which three families of invariant tori are present around $L_1$.
The left panel displays $\{\Curvejhk{2}{1}{h_2}{6}\}_i$ for selected
$\Tcal_i^{(2)}\subset\Wcal_{h_2}^c(L_2)$.
The tori emanating from the Lyapunov p.o., each one giving a pair of closed
curves surrounding the purple points,
generate non-contractible curves in $\theta$-parameters.
Their associated values $\alpha^{(2)}$ range from 1, near the planar
connection, to $\alpha^{(2)}\approx0.14$, where the $\theta$-trace of the
intersections evolves to two contractible curves as the tori approach
tangency with the CS. This is reflected in a change in the orientation of the intersection
curves on the CS. 

The center panel presents the translated intersections
$\{\tCurvemjhk{2}{1}{h_2}{6}\}_j$ for selected
$\Tcal^{(1)}_j\subset\Wcal_{h_2}^c(L_1)$. In this case, it is necessary to
distinguish between the different families. The blue curves correspond to the
intersections of Lissajous tori. Notably, these curves do not surround the
planar connections (purple points, recall that one is hidden behind the CS),
and each torus intersection traces two contractible curves in
$\theta$-parameters similar to the ones presented in
Figure~\ref{fig:torpardinsesfera1}. The largest curves correspond to Lissajous
tori with $\alpha^{(1)}$ close to 1, and, as $\alpha^{(1)}$ decreases, they
progressively shrink into smaller pairs. They eventually intersect the CS
tangentially at $\alpha^{(1)}\approx0.92$, that corresponds to the last
connected torus. The orange curves are the intersections of the quasi-Halo
tori. Unlike the Lissajous case, each quasi-Halo torus produces a single
contractible curve in $\theta$-parameters, as in
Figure~\ref{fig:torqhpardinsesfera1} right. The upper region of the CS
corresponds to the southern Halo family, while the lower region corresponds to
the northern one. The largest curves again correspond to tori with
$\alpha^{(1)}$ near 1, and they gradually contract into smaller curves as
$\alpha^{(1)}$ approaches $d_h^{(1)}$. They eventually collapse to a point,
that corresponds to a tangential intersection of the $L_1$ Halo periodic orbit
with the arrival CS $\CalA{2}{1}{h_2}{6}$.
The boundary separating the two families is the pre-image through the
transition map $\tau^{2,1}_{h_2,6}$ of the separatrix between Lissajous and
quasi-Halo motion in $\Wcal^c_{h_2}(L_1)$ (see Section~\ref{ssec:dynwc}).

The right panel shows the superposition of both sets of intersections.  It can
be seen that the tori around $L_2$ with $\alpha^{(2)}$ small intersect the
largest range of tori around $L_1$, as was observed in
Figure~\ref{fig:kvk_tl_-1.585}.
By plotting just the curves $\CalD21{h_2}6(\Tcal^{(2)}_{i_0})$ and
$\tCalD21{h_2}6(\Tcal^{(1)}_j)$, for $i_0$ fixed and varying
$j$\footnote{These plots have not been included for shortness}, one can find
that, except for the cases in which
$\CalD21{h_2}6(\Tcal^{(2)}_{i_0})$, $\tCalD21{h_2}6(\Tcal^{(1)}_j)$
are tangent, if $\Tcal^{(2)}_{i_0}$ and $\Tcal^{(1)}_j$ are heteroclinically
connected they have either 8 connections, if $\Tcal^{(1)}_j$ is a Lissajous
torus, or 4, if $\Tcal^{(1)}_j$ is a quasi-Halo one (either Northern or
Southern). This is in agreement with some of the examples presented in
\cite{2023McCaHo}.

\subsubsection{Earth-Moon, Lissajous and quasi-Halo--Lissajous and quasi-Halo}

\begin{figure}[tb]
   \centering
   \includegraphics[trim=1.5cm 1.5cm 1.5cm 1cm,clip]{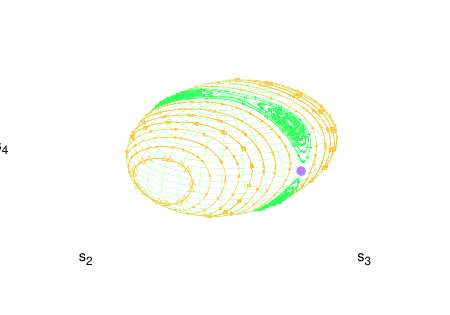}
   \includegraphics[trim=1.5cm 1.5cm 1.5cm 1cm,clip]{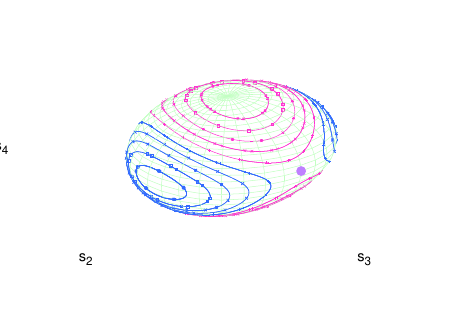}
   \includegraphics[trim=1.5cm 1.5cm 1.5cm 1cm,clip]{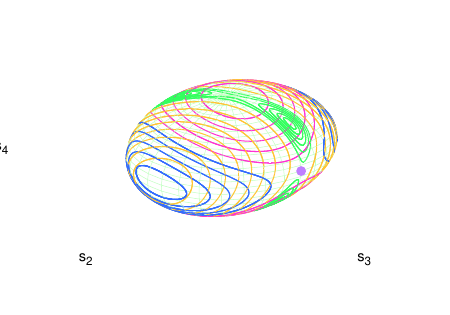}
   \caption{
      For Case 3 in Table~\ref{tab:scen},
      $(s_2,s_3,s_4)$ projections of,
      from left to right,
      $\{\Curvejhk{2}{1}{h_3}{2}\}_i$, $\{\tCurvemjhk{2}{1}{h_3}{2}\}_j$ and
      $\{\Curvejhk{2}{1}{h_3}{2}\}_i \cup \{\tCurvemjhk{2}{1}{h_3}{2}\}_j$, all
      represented on the $\CalD{2}{1}{h_3}{2}$ CS, shown in green and rendered
      as opaque. Orange and green curves correspond to Lissajous and
      quasi-Halo tori from $L_2$, respectively, while blue and pink curves
      correspond to Lissajous and quasi-Halo tori from $L_1$, respectively.
      Different marker styles correspond to different tori.  The heteroclinic
      connections between the planar Lyapunov p.o.~are shown in purple.
   }
   \label{fig:esferacas3}
\end{figure}

Figure~\ref{fig:esferacas3} presents the results for Case 3, which involves
three families of invariant tori in both center manifolds. In the left panel,
the orange curves correspond to $L_2$ Lissajous tori, whose largest
intersections start
with two contractible $\theta$-traces for $\alpha^{(2)}$ near 1, shrink with
decreasing $\alpha^{(2)}$, and disappear at tangency with the CS, as in
previous cases. The green curves correspond to $L_2$ quasi-Halo tori: for
$\alpha^{(2)}$ close to 1 they intersect as a large, single curve with
contractible $\theta$-trace (as in Figure~\ref{fig:torqhpardinsesfera1}).
Around $\alpha^{(2)}\approx1.003$, this intersection bifurcates into two
curves with non-contractible $\theta$-trace (as in
Figure~\ref{fig:tordinsesfera1}).  These curves then shrink and collapse into
two points, corresponding to transverse heteroclinic connections of the Halo
periodic orbit with the center manifold of $L_1$.

In the center panel, the blue curves represent the translated intersections
of the $L_1$ Lissajous tori, exhibiting the same behavior as before.  The pink
curves correspond to quasi-Halo intersections,
which start as before with a large curve, with the distinction that in this
case no bifurcation into two intersection curves occurs. Instead, they
contract progressively and ultimately collapse into a single point,
corresponding to the tangential intersection of the $L_1$ Halo periodic orbit
with arrival CS $\CalA21{h_2}6$.

In the right panel, we display the superposition of both figures. Both in the
$L_2$ and $L_1$ cases, the resulting pattern reveals, as before, that the
number of intersections between two tori is governed by their family type. The
shift in orientation across families is dictated by the separatrix, whose
intersection with the CS can be identified in all panels.

\subsubsection{Example: evolution of connections between invariant tori}

We end this section by giving an example of Case 1 of Table~\ref{tab:scen} on
the evolution of the number of connections between a fixed starting torus
$\Tcal^{(2)}\subset\Wcal_{h_1}^c(L_2)$, with
parameter $\alpha^{(2)}=0.945$, and a varying ending torus
$\Tcal^{(1)}\subset\Wcal_{h_1}^c(L_1)$,
$\Tcal^{(1)}\in\{\Tcal^{(1)}_i\}_{i=0}^2$, with parameters $\alpha_i^{(1)} =
0.960-0.015i$. All these tori are thus Lissajous.
This example has been chosen in order to illustrate the fact that the number of
connections between individual pairs of tori can vary within fixed starting
and ending families,
and that no symmetry can be assumed for the $\Curveju$ and $\tCurvejum$
curves.

\begin{figure}[tb]
   \centering
   \includegraphics[trim=1cm 0cm 1cm 0cm,clip]{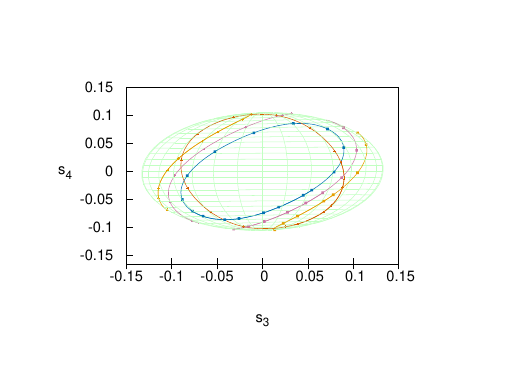}
   \includegraphics[trim=1cm 0cm 1cm 0cm,clip]{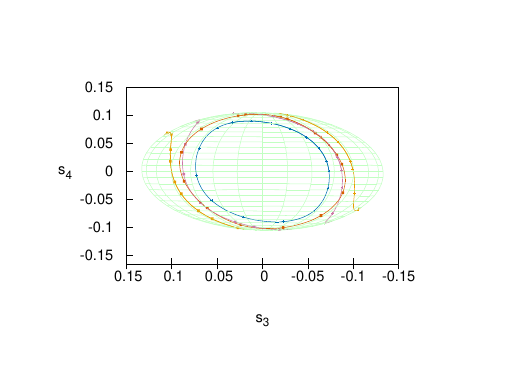}
   \caption{
      $(s_3,s_4)$ projections of $\Curvejh{2}{1}{h_1}{4}$ (red) and
      $\{\tCurvemihk{2}{1}{h_1}{4}\}_{i=0}^{2}$ (blue, pink, orange) on the
      $\CalD 21{h_1}4$ CS, rendered as opaque.
      Left: $\{s_2>0\}$ hemisphere of $\CalD 21{h_1}4$.
      Right: $\{s_2<0\}$ hemisphere. Different marker styles
      correspond to different connected components of $\Curvejh{2}{1}{h_1}{4}$
      and each $\tCurvemihk{2}{1}{h_1}{4}$.
   }
   \label{fig:tallstors}
\end{figure}

Following Appendix~\ref{app:int2tori}, the sets 
$\Curvejh{2}{1}{h_1}{4}$, $\{\tCurvemihk{2}{1}{h_1}{4}\}_i$, 
are computed and visualized in Figure~\ref{fig:tallstors}, where projections
of two hemispheres of the $\CalD{2}{1}{h_1}{4}$ CS are shown with the previous
sets on them. The intersections
$\Curvejh{2}{1}{h_1}{4}\cap\tCurvemihk{2}{1}{h_1}{4}$ provide the exact number
of heteroclinic connections between the fixed torus $\Tcal^{(2)}$ and each
torus $\Tcal^{(1)}_i$.
In Figure~\ref{fig:tallstors}, the red curves correspond to the fixed torus
$\Tcal^{(2)}$, while the blue, pink and orange curves correspond to the family
$\{\Tcal_i^{(1)}\}_i$, with parameters $\alpha_i^{(1)}$, for $i=0,1,2$,
respectively.
Different mark styles are used to represent the different connected components
of the intersection of each torus with the CS. The red curve in
Figure~\ref{fig:tallstors} left gives rise to four heteroclinic connections
with each torus $\Tcal_i^{(1)}$.  However, the red curve in
Figure~\ref{fig:tallstors} right intersects only $\Tcal_1^{(1)}$, generating
eight additional heteroclinic connections.
Table~\ref{tab:restalls} summarizes the number of connections obtained between
each pair $\Tcal^{(2)},\Tcal_i^{(1)}$.

\begin{table}
   \centering
   \begin{tabular}{ c | c | c | c }
   $i$ & $\alpha_i^{(1)}$ & Color & Number of connections \\ \hline
   $0$ & $0.960$ & Blue & $4$ \\
   $1$ & $0.945$ & Pink & $4+8$ \\
   $2$ & $0.930$ & Orange & $4$ \\
   \end{tabular}
   \caption{
   Number of heteroclinic connections between the fixed torus $\Tcal^{(2)}$ and
   each torus $\Tcal_i^{(1)}$, obtained from the intersections of 
   the curves $\Curvejh{2}{1}{h_1}{4}$ and $\{\tCurvemihk{2}{1}{h_1}{4}\}_i$.
   }
   \label{tab:restalls}
\end{table}

\begin{figure}[tbp]
   \centering
   \includegraphics{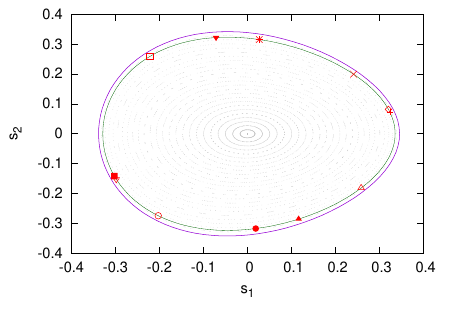}
   \includegraphics{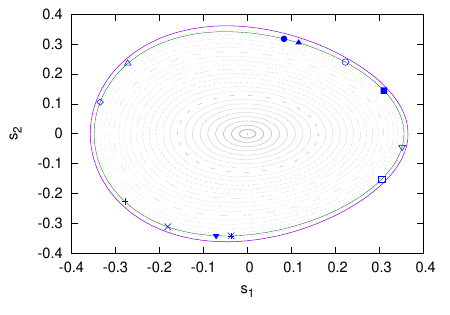}
   \includegraphics{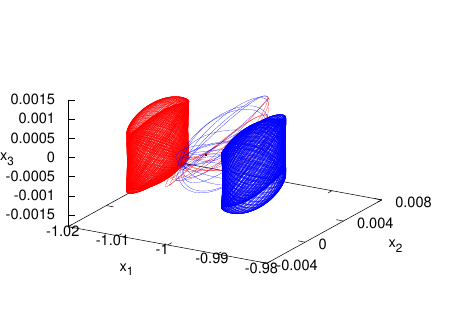}
   \includegraphics{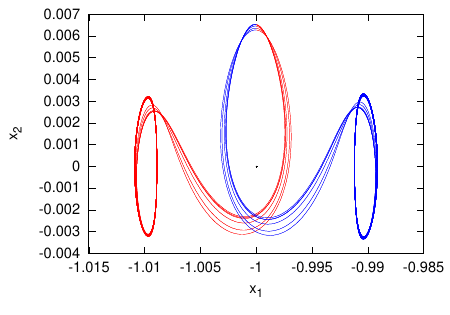}
   \caption{
      Heteroclinic connections between $\Tcal^{(2)}$ and $\Tcal^{(1)}_1$ (see
      Table~\ref{tab:restalls}). Top: departure points on
      $\Curvejh{2}{1}{h_1}{4}$ and arrival points on
      $\CalA{2}{1}{h_1}{4}(\Tcal^{(1)}_1)$ after applying the $\poincpi$
      Poincaré map; points with the same marker form a connection.  Bottom:
      corresponding trajectories in phase space.
   }
   \label{fig:torvtor}
\end{figure}

Figure~\ref{fig:torvtor} shows the $12$ heteroclinic connections
obtained for the case $i=1$. The top panel displays
the $\poincpi$ Poincaré sections of the corresponding departure and arrival
points on $\Curvejh{2}{1}{h_1}{4}$ and $\CalA{2}{1}{h_1}{4}(\Tcal^{(1)}_1)$.
The different marker types correspond to the different connections.
The bottom panel displays the resulting trajectories in phase space. 

\section{Conclusions}

In this paper we have introduced a framework to organize the dynamics of
heteroclinic and homoclinic connections associated with
the collinear libration points $L_1$ and $L_2$ of the spatial, circular RTBP.
Starting from the computation of global sets of connection surfaces and transition
maps between isoenergetic slices of the corresponding center manifolds, we have
defined and computed \emph{connectivity maps} between families of invariant tori.
These maps encode, at the level of invariant tori, which objects are connected
and provide a practical way to recover qualitative and quantitative information,
including the number of connections between pairs of tori and how these
connections vary along families.

The methodology has been illustrated on both the Sun-Earth and Earth-Moon
systems for representative energy levels covering the regimes where only
Lissajous tori are present as well as those where quasi-Halo tori coexist with
Lissajous ones. In addition to computing connectivity maps, we have used them to
construct pseudo-trajectories obtained by concatenating heteroclinic and
homoclinic transitions with trajectories inside the torus, showing how the maps
can be exploited to produce systematic transfers (e.g.~large variations of the
torus parameter and transitions between Lissajous and quasi-Halo families)
without solving full-scale boundary value problems at every step. We have also
explored how the topology of the intersections between a torus and a connection
surface is reflected in traces in torus parameter space, and how this
information can be used to estimate return times and to locate explicit
torus-to-torus connections.

Several natural extensions of the present work arise from the results
presented here.
For instance, the connectivity maps introduced here as an organizational tool
could be combined with optimization and boundary-value solvers to select, refine,
and optimize individual transfers, incorporating constraints such as
time-of-flight and/or geometric requirements in a systematic way.
Another possible extension could be to broaden this methodology
to higher-energy regimes farther from the libration points. This would
require purely numerical approaches for obtaining local approximations of the
center-stable and
center-unstable manifolds that can be obtained by novel
techniques such as \cite{2021aHaMo, 2022KuAndLla, 2026FeHaLlaMo}.
Along the same lines, an extension of these results is the application of
scattering maps, in the spirit of \cite{2026aDeGiRo}, in the new energy regimes
that can be reached.
Finally,
another promising direction would be the development of computer-assisted proofs
for the existence of the computed heteroclinic and homoclinic connection surfaces.
This would require combining validated computations of normally hyperbolic
invariant manifolds and their stable and unstable foliations with a posteriori
validation of the nonlinear equations for the connections and the involved
numerical integrations. It would also be interesting to use these techniques to
rigorously characterize the topology of the connection surfaces and to prove the
topological changes observed in the different bifurcation scenarios.

\section*{Acknowledgments}
The authors also would like to thank the Departament de Matemàtiques at UAB for
maintaining Antz, the computing cluster where most of the calculations presented
here were performed.

\appendix

\section{On the intersections of the manifold tubes of planar Lyapunov p.o.}
\label{app:primarity}

\begin{figure}[tbp]
   \centering
   \includegraphics{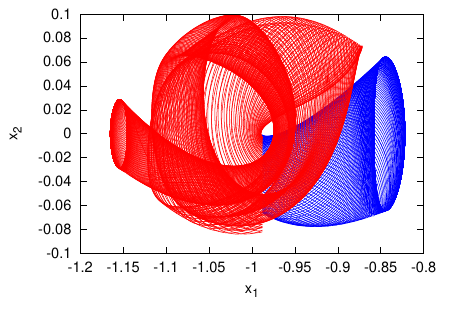}
   \includegraphics{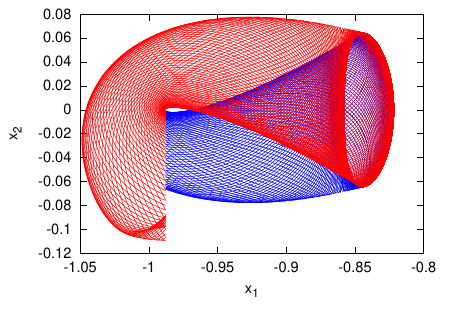}
   \caption{Intersection of the invariant manifold tubes of planar Lyapunov
   p.o.~with $\Sigma$. Left: heteroclinic case with $k=6$. Right: homoclinic
   $L_1$ case with $k=3$.}
   \label{fig:tubs}
\end{figure}

This appendix explores the relation between the geometry of the intersections
of the manifold tubes of the planar Lyapunov p.o.~and the topology of
connection surfaces. The notion of {\em primarity} of homoclinic and
heteroclinic connections plays a central role in this relation. Primary
connections correspond to the first transverse intersections of the invariant
manifolds of the connected objects, while non-primary connections arise from
subsequent intersections.
More precisely, for a fixed energy level, connections between invariant objects
only occur
after a minimal number of crossings with the Poincaré section
$\Sigma=\{x_1=\mu-1\}$. This minimal number is encoded in the index $k$.
Recall that we define $k:=k_u+k_s$, where $k_u$ (resp.~$k_s$) denotes the
number of crossings of the unstable (resp.~stable) manifold with $\Sigma$.
From $k$ we can recover the number of revolutions performed by the
corresponding connection around the second primary: it equals $(k-2)/2$ in the
heteroclinic case and $(k-1)/2$ in the homoclinic case, see
Figure~\ref{fig:tubs}. This index determines
not only whether connections exist, but also constrains the geometry of the
associated connection surfaces.

We describe below the dynamical mechanism underlying this behavior, relating the
appearance of planar connections with 
the geometry of
the connection surfaces. To visualize these surfaces, we consider sections of
the form
$\xi=\xi_0$ along the hyperbolic direction $s_5$; see Figs.~\ref{fig:somehet}
and \ref{fig:somehom}. In
these sections, the connections associated with planar
Lyapunov p.o.~appear as isolated points. 
Whenever such
planar connections exist, connections involving higher-dimensional invariant
objects inside the same isoenergetic center manifold are expected to occur
nearby.
As additional connections
associated with quasi-periodic motions are computed, 
these sets organize into continuous surfaces of
homoclinic or
heteroclinic connections. In this sense, the planar Lyapunov connections act as the
poles of the corresponding surfaces (see Figures~\ref{fig:somehet} and \ref{fig:somehom}).

Since the planar Lyapunov p.o.~and their invariant manifolds are
contained in the invariant submanifold $\{x_3=x_6=0\}$, for fixed values of
$k_u$ and $k_s$ the existence and exact number of planar connections can be
determined directly from their Poincaré sections. More precisely, one globalizes the
stable and unstable manifold tubes up to the section $\Sigma$ and considers
their projections onto the $(x_2,x_5)$ plane; see, for instance,
Figure~\ref{fig:lyapsl1movingen}, where the red and blue curves represent the
unstable and stable manifold tubes, respectively. Since every point on the
section uniquely determines a state on the corresponding energy level, each
intersection between these curves corresponds to a planar Lyapunov connection.
By repeating this procedure for different values of $k=k_u+k_s$, one can
identify the minimal value of $k$ for which intersections occur.

Once this minimal value has been identified, the corresponding homoclinic and
heteroclinic connections between planar Lyapunov p.o.~can be computed by solving
equation~\eqref{eq:conncond} while imposing $s_3^{(j)}=s_4^{(j)}=0$ for $j \in
\{1,2\}$ at both the departure and arrival periodic orbits.

We now explore several representative configurations. These examples have been
organized according to two complementary mechanisms: varying the energy while
keeping the number of crossings fixed, and varying the number of crossings while
keeping the energy fixed.

\subsection{Emergence of primary homoclinic connections for fixed $k$ and
increasing energy}

\begin{figure}[tbp]
   \centering
   \begin{subfigure}{.45\textwidth}
      \includegraphics[width=\textwidth]{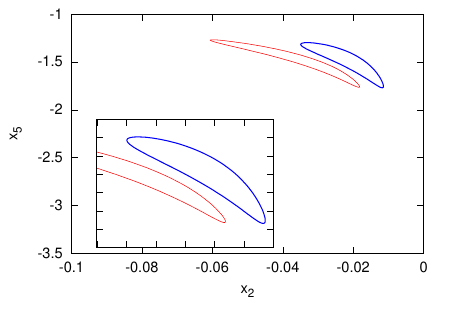}
      \caption{\label{sfig:res2esf1}$h=-1.593$}
   \end{subfigure}
   \begin{subfigure}{.45\textwidth}
      \includegraphics[width=\textwidth]{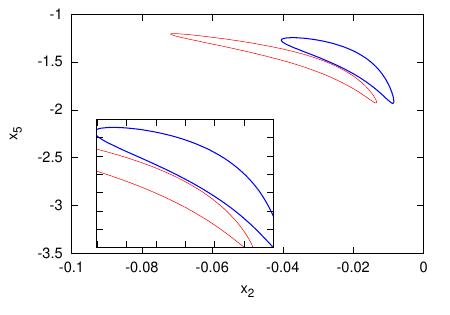}
      \caption{\label{sfig:res2esf2}$h=-1.592$}
   \end{subfigure}
   \begin{subfigure}{.45\textwidth}
      \includegraphics[width=\textwidth]{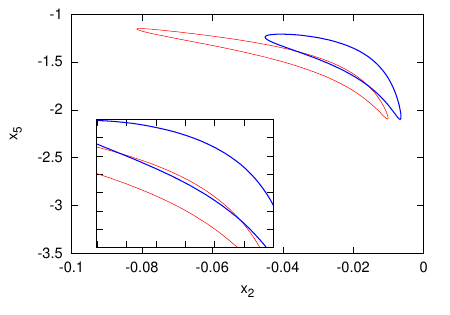}
   \caption{\label{sfig:res2esf3}$h=-1.591$}
   \end{subfigure}
   \begin{subfigure}{.45\textwidth}
      \includegraphics[width=\textwidth]{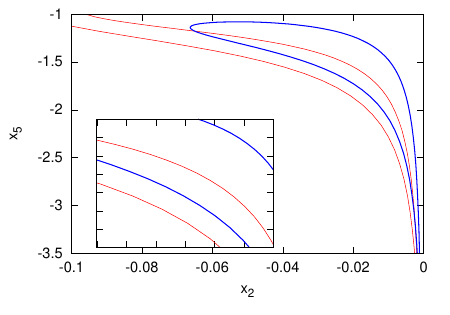}
      \caption{\label{sfig:res2esf4}$h=-1.585$}
   \end{subfigure}
   \caption{
      Poincaré sections (projection onto the $(x_2,x_5)$ plane) of the stable (blue)
   and unstable (red) manifolds of the planar Lyapunov p.o.~around $L_1$,
   for $k=3$ with $(k_u,k_s)=(2,1)$ at several energy levels.
   }
   \label{fig:lyapsl1movingen}
\end{figure}

We begin by exploring planar homoclinic connections associated with the $L_1$
planar Lyapunov p.o.~for the minimal admissible value $k=3$. Results are shown
in Figure~\ref{fig:lyapsl1movingen} where we keep $(k_u,k_s)=(2,1)$ fixed and
study how the intersections of the stable and unstable manifold tubes evolve as
the energy increases.

At the lowest energy shown, Fig.~\ref{sfig:res2esf1} ($h=-1.593$), the projected
manifolds tubes do not intersect, and therefore the Lyapunov p.o.~at this energy
level does not admit homoclinic connections with $k=3$. As the energy increases
to $h=-1.592$ (Fig.~\ref{sfig:res2esf2}), the two projected tubes become nearly
tangent. When they become tangent, a non-transverse homoclinic connection
appears and marks the threshold at which connections first become possible for
this value of $k$. 

For slightly larger energies (Figs.~\ref{sfig:res2esf3} and
\ref{sfig:res2esf4}), the tangency unfolds into two transverse intersection
points. These give rise to two transverse homoclinic connections, which are
\emph{primary} because they occur for the minimal admissible value of $k$.
These two planar primary connections act as the poles of the full
the homoclinic connection surface, which includes the connections associated
with the remaining invariant objects in the same isoenergetic center manifold.
Topologically, this surface is homeomorphic to $S^2$, as illustrated in
Figure~\ref{sfig:homoesf} for $h=-1.585$.

\subsection{Infinitely many non-primary homoclinic intersections for increasing $k$}

\begin{figure}[tbp]
   \centering
   \begin{subfigure}{.45\textwidth}
   \includegraphics[width=\textwidth]{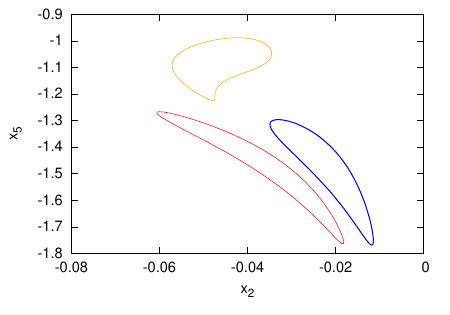}
      \caption{\label{sfig:multi1}$h=-1.593$}
   \end{subfigure}
   \begin{subfigure}{.45\textwidth}
      \includegraphics[width=\textwidth]{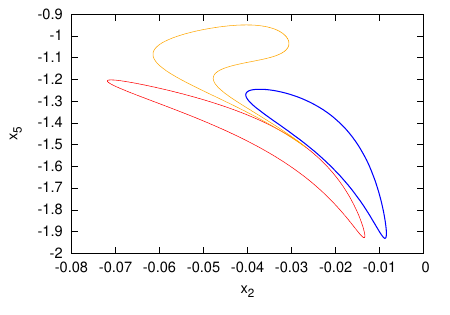}
      \caption{\label{sfig:multi2}$h=-1.592$}
   \end{subfigure}
   \begin{subfigure}{.45\textwidth}
      \includegraphics[width=\textwidth]{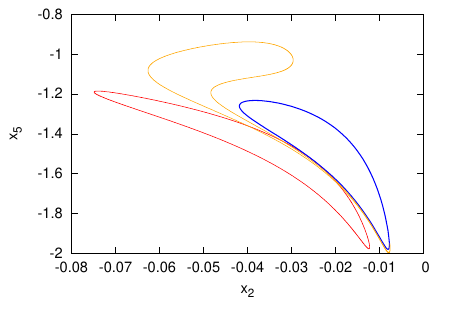}
      \caption{\label{sfig:multi3}$h=-1.5917$}
   \end{subfigure}
   \begin{subfigure}{.45\textwidth}
      \includegraphics[width=\textwidth]{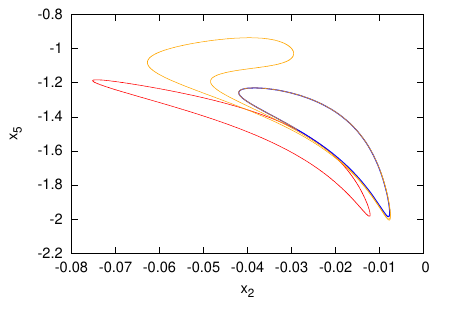}
      \caption{\label{sfig:multi4}$h=-1.5916668$}
   \end{subfigure}
   \caption{
      Poincaré sections (projection onto the $(x_2,x_5)$ plane) of the stable
      (blue/orange) and unstable (red) manifold tubes of the planar Lyapunov
      p.o.~around $L_1$, for $k=5$ with $k_u=2$ and $k_s\in\{1,3\}$. The stable
      manifold with $k_s=1$ (blue) is included as a reference.
   }
   \label{fig:lyapsl1movingeninck}
\end{figure}

We now consider the same homoclinic setting associated with the $L_1$ planar
Lyapunov p.o., but increase the total number of crossings to $k=5$. In this
case, we keep $k_u=2$ fixed and increase the number of stable crossings to
$k_s=3$ in order to examine the emergence of non-primary homoclinic
intersections.
Figure~\ref{fig:lyapsl1movingeninck} shows the corresponding projection of the
tubes onto $\Sigma$.

At the lower energy $h=-1.593$, the projected tubes (orange and red) still do
not intersect, exactly as in the case $k=3$. As the energy increases and the
projected tubes for $k=3$ (blue and red) approach tangency ($h=-1.592$), the
stable manifold after $k_s=3$ crossings (orange) begins to wind around the
reference stable tube corresponding to $k_s=1$ (blue), Figs.~\ref{sfig:multi3}
and \ref{sfig:multi4}. Then, after the first intersection of the $k=3$ tubes,
the orange tube accumulates onto its reference one, producing infinitely many
turns. 

As a consequence, transverse intersections between the unstable tube (red)
and the reference stable tube (blue) at $k=3$ produce nearby, for $k=5$, infinitely
many transverse intersections between the unstable tube and the orange curve.
Each of these intersections corresponds to a distinct \emph{non-primary}
homoclinic connection with $k=5$ crossings.  This dynamical explanation for the
progressive entangling of the invariant manifold tubes and the proliferation of
infinitely many non-primary homoclinic connections as $k$ increases is also
explored in \cite{2007GiMa,2009BaMoOll}.

From a practical numerical perspective, 
this also defines a maximal range of
values of $k_s$ and $k_u$ for which the projected tubes have not yet started to
wrap around themselves, and the number of homoclinic connections therefore
remains finite. In our computations, this threshold is determined empirically.

\subsection{Topology change of primary homoclinic connection surfaces}

\begin{figure}[tbp]
   \centering
   \begin{subfigure}{.45\textwidth}
      \includegraphics[width=\textwidth]{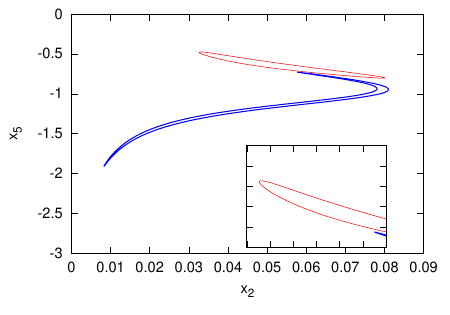}
      \caption{\label{sfig:esf2tor1}$h=-1.5853$}
   \end{subfigure}
   \begin{subfigure}{.45\textwidth}
      \includegraphics[width=\textwidth]{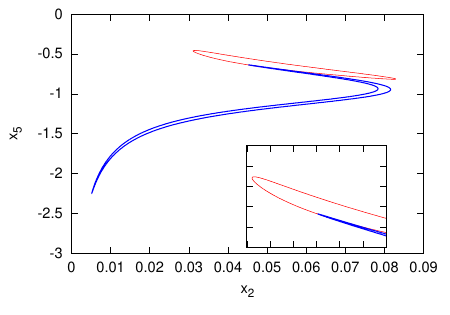}
   \caption{\label{sfig:esf2tor3}$h=-1.58515$}
   \end{subfigure}
   \begin{subfigure}{.45\textwidth}
      \includegraphics[width=\textwidth]{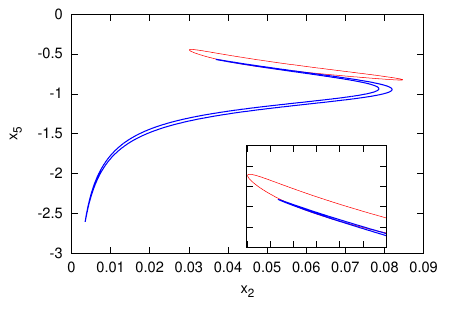}
      \caption{\label{sfig:esf2tor5}$h=-1.58505$}
   \end{subfigure}
   \begin{subfigure}{.45\textwidth}
      \includegraphics[width=\textwidth]{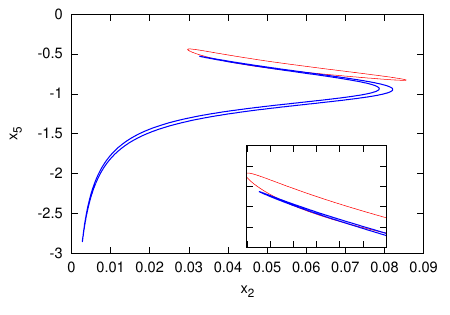}
      \caption{\label{sfig:esf2tor6}$h=-1.585$}
   \end{subfigure}
   \caption{
      Poincaré sections (projection onto the $(x_2,x_5)$ plane) of the stable
      (blue) and unstable (red) manifold tubes of the planar Lyapunov p.o.~around
      $L_2$, for $k=7$ with $(k_u,k_s)=(5,2)$ at several energy levels.
   }
   \label{fig:lyapsl2movingen}
\end{figure}

We now examine planar homoclinic connections associated with the $L_2$ planar
Lyapunov p.o.~for fixed $k=7$. The corresponding results are shown in 
Figure~\ref{fig:lyapsl2movingen}, where we take $(k_u,k_s)=(5,2)$.

At the lowest energy shown, the projected tubes do not intersect. As the energy
increases, the tubes first become tangent and then intersect
transversely at $h=-1.58515$, producing two transverse intersection points.
These intersections are primary at this energy level because they occur for the
smallest admissible value of $k$ for which connections exist. They give rise to
the primary homoclinic connection surface in $\Wcal_h^c(L_2)$, which is
homeomorphic to $S^2$. 

A second qualitative change occurs as the 
energy increases further. The projected tubes become tangent again
and this tangency subsequently unfolds into two additional transverse intersections.
As a result, from this energy level, the total number of planar intersections
increases from two to four.

This second pair of planar connections induces a topological transition of the
full homoclinic connection surface: the surface changes from being homeomorphic
to a sphere to being homeomorphic to a torus. Geometrically, two of the planar
connections lie on the exterior region of the resulting toroidal surface, while
the other two lie in its interior; see Figure~\ref{sfig:tor}.

\subsection{Non-primary homoclinic tori for increasing $k$}

\begin{figure}[tbp]
   \centering
   \begin{subfigure}{.45\textwidth}
      \includegraphics[width=\textwidth]{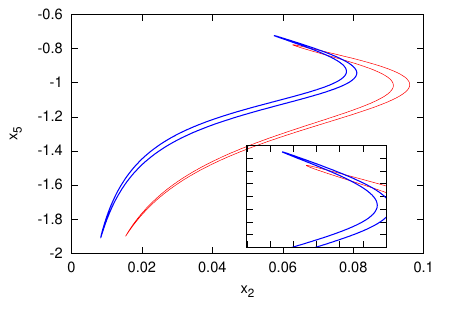}
      \caption{\label{sfig:esf2tor1_9}$h=-1.5853$}
   \end{subfigure}
   \begin{subfigure}{.45\textwidth}
      \includegraphics[width=\textwidth]{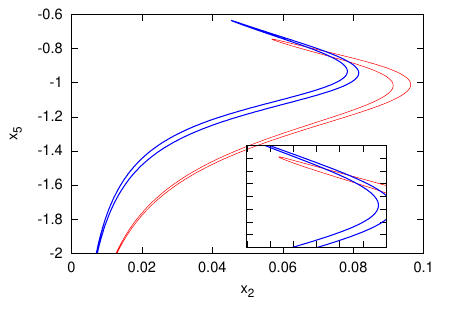}
   \caption{\label{sfig:esf2tor3_9}$h=-1.58515$}
   \end{subfigure}
   \begin{subfigure}{.45\textwidth}
      \includegraphics[width=\textwidth]{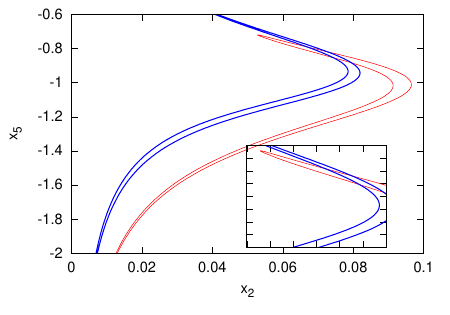}
      \caption{\label{sfig:esf2tor5_9}$h=-1.58505$}
   \end{subfigure}
   \begin{subfigure}{.45\textwidth}
      \includegraphics[width=\textwidth]{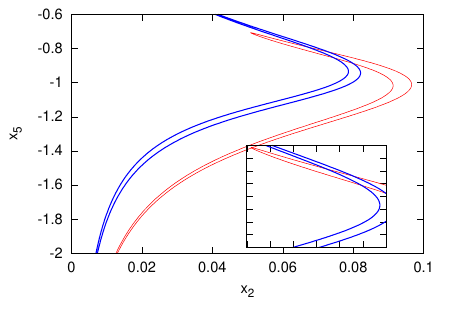}
      \caption{\label{sfig:esf2tor6_9}$h=-1.585$}
   \end{subfigure}
   \caption{
      Poincaré sections (projection onto the $(x_2,x_5)$ plane) of the unstable
      (red) and stable (blue) manifold tubes of the planar Lyapunov p.o.~around
      $L_2$ at several energy levels, for $k=9$ with $(k_u,k_s)=(5,4)$.
   }
   \label{fig:lyapsmovingen9}
\end{figure}

We next consider the same homoclinic setting associated with the $L_2$ planar
Lyapunov p.o., but increase the total number of crossings from $k=7$ to $k=9$.
The corresponding projections are shown in
Figure~\ref{fig:lyapsmovingen9}, where we keep $k_u=5$ fixed and increase the
number of stable crossings to $k_s=4$.

Because
$k_u$ remains fixed, the unstable tube has not yet started to wrap around
itself. As a consequence, the number of intersections between the projected
tubes remains finite in this regime.

Observe that only the intersections at $h=-1.5853$ are primary according to our
definition.  Indeed, for this energy there are no intersections for $k=7$
(Figure~\ref{sfig:esf2tor1}), whereas intersections already appear for $k=9$
(Figure~\ref{sfig:esf2tor1_9}).  For all larger energies shown in
Figure~\ref{fig:lyapsmovingen9}, the four transverse homoclinic intersections
are non-primary, since connections already exist for smaller values of $k$.
Numerically, these four non-primary planar intersections generate non-primary
homoclinic connection surfaces that are again homeomorphic to tori comparable to
Figure~\ref{sfig:tor}.

The next admissible odd value is $k=11$, obtained by keeping $k_u=5$ and taking
$k_s=6$. At this stage, the projected stable manifold tube begins to wrap around
the reference stable tube  corresponding to $k_s=2$. As discussed above, this
wrapping mechanism produces infinitely many non-primary homoclinic connections
between the planar Lyapunov p.o. Consequently, in this example $k_u=5$
(equivalently, $k_s=5$) is the largest value for which the number of homoclinic
connections remains finite.

\section{Construction of the continuous surfaces of connections}
\label{app:secs2sphere}

This appendix describes the numerical procedure used to transform the discrete
sets of computed connections into continuous connection surfaces and to use
these surfaces as a global Poincaré section for trajectories evolving in the
isoenergetic center manifolds.

We focus on the case in which the departure and arrival connection surfaces,
$\Dhk{l}{m}$ and $\Ahk{l}{m}$, are homeomorphic to a sphere\footnote{The
methodology is not restricted to this topology and can be applied whenever the
function defining the Poincaré section can be evaluated almost globally}.  The
central idea is to remesh the discrete connection cloud using a spherical
representation with uniformly spaced angular coordinates, from which the
continuous surface can be reconstructed by two-dimensional interpolation.  For
clarity, we present the construction only for the departure connection surface
$\Dhk{l}{m}$.  The same methodology applies to the arrival surface $\Ahk{l}{m}$.

\subsection{The equally spaced sphere of connections}
\label{app:remeshed}

To transition from the discrete set of computed connections $\hat{\mathcal
D}_{h,k}^{l,m}$ to its continuous representation $\Dhk{l}{m}$, we first remesh
the data so that the resulting surface is uniformly distributed in angular
coordinates.
We write $c=(s_1,s_2,s_3,s_4)$ for the points in $\hat{\mathcal
D}_{h,k}^{l,m}$. We then choose 
a suitable triple of coordinates
$\{s_{i_1},s_{i_2},s_{i_3}\}\subset\{s_1,s_2,s_3,s_4\}$ 
such that the connection surface is observed to be approximately spherical and
we rescale it with respect to its bounding box.
For simplicity, we assume that this
triple is $\{s_1,s_2,s_3\}$. Then, for a point $c \in \hat{\mathcal
D}_{h,k}^{l,m}\subset S_h^{(l)}$, we normalize its projected coordinates by
\begin{equation}\label{eq:bbox}
\tilde s_i = \Delta_i^{-1}(s_i-s_i^-)-\frac12,\qquad i=1,2,3,4
\end{equation}
where 
\[
   s_i^{-}:= \min_{c \in \hat{\mathcal D}_{h,k}^{l,m}} s_i, \qquad
   s_i^{+}:= \max_{c \in \hat{\mathcal D}_{h,k}^{l,m}} s_i, \qquad
   \Delta_i := s_i^+ - s_i^-.
\]
This normalization translates the mesh of the connection surface into a rescaled
box centered at the origin.

Using the rescaled coordinates,
we introduce spherical coordinates $(r,\psi,\eta)$ for the chosen triple through
\begin{equation}\label{eq:spherical}
   \begin{aligned}
      \tilde s_1 &= r \cos \psi \cos \eta, \\
      \tilde s_2 &= r \sin \psi \cos \eta, \\
      \tilde s_3 &= r \sin \eta, 
   \end{aligned}
\end{equation}
where $\psi \in [-\pi,\pi]$ denotes the longitude and $\eta \in [-\pi/2,\pi/2]$
denotes the latitude. 
In particular, a point $c$ on the departure CS $\Dlmhk$ with prescribed 
longitude $\bar \psi$ and latitude $\bar \eta$ must satisfy
\[
   \begin{aligned}
         (-\sin \bar \psi, \cos \bar \psi) \cdot (\tilde s_1, \tilde s_2) &= 0, \\
      (-\sin\bar \eta, \cos \bar \eta) \cdot (\sqrt{(\tilde s_1)^2+(\tilde s_2)^2},
      \tilde s_3) &=0.
   \end{aligned}
\]

Incorporating these angular constraints into the connection system of equations
given by \eqref{eq:conncond}, leads to the following augmented nonlinear system
\begin{equation}
   \begin{aligned}
      \Phi_{T^{(m)}}\bigl(W^{(m)}(c^{(m)},\xi^{(m)})\bigr) -
      \Phi_{T^{(l)}}\bigl(W^{(l)}(c^{(l)},\xi^{(l)})\bigr) &= 0, \\
      (-\sin \psi, \cos \psi) \cdot (\tilde s_1^{(l)}, \tilde
      s_2^{(l)}) &= 0, \\
      (-\sin \eta, \cos \eta) \cdot (\sqrt{(\tilde s_1^{(l)})^2+(\tilde
         s_2^{(l)})^2},
      \tilde s_3^{(l)}) &= 0,\\
      g\left( \Phi_{T^{(m)}}\bigl(W^{(m)}(c^{(m)},\xi^{(m)})\bigr)\right) &= 0, \\
      H \left( \Phi_{T^{(m)}}\bigl(W^{(m)}(c^{(m)},\xi^{(m)})\bigr) \right) -h
      &=0,
   \end{aligned}
   \label{eq:sysconncond}
\end{equation}
where $g(x):=x_1-\mu+1$, and the point $c^{(l)} \in \Dlmhk$ is constrained to
have fixed values of the spherical angles $\psi$ and $\eta$\footnote{An actual
implementation of \eqref{eq:sysconncond} would involve a multiple-shooting
strategy in order to deal with the hyperbolicity of the problem due to large
$T^{(l)}$ and $T^{(m)}$ integration times.}.  These two angular constraints can
be added because the original connection problem has two free degrees of freedom
along the surface $\Dlmhk$ (see \cite{2024BaHaMo}).  After fixing the longitude
and latitude, the resulting system \eqref{eq:sysconncond} becomes square.  In
all our numerical experiments, its differential was observed to have full rank.

Denoting system
\eqref{eq:sysconncond} by $G_{\psi,\eta}^{(l)}(X)$, where the vector of unknowns
is given by
\[
   X = \bigl( c^{(l)} , T^{(l)} , c^{(m)} , T^{(m)} \bigr),
\]
we solve the collection of systems
$\{G_{\psi_n,\eta_m}^{(l)}(X)\}_{n,m=0}^{N,M}$, on a prescribed grid of angles,
with $\psi_n=-\pi + \frac{2\pi n}{N}$ and $\eta_m=-\pi/2+\frac{m\pi}{M}$, in
order to obtain a remeshed set $\hDlmhk$ that is uniformly spaced in the angular
variables of the rescaled CS. For each pair
$(\psi_n,\eta_m)$, an initial
guess $\bar c$ is obtained from the original discrete set by selecting the point whose
spherical coordinates are closest to the target angles.
More precisely, we choose
\[
   \bar c(\psi_n, \eta_m)=\arg \min_{c \in \hDlmhk}
   \sqrt{(\psi_n-\psi(c))^2+(\eta_m-\eta(c))^2},
\]
where $\psi(c)$ and $\eta(c)$ are obtained by inverting \eqref{eq:spherical}.

It is worth noting that the spherical coordinates are not defined at the
poles, i.e.~$\eta=\pm\pi/2$, where the longitude $\psi$ is undefined. 
To avoid this issue, we employ two complementary coordinate charts.
Besides the projection $(s_1,s_2,s_3)$, we use the alternative ordering
$(s_1,s_3,s_2)$, which changes the north--south axis and moves the singular set
away from the region currently being meshed; in each chart, the projected
coordinates are rescaled accordingly to \eqref{eq:bbox} and
\eqref{eq:spherical}.
Then, for both projections, we solve
$G^{(l)}_{\psi,\eta}(X)$ for $\psi \in [-\pi, \pi]$ and $\eta \in [-\pi/2+\epsilon,
\pi/2-\epsilon]$ with a small $\epsilon>0$. This approach
ensures complete coverage of $\hDlmhk$ by the two charts while avoiding
the coordinate singularities at the poles.

By selecting the appropriate projection and using a two-dimensional
interpolation strategy taking as nodes the two equally spaced angular grids of
CS, we obtain a continuous map
\begin{equation}
   \begin{aligned}
   [-\pi,\pi] \times [-\pi/2, \pi/2 ] & \longrightarrow \Dlmhk \\
   (\psi,\eta) &\longmapsto c^{(l)}(\psi,\eta),
   \end{aligned}
   \label{eq:mapcont}
\end{equation}
which allows us to obtain any point on $\Dlmhk$. By interpolating the full
tuples in $\hRlmhk$, we 
analogously obtain the complete connection described by
$(c^{(l)}, T^{(l)}, c^{(m)}, T^{(m)}) \in
\Rlmhk$.
In this way, we obtain the continuous
extensions of $\Dlmhk$, $\Almhk$, and $\Rlmhk$.

\subsection{A global Poincaré section with the sphere}
\label{app:globalpoinc}

The continuous map \eqref{eq:mapcont}, together with the uniform angular meshes
constructed above, also provides a convenient way to represent the departure and
arrival connection surfaces implicitly. Such a representation is particularly
useful to compute Poincaré sections with these surfaces.

Consider the departure CS $\Dlmhk \subset \calV^{(l)}$.  Our goal is to
construct a function $g(c) : \R^4 \rightarrow \R$ such that the CS is
represented
(numerically) as $\Dlmhk:=\{c \in \R^4 : g(c) =0 \}$.

To this end, we first define a deviation function measuring the difference
between the rescaled radius of a point $c \in \R^4$, denoted by $r(c)$ (obtained
by inverting \eqref{eq:spherical}), and the radius of a point on $\Dlmhk$ having
the same angular coordinates $(\psi(c),\eta(c))$ of the point $c$. The latter is
obtained by interpolating the radial coordinate on the uniform angular mesh,
yielding the interpolant $\widehat r(\psi,\eta)$, obtained by evaluating
\eqref{eq:mapcont} at the angular coordinates $(\psi,\eta)$ and inverting
\eqref{eq:spherical} (using the appropriate
chart). Thus, we define 
\begin{equation}
\begin{aligned}
   g_r : \R^4 &\longrightarrow \R  \\
   c & \longmapsto r(c) - \widehat r\bigl(\psi(c),\eta(c)\bigr).
\end{aligned}
\label{eq:poincesf}
\end{equation}

Since spherical coordinates are computed only from the projected
three-dimensional coordinates, they do not uniquely determine a point $c \in
\R^4$.
For this reason, in addition to the radial deviation,
we also require consistency in the remaining free coordinate:
\begin{equation}
\begin{aligned}
   g_f : \R^4 &\longrightarrow \R  \\
   c & \longmapsto c_4 - \widehat c_{4}(\psi(c),\eta(c)), 
\end{aligned}
\label{eq:freevar}
\end{equation}
where $\widehat c_{4}(\psi(c),\eta(c))$ is obtained by evaluating
\eqref{eq:mapcont} at the angular coordinates of the point $c$.

Hence, the function $g(c)$ can be defined as $g(c)=g_r(c)^2+g_f(c)^2$.

\section{Computation of heteroclinic connections using frequency analysis}
\label{app:int2tori}

This appendix addresses the computation of parameterizations of invariant
tori using frequency analysis and describes the procedure
to compute heteroclinic connections
between invariant objects when one or both objects are two-dimensional invariant
tori described by these parameterizations.

\subsection{Torus reconstruction}

Consider a two-dimensional invariant torus
$\Tcal^{(j)}\subset\Vcalj$ in the isoenergetic slice
of the center manifold of $L_j$, with frequency vector
$\omega\in\R^2$.
Our goal is to reconstruct the parameterization
\[
   \Psi^{(j)} : \T^2 \longrightarrow \R^4,
\]
satisfying the invariance equation \eqref{eq:inv2tori}.

We will obtain the parameterization $\Psi$ from a real-analytic quasi-periodic
trajectory 
$f:\R\to\R^4$ of the form
\[
   f(t) = \sum_{\stackrel{\scriptstyle k \in
   \Z^2}{\langle k, \omega \rangle >0 }} a_k e^{\bmi 2\pi \langle
   k, \omega \rangle t},
\]
where $a_k\in\C^4$ are the complex amplitudes, 
$\langle k,\omega \rangle = k_1 \omega_1 + k_2 \omega_2$, and
$\omega=(\omega_1,\omega_2)$ is the frequency vector.
We assume the frequency vector satisfies a Diophantine condition and that the
Fourier coefficients satisfy the usual Cauchy estimates, ensuring analyticity of
$f$.

Following \cite{2010GoMoSi,2010GoMoSi-2}, we numerically compute the frequencies
$\{\langle k, \omega \rangle\}_{1 \le |k| \le M}$ of $f$ up to a prescribed
maximum order $M$, as well as the associated amplitudes $a_k$ (i.e., a truncated
Fourier series representation).  To reconstruct this Fourier representation
numerically, we first sample a single trajectory that densely fills the
torus.
Specifically, we choose an initial point $c\in\Tcal^{(j)}$ and integrate the
flow $\varphi^{(j)}$ over a sufficiently long time interval $T$, long enough to
guarantee both
adequate coverage of the torus and sufficient spectral resolution.
As a result, we obtain a time series given by
\[
   \left\{\bigl(t_k,
   \varphi_{t_k}^{(j)}(c)\bigr)\right\}_{k=0}^{N-1},
\]
where $t_k=kT/N$ and $\varphi_{t_k}^{(j)}(c)$
corresponds to a sample of $f(t_k)$.

We then apply frequency analysis to the sampled signal 
and retain the $N_f$
modes with largest amplitudes, which provide the dominant contribution to the
signal.
We denote the extracted integer vectors by
$k_l\in\Z^2$ and the associated frequencies by $\nu_l=\langle
k_l,\omega\rangle$.
This yields a truncated Fourier parameterization of the torus.
In angular variables $\theta=\omega t\in\T^2$, the reconstructed embedding is
given by
\[
   \Psi^{(j)}(\theta) = A_0 + \sum_{l=1}^{N_f} \left(
   A_l^c \cos \left( 2\pi \langle k_l, \theta \rangle \right)
   + A_l^s \sin \left( 2 \pi \langle k_l, \theta \rangle \right)
   \right),
\]
where $A_0\in\R^4$ is the constant term.
Up to numerical accuracy, the extracted set of frequencies is contained in the
finite set
\[
   \{\nu_1, \dots, \nu_{N_f} \} = \{\langle k , \omega \rangle : k \in \Z^2, 1
   \leq |k| \leq M, \langle k , \omega \rangle>0\},
\]
and the real coefficients are related to the complex Fourier ones by
$A_l^c=2 \Re (a_{k_l})$ and $A_l^s=-2\Im (a_{k_l})$.

The resulting truncated Fourier series provides the numerical approximation of
the torus embedding $\Psi^{(j)}$.

\subsection{Connection refinements and continuation}

Assume we have obtained parameterizations $\Psi^{(l)}$ and $\Psi^{(m)}$ of
two tori 
$\Tcal^{(l)} \subset \calV^{(l)}$ and
$\Tcal^{(m)} \subset \calV^{(m)}$
within the same energy level.
The connection
condition \eqref{eq:conncond} can be rewritten as a nonlinear system whose
solutions correspond to heteroclinic intersections between these two tori
\[
   \Wcal^u(\Tcal^{(l)}) \cap \Wcal^s(\Tcal^{(m)}) \cap \{H = h\} \cap \Sigma.
\]
Explicitly, this system can be written as
\begin{equation}
\begin{aligned}
   \Phi_{T^{(m)}}\bigl(W^{(m)}(\Psi^{(m)}(\theta^{(m)}),\xi^{(m)})\bigr) -
   \Phi_{T^{(l)}}\bigl(W^{(l)}(\Psi^{(l)}(\theta^{(l)}),\xi^{(l)})\bigr) &= 0, \\
   g\left(
   \Phi_{T^{(m)}}\bigl(W^{(m)}(\Psi^{(m)}(\theta^{(m)}),\xi^{(m)})\bigr)\right)
   &= 0, \\
\end{aligned}
\label{eq:systor2tor}
\end{equation}
with the set of unknowns
\[
   X = \left( T^{(l)}, \theta^{(l)}, T^{(m)}, \theta^{(m)}\right),
\]
where $g(x):=x_1-\mu+1$ and $\Sigma:=\{x\in\R^6: g(x)=0\}$.  System
\eqref{eq:systor2tor} is actually solved using a multiple-shooting strategy.  A
simple dimension count shows that the intersections are expected to be locally
isolated whenever they exist.  This is consistent with the system being
generically square and full rank.

Initial guesses for solving \eqref{eq:systor2tor} are obtained by discretizing
both tori, globalizing their local invariant manifolds under the flow $\Phi_t$ from
\eqref{eq:seclocws} with $c\in\{\Psi^{(j)}(\theta)\}_{\theta \in \T^2}$, for
$j=l,m$, and searching
for near intersections on the section $\Sigma$ through a minimum-distance
criterion. Solutions to this system provide the description of the sets
\eqref{eq:intcurves}.

When one torus is kept fixed while the other invariant object is allowed to vary
over the full isoenergetic center-stable or center-unstable manifold, the same
construction naturally leads to one-parameter families of intersections. System
\eqref{eq:systor2tor} can be adapted to search for intersections of the form
\begin{equation}
   \Wcal^u(\Tcal^{(l)}) \cap \Wcal^{cs}(L_m) \cap \{H = h\} \cap \Sigma,
   \label{eq:tr2wc}
\end{equation}
or 
\begin{equation}
   \Wcal^{cu}(L_l) \cap \Wcal^{s}(\Tcal^{(m)}) \cap \{H = h\} \cap \Sigma,
   \label{eq:wc2tr}
\end{equation}
which correspond to \eqref{eq:conttorWc}. 
In the first case, we solve
\begin{equation}
\begin{aligned}
   \Phi_{T^{(m)}}\bigl(W^{(m)}(c^{(m)},\xi^{(m)})\bigr) -
   \Phi_{T^{(l)}}\bigl(W^{(l)}(\Psi^{(l)}(\theta^{(l)}),\xi^{(l)})\bigr) &= 0, \\
   g\left(
      \Phi_{T^{(m)}}\bigl(W^{(m)}(c^{(m)},\xi^{(m)})\bigr)\right)
      &= 0,
\end{aligned}
\label{eq:systor2wc}
\end{equation}
for the set of unknowns
\[
   X = \left ( T^{(l)}, \theta^{(l)}, c^{(m)}, T^{(m)}\right),
\]
where $c^{(m)} \in S_h^{(m)}$ represents a point on the isoenergetic slice of
the center manifold
of $L_m$ with energy $h$. The second case is analogous.

Note that system \eqref{eq:systor2wc} has one more unknown than equations. It
has also been checked numerically to have full rank. Therefore, 
solutions can be continued locally as a one-parameter family. This is
consistent with the expected dimension of the sets \eqref{eq:tr2wc}
(resp.~\eqref{eq:wc2tr}), which
correspond to the family of heteroclinic connections between the invariant
manifold of a torus and the center-stable (resp.~center-unstable manifold) of a
libration point. 

Consequently, these systems are well suited for continuation,
allowing the computation of the full one-parameter family of heteroclinic
connections.
The initial seed for the continuation is obtained from an intersection of
$\Psi^{(l)}$ with $\Dlmhk$ (resp.~$\Psi^{(m)}$ with $\Almhk$), after
which the corresponding transition map provides the full initial state vector
$X$.

\bibliographystyle{abbrv}
\bibliography{references}

\end{document}

%% file: macros.tex
\definecolor{darkgreen}{rgb}{0,0.4,0}
\definecolor{orange}{rgb}{0.5020,0.0784,0}
\usetikzlibrary{decorations.markings}

\newcommand{\darkercolor}[3]{
    \colorlet{#3}{#1!#2!black}
}
\darkercolor{blue}{50}{DarkBlue}

\newcommand{\calW}{{\mathcal W}}
\newcommand{\calU}{{\mathcal U}}
\newcommand{\calV}{{\mathcal V}}
\newcommand{\calD}{{\mathcal D}}
\renewcommand{\Re}{{\mathop{\rm Re}\nolimits}}
\renewcommand{\Im}{{\mathop{\rm Im}\nolimits}}

\newcommand{\Vcalj}{{\mathcal V^{(j)}}}

\newcommand{\Wcal}{{\mathcal W}}

\newcommand{\Tcal}{{\mathcal T}}
\newcommand{\Dcal}{{\mathcal D}}

\newcommand{\R}{{\mathbb R}}
\newcommand{\C}{{\mathbb C}}
\newcommand{\T}{{\mathbb T}}
\newcommand{\Z}{{\mathbb Z}}
\newcommand{\N}{{\mathbb N}}
\newcommand{\bmi}{\bm{i}}

\newcommand{\Spec}{\mathop{\rm Spec}\nolimits}

\newcommand{\hRlmhk}{{\hat R_{h,k}^{l,m}}}
\newcommand{\Rlmhk}{{R_{h,k}^{l,m}}}

\newcommand{\Almhk}{{\mathcal{A}_{h,k}^{l,m}}}
\newcommand{\hDlmhk}{{\hat{\mathcal{D}}_{h,k}^{l,m}}}
\newcommand{\Dlmhk}{{\mathcal{D}_{h,k}^{l,m}}}

\newcommand{\Ahk}[2]{{\mathcal{A}_{h,k}^{#1,#2}}}
\newcommand{\Dhk}[2]{{{\mathcal{D}}_{h,k}^{#1,#2}}}
\newcommand{\tDhk}[2]{{\tilde{\mathcal{D}}_{h,k}^{#1,#2}}}
\newcommand{\tAhk}[2]{{\tilde{\mathcal{A}}_{h,k}^{#1,#2}}}
\newcommand{\CalD}[4]{{\mathcal{D}_{#3,#4}^{#1,#2}}}
\newcommand{\tCalD}[4]{{\tilde{\mathcal{D}}_{#3,#4}^{#1,#2}}}
\newcommand{\CalA}[4]{{\mathcal{A}_{#3,#4}^{#1,#2}}}

\newcommand{\Curveju}{{\Dhk{l}{m} ( \Tcal^{(l)} )}}

\newcommand{\tCurvejum}{{\tDhk{l}{m} ( \Tcal^{(m)} )}}
\newcommand{\Curveja}{{\Ahk{l}{m} ( \Tcal^{(m)} )}}

\newcommand{\Curvejuk}{{\Dhk{l}{m} ( \Tcal^{(l)}_i )}}

\newcommand{\Curvejh}[4]{{\CalD{#1}{#2}{#3}{#4} ( \Tcal^{(#1)} )}}
\newcommand{\Curvejhk}[4]{{\CalD{#1}{#2}{#3}{#4} ( \Tcal^{(#1)}_i)}}

\newcommand{\tCurvemjhk}[4]{{\tCalD{#1}{#2}{#3}{#4} ( \Tcal^{(#2)}_j)}}

\newcommand{\poincpi}{\{s_4=0, \allowbreak \,s_3>0\}}

\newcommand{\tCurvemihk}[4]{{\tCalD{#1}{#2}{#3}{#4} ( \Tcal^{(#2)}_i)}}
\newcommand{\toralfa}[3]{{\mathcal T}^{(#1)}_{#2}(#3)}
\newcommand{\round}{\mathop{\rm round}\nolimits}

%% file: rtbp.tex
\tikzset{every picture/.style={line width=0.75pt}} 

\begin{tikzpicture}[x=0.75pt,y=0.75pt,yscale=-1,xscale=1]

\draw    (378.28,224.5) -- (151,112) ;
\draw    (225.52,226) -- (151,112) ;
\draw  [dash pattern={on 0.84pt off 2.51pt}]  (301.9,86) -- (225.52,226) ;
\draw  [dash pattern={on 0.84pt off 2.51pt}]  (226.12,224.82) -- (302.5,363.32) ;
\draw  [dash pattern={on 0.84pt off 2.51pt}]  (378.28,224.5) -- (302.5,363.32) ;
\draw  [dash pattern={on 0.84pt off 2.51pt}]  (301.9,86) -- (378.28,224.5) ;
\draw    (114,226) -- (561,225) ;
\draw [shift={(563,225)}, rotate = 179.87] [color={rgb, 255:red, 0; green, 0; blue, 0 }  ][line width=0.75]    (10.93,-3.29) .. controls (6.95,-1.4) and (3.31,-0.3) .. (0,0) .. controls (3.31,0.3) and (6.95,1.4) .. (10.93,3.29)   ;
\draw    (332,373) -- (332,59) ;
\draw [shift={(332,57)}, rotate = 90] [color={rgb, 255:red, 0; green, 0; blue, 0 }  ][line width=0.75]    (10.93,-3.29) .. controls (6.95,-1.4) and (3.31,-0.3) .. (0,0) .. controls (3.31,0.3) and (6.95,1.4) .. (10.93,3.29)   ;
\draw  [fill={rgb, 255:red, 74; green, 144; blue, 226 }  ,fill opacity=1 ] (358.78,224.5) .. controls (358.35,213.73) and (366.73,205) .. (377.5,205) .. controls (388.27,205) and (397.35,213.73) .. (397.78,224.5) .. controls (398.21,235.27) and (389.83,244) .. (379.06,244) .. controls (368.29,244) and (359.21,235.27) .. (358.78,224.5) -- cycle ;
\draw  [fill={rgb, 255:red, 155; green, 155; blue, 155 }  ,fill opacity=1 ] (212.52,226) .. controls (212.23,218.82) and (217.82,213) .. (225,213) .. controls (232.18,213) and (238.23,218.82) .. (238.52,226) .. controls (238.81,233.18) and (233.22,239) .. (226.04,239) .. controls (218.86,239) and (212.81,233.18) .. (212.52,226) -- cycle ;
\draw  [fill={rgb, 255:red, 126; green, 211; blue, 33 }  ,fill opacity=1 ] (165,226) .. controls (165,223.24) and (167.24,221) .. (170,221) .. controls (172.76,221) and (175,223.24) .. (175,226) .. controls (175,228.76) and (172.76,231) .. (170,231) .. controls (167.24,231) and (165,228.76) .. (165,226) -- cycle ;
\draw  [fill={rgb, 255:red, 126; green, 211; blue, 33 }  ,fill opacity=1 ] (274,226) .. controls (274,223.24) and (276.24,221) .. (279,221) .. controls (281.76,221) and (284,223.24) .. (284,226) .. controls (284,228.76) and (281.76,231) .. (279,231) .. controls (276.24,231) and (274,228.76) .. (274,226) -- cycle ;
\draw  [fill={rgb, 255:red, 126; green, 211; blue, 33 }  ,fill opacity=1 ] (477,225) .. controls (477,222.24) and (479.24,220) .. (482,220) .. controls (484.76,220) and (487,222.24) .. (487,225) .. controls (487,227.76) and (484.76,230) .. (482,230) .. controls (479.24,230) and (477,227.76) .. (477,225) -- cycle ;
\draw  [fill={rgb, 255:red, 126; green, 211; blue, 33 }  ,fill opacity=1 ] (297.5,363.32) .. controls (297.5,360.56) and (299.74,358.32) .. (302.5,358.32) .. controls (305.26,358.32) and (307.5,360.56) .. (307.5,363.32) .. controls (307.5,366.08) and (305.26,368.32) .. (302.5,368.32) .. controls (299.74,368.32) and (297.5,366.08) .. (297.5,363.32) -- cycle ;
\draw  [fill={rgb, 255:red, 126; green, 211; blue, 33 }  ,fill opacity=1 ] (296.9,86) .. controls (296.9,83.24) and (299.14,81) .. (301.9,81) .. controls (304.66,81) and (306.9,83.24) .. (306.9,86) .. controls (306.9,88.76) and (304.66,91) .. (301.9,91) .. controls (299.14,91) and (296.9,88.76) .. (296.9,86) -- cycle ;
\draw  [fill={rgb, 255:red, 0; green, 0; blue, 0 }  ,fill opacity=1 ] (146,112) .. controls (146,109.24) and (148.24,107) .. (151,107) .. controls (153.76,107) and (156,109.24) .. (156,112) .. controls (156,114.76) and (153.76,117) .. (151,117) .. controls (148.24,117) and (146,114.76) .. (146,112) -- cycle ;

\draw (574,222.07) node [anchor=north west][inner sep=0.75pt]    {$x_1$};
\draw (327,30.4) node [anchor=north west][inner sep=0.75pt]    {$x_2$};
\draw (125,97.4) node [anchor=north west][inner sep=0.75pt]  [font=\footnotesize]  {$m_{3}$};
\draw (293,59.4) node [anchor=north west][inner sep=0.75pt]  [font=\footnotesize]  {$L_{5}$};
\draw (296,376.4) node [anchor=north west][inner sep=0.75pt]  [font=\footnotesize]  {$L_{4}$};
\draw (161,239.4) node [anchor=north west][inner sep=0.75pt]  [font=\footnotesize]  {$L_{2}$};
\draw (271,236.4) node [anchor=north west][inner sep=0.75pt]  [font=\footnotesize]  {$L_{1}$};
\draw (475,236.4) node [anchor=north west][inner sep=0.75pt]  [font=\footnotesize]  {$L_{3}$};
\draw (174,166.4) node [anchor=north west][inner sep=0.75pt]  [font=\footnotesize]  {$r_{2}$};
\draw (211,126.4) node [anchor=north west][inner sep=0.75pt]  [font=\footnotesize]  {$r_{1}$};
\draw (216,217.4) node [anchor=north west][inner sep=0.75pt]    {$M$};
\draw (372,216.4) node [anchor=north west][inner sep=0.75pt]    {$E$};

\end{tikzpicture}

%% file: paramsinvcurve.tex
\tikzset{every picture/.style={line width=0.75pt}} 

\begin{tikzpicture}[x=0.75pt,y=0.75pt,yscale=-1,xscale=1]

\draw (168,157.25) node  {\includegraphics[width=238.5pt,height=178.88pt]{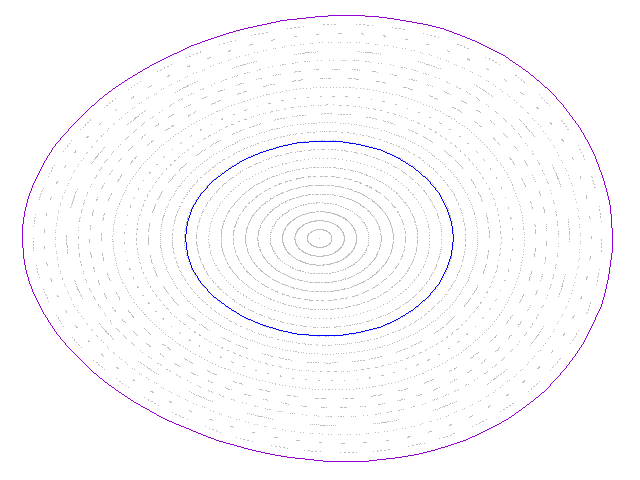}};
\draw  [dash pattern={on 0.84pt off 2.51pt}]  (168,46.5) -- (313,48.5) ;
\draw    (218,125.5) -- (167.5,156.5) ;
\draw  [dash pattern={on 0.84pt off 2.51pt}]  (167.5,156.5) -- (320,158) ;
\draw    (167.5,156.5) -- (168,46.5) ;
\draw    (164,112) -- (164,154) ;
\draw [shift={(164,157)}, rotate = 270] [fill={rgb, 255:red, 0; green, 0; blue, 0 }  ][line width=0.08]  [draw opacity=0] (7.14,-3.43) -- (0,0) -- (7.14,3.43) -- cycle    ;
\draw [shift={(164,109)}, rotate = 90] [fill={rgb, 255:red, 0; green, 0; blue, 0 }  ][line width=0.08]  [draw opacity=0] (7.14,-3.43) -- (0,0) -- (7.14,3.43) -- cycle    ;
\draw    (211.43,123.54) -- (169.07,148.96) ;
\draw [shift={(166.5,150.5)}, rotate = 329.04] [fill={rgb, 255:red, 0; green, 0; blue, 0 }  ][line width=0.08]  [draw opacity=0] (5.36,-2.57) -- (0,0) -- (5.36,2.57) -- cycle    ;
\draw [shift={(214,122)}, rotate = 149.04] [fill={rgb, 255:red, 0; green, 0; blue, 0 }  ][line width=0.08]  [draw opacity=0] (5.36,-2.57) -- (0,0) -- (5.36,2.57) -- cycle    ;
\draw  [dash pattern={on 0.84pt off 2.51pt}]  (204,134) .. controls (208,130) and (222,144) .. (210,157) ;
\draw (491.5,153.13) node  {\includegraphics[width=246.75pt,height=185.06pt]{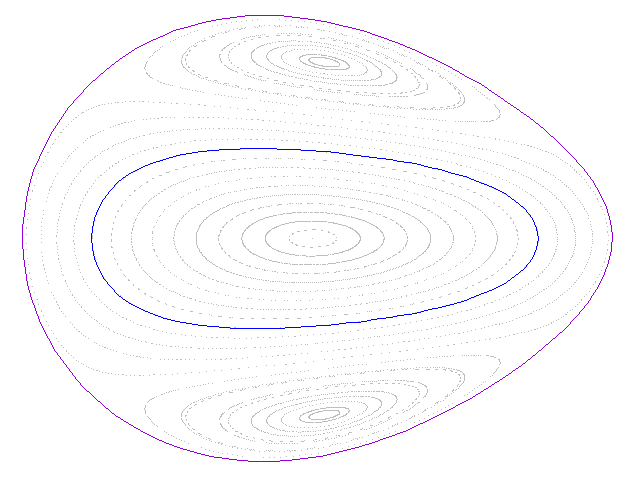}};
\draw  [dash pattern={on 0.84pt off 2.51pt}]  (488,153.5) -- (640.5,155) ;
\draw    (488,153.5) -- (494,60.5) ;
\draw  [dash pattern={on 0.84pt off 2.51pt}]  (348.5,80) -- (492.5,83) ;
\draw    (487.28,108.99) -- (484.22,150.01) ;
\draw [shift={(484,153)}, rotate = 274.26] [fill={rgb, 255:red, 0; green, 0; blue, 0 }  ][line width=0.08]  [draw opacity=0] (6.25,-3) -- (0,0) -- (6.25,3) -- cycle    ;
\draw [shift={(487.5,106)}, rotate = 94.26] [fill={rgb, 255:red, 0; green, 0; blue, 0 }  ][line width=0.08]  [draw opacity=0] (6.25,-3) -- (0,0) -- (6.25,3) -- cycle    ;
\draw    (555.24,120.68) -- (493.26,147.32) ;
\draw [shift={(490.5,148.5)}, rotate = 336.75] [fill={rgb, 255:red, 0; green, 0; blue, 0 }  ][line width=0.08]  [draw opacity=0] (5.36,-2.57) -- (0,0) -- (5.36,2.57) -- cycle    ;
\draw [shift={(558,119.5)}, rotate = 156.75] [fill={rgb, 255:red, 0; green, 0; blue, 0 }  ][line width=0.08]  [draw opacity=0] (5.36,-2.57) -- (0,0) -- (5.36,2.57) -- cycle    ;
\draw    (564,120) -- (488,153.5) ;
\draw  [dash pattern={on 0.84pt off 2.51pt}]  (544,130) .. controls (548,126) and (562,140) .. (550,153) ;
\draw  [dash pattern={on 0.84pt off 2.51pt}]  (350,57.5) -- (494,60.5) ;

\draw (281,143.9) node [anchor=north west][inner sep=0.75pt]  [font=\scriptsize]  {$\alpha =0$};
\draw (278.83,53.57) node [anchor=north west][inner sep=0.75pt]  [font=\scriptsize]  {$\alpha =1$};
\draw (165.74,124.92) node [anchor=north west][inner sep=0.75pt]  [rotate=-328.26]  {$r_{\alpha }( \delta )$};
\draw (149.95,138.33) node [anchor=north west][inner sep=0.75pt]  [font=\scriptsize,rotate=-267.96]  {$\alpha $};
\draw (196.76,139.64) node [anchor=north west][inner sep=0.75pt]  [font=\scriptsize,rotate=-359.84]  {$\delta $};
\draw (469.86,132.86) node [anchor=north west][inner sep=0.75pt]  [font=\scriptsize,rotate=-273.66]  {$\alpha $};
\draw (608,141.4) node [anchor=north west][inner sep=0.75pt]  [font=\scriptsize]  {$\alpha =0$};
\draw (339.83,69.57) node [anchor=north west][inner sep=0.75pt]  [font=\scriptsize]  {$\alpha =1$};
\draw (530.76,141.64) node [anchor=north west][inner sep=0.75pt]  [font=\scriptsize,rotate=-359.84]  {$\delta $};
\draw (491.62,124.55) node [anchor=north west][inner sep=0.75pt]  [rotate=-337.73]  {$r_{\alpha }( \delta )$};
\draw (339.33,46.57) node [anchor=north west][inner sep=0.75pt]  [font=\scriptsize]  {$\alpha =d_{h}$};

\end{tikzpicture}

%% file: sketch.tex
\begin{tikzpicture}[
    >=stealth,
    font=\small,
    x=1.8cm,   
    y=1.6cm,   
    every node/.style={inner sep=1pt}
]

\node (C0) at (0,0) {$C_0^{(2)}$};
\node (C1) at (1,0) {$C_1^{(2)}$};
\node (C2) at (2,0) {$C_2^{(1)}$};
\node (C3) at (3,0) {$C_3^{(1)}$};
\node (C4) at (4,0) {$C_4^{(2)}$};
\node (C5) at (5,0) {$C_5^{(2)}$};
\node (C6) at (6,0) {$C_6^{(1)}$};
\node (C7) at (7,0) {$C_7^{(1)}$};
\node (C8) at (8,0) {$C_8^{(2)}$};

\node (A0) at (0,1) {$\alpha_0^{(2)}$};
\node (A1) at (2,1) {$\alpha_1^{(1)}$};
\node (A2) at (4,1) {$\alpha_2^{(2)}$};
\node (A3) at (6,1) {$\alpha_3^{(1)}$};
\node (A4) at (8,1) {$\alpha_4^{(2)}$};

\draw[->,thick] (A0) -- node[above] {$f_0=4$} (A1);
\draw[->,thick] (A1) -- node[above] {$f_1=4$} (A2);
\draw[->,thick] (A2) -- node[above] {$f_2=4$} (A3);
\draw[->,thick] (A3) -- node[above] {$f_3=4$} (A4);

\draw[->,thick] (C0) -- node[above] {$x_0$} node[below] {$k_0=0$} (C1);
\draw[->,thick] (C1) -- node[above] {$x_1$} node[below] {$k_1=4$} (C2);
\draw[->,thick] (C2) -- node[above] {$x_2$} node[below] {$k_2=0$} (C3);
\draw[->,thick] (C3) -- node[above] {$x_3$} node[below] {$k_3=2$} (C4);
\draw[->,thick] (C4) -- node[above] {$x_4$} node[below] {$k_4=0$} (C5);
\draw[->,thick] (C5) -- node[above] {$x_5$} node[below] {$k_5=4$} (C6);
\draw[->,thick] (C6) -- node[above] {$x_6$} node[below] {$k_6=0$} (C7);
\draw[->,thick] (C7) -- node[above] {$x_7$} node[below] {$k_7=4$} (C8);

\end{tikzpicture}